\documentclass[11pt]{article}

\usepackage{amssymb}
\usepackage{amscd}
\usepackage{amsmath}
\usepackage{amsfonts}
\usepackage{theorem}
\usepackage{mathrsfs}
\usepackage{tikz}
\usepackage{graphicx}
\usepackage{float}
\usepackage{color}
\definecolor{red}{rgb}{.7,0,0} 
\definecolor{blue}{rgb}{0,0,1}

{\theorembodyfont{\slshape}
\newtheorem{theorem}{Theorem}
\numberwithin{theorem}{section}
\newtheorem{proposition}[theorem]{Proposition}
\newtheorem{lemma}[theorem]{Lemma}
\newtheorem{corollary}[theorem]{Corollary}

\newtheorem{hypothesis}[theorem]{Hypothesis}
}
{\theorembodyfont{\rmfamily}
\newtheorem{definition}[theorem]{Definition}
\newtheorem{problem}[theorem]{Problem}
\newtheorem{remark}[theorem]{Remark}

}
\def\proof{{\noindent\sc Proof. \quad}}
\newcommand{\proofof}[1]{{\noindent\sc Proof of #1. \quad}}
\def\eproof{{\mbox{}\hfill\qed}\medskip}
\newcommand\qed{{\unskip\nobreak\hfil\penalty50\hskip2em\vadjust{}
\nobreak\hfil$\Box$\parfillskip=0pt\finalhyphendemerits=0\par}}

\def\R{{\mathbb{R}}}
\def\Z{{\mathbb{Z}}}

\def\N{{\mathbb{N}}}
\def\C{{\mathbb{C}}}

\def\CU{{\mathcal{U}}}
\renewcommand{\P}{\mathbb{P}}
\def\E{\mathop{\mathbb{E}}}
\newcommand{\Exp}{\mathop{\mathbb{E}}}
\def\Prob{\mathop{\mathsf{Prob}}}

\def\mcU{\mathcal{U}}

\def\mcD{\mathcal{D}}

\def\Cnn{{\C^{n\times n}}}
\def\mcS{{\mathcal{S}}}
\def\Oh{{\mathcal{O}}}
\def\msD{\mathscr D}

\def\msU{\mathscr U}
\def\dist{\mathsf{dist}}
\def\diag{\mathsf{diag}}
\def\NJ{\mathrm{NJ}}
\def\rank{\mathsf{rank}}
\def\cond{\mathsf{cond}}

\def\Id{\mathsf{Id}}

\def\uno{\mbox{1\hspace*{-2.5pt}l}}

\renewcommand{\tilde}{\widetilde}
\def\a{\alpha}
\def\e{\varepsilon}

\def\s{\sigma}
\def\dS{d_{\SS}}
\def\dpr{d_{\P}}

\def\EC{{\sf Path-follow}}

\def\hA{\hat{A}}
\def\hB{\hat{B}}

\def\Cnn{\C^{n\times n}}
\def\Pn{\P(\C^n)}

\def\SS{\mathbb{S}}
\def\av{\mathsf{av}}
\def\mum{\mu_{\max}}
\def\muFa{\mu_{F,\av}}

\newcommand{\algoritmo}{\begin{minipage}{0.87\hsize}\linea}
\newcommand{\falgoritmo}{\linea\end{minipage}\bigskip}
\newcommand{\linea}{\vspace*{-5pt}\hrule\vspace*{5pt}}
\newtheorem{algorithm}{Algorithm}
\def\espacio{\hspace*{1cm}}

\def\eeespacio{\hspace*{2cm}}
\newcommand{\inputalg}[1]{\linea\bf Input:\quad\rm #1\vspace*{3pt}}
\newcommand{\specalg}[1]{\bf Preconditions:\quad\rm #1}
\newcommand{\Output}[1]{\linea\bf Output:\quad\rm #1\vspace*{2pt}}
\newcommand{\postcond}[1]{\bf Postconditions:\quad\rm #1\vspace*{3pt}}
\newcommand{\bodyalg}[1]{\linea\tt #1\vspace*{3pt}}

\def\la{\lambda}
\def\avcost{{\mathsf{Avg\_Cost}}}
\def\smcost{{\mathsf{Smd\_Cost}}}
\def\aviter{{\mathsf{Avg\_Num\_Iter}}}

\def\siter{{\mathsf{Smd\_Num\_Iter}}}
\def\siterall{{\mathsf{Smd\_Num\_Iter\_All}}}
\def\aviterall{{\mathsf{Avg\_Num\_Iter\_All}}}
\medskip

\def\mcN{\mathcal{N}}

\def\Vol{\mathsf{vol}}
\newcommand{\vol}{\mathsf{vol}}

\def\pil{\pi^{\mathsf{lin}}}
\def\phil{\phi^{\mathsf{lin}}}
\newcommand{\mnc}{\C^{n\times n}}

\newcommand{\pes}[2]{\langle #1,#2\rangle}%{\left\langle #1,#2\right\rangle}

\newcommand{\Pnn}{\P(\C^n)}

\newcommand{\Un}{\mathcal{U}_n}%{\mathcal{U}(n+1)}

\newcommand{\ambient}{\Cnn\times\C\times\prc}
\newcommand{\ambientu}{\Cnn\times\C}
\newcommand{\prc}{\P(\mathbb{C}^n)}

\newcommand{\V}{\mathcal{V}}

\newcommand{\wW}{\widetilde{\W}}
\newcommand{\A}{\mathcal{A}}
\newcommand{\Vlin}{\mathcal{V}^{\mathsf{lin}}}
\newcommand{\W}{\mathcal{W}}

\newcommand{\then}{\Rightarrow}

\begin{document}

\bibliographystyle{plain}

\makeatletter

%%%%%%%
% Some abreviations for the bib database
%%%%%%%

\def\JACM{Journal of the ACM}
\def\CACM{Communications of the ACM}
\def\ICALP{International Colloquium on Automata, Languages
            and Programming}
\def\STOC{annual ACM Symp. on the Theory
          of Computing}
\def\FOCS{annual IEEE Symp. on Foundations of Computer Science}
\def\SIAM{SIAM Journal on Computing}
\def\SIOPT{SIAM Journal on Optimization}
\def\MOR{Math. Oper. Res.}
\def\BSMF{Bulletin de la Soci\'et\'e Ma\-th\'e\-ma\-tique de France}
\def\CRAS{C. R. Acad. Sci. Paris}
\def\IPL{Information Processing Letters}
\def\TCS{Theoretical Computer Science}
\def\BAMS{Bulletin of the Amer. Math. Soc.}
\def\TAMS{Transactions of the Amer. Math. Soc.}
\def\PAMS{Proceedings of the Amer. Math. Soc.}
\def\JAMS{Journal of the Amer. Math. Soc.}
\def\LNM{Lect. Notes in Math.}
\def\LNCS{Lect. Notes in Comp. Sci.}
\def\JSL{Journal for Symbolic Logic}
\def\JSC{Journal of Symbolic Computation}
\def\JCSS{J. Comput. System Sci.}
\def\JoC{J. of Complexity}
\def\MP{Math. Program.}
\sloppy

\begin{title}
{{\bf  A stable, polynomial-time algorithm for the eigenpair problem}} 
\end{title}
\author{
Diego Armentano\thanks{Partially supported by Agencia Nacional de Investigaci\'on e Innovaci\'on (ANII), Uruguay, and by CSIC group 618}\\
Universidad de La Rep\'ublica\\
URUGUAY\\
e-mail: {\tt diego@cmat.edu.uy}
\and
Carlos Beltr\'an\thanks{ partially suported by
the research projects MTM2010-16051 and MTM2014-57590 from Spanish Ministry of Science MICINN}\\
Universidad de Cantabria\\ 
SPAIN\\
e-mail: {\tt beltranc@unican.es}
\and
Peter B\"urgisser\thanks{Partially funded by DFG research grant BU 1371/2-2}\\
Technische Universit\"at Berlin\\
%Institut f\"ur Mathematik\\
%Sekretariat MA 3-2\\
%Stra{\sz}e des 17. Juni 136
%10623 Berlin, 
GERMANY\\
e-mail: {\tt pbuerg@math.tu-berlin.de}
\and
Felipe Cucker\thanks{Partially funded by
a GRF grant from the Research Grants Council of the
Hong Kong SAR (project number CityU 100813).}\\
%Department of Mathematics\\
City University of Hong Kong\\
HONG KONG\\
e-mail: {\tt macucker@cityu.edu.hk}
\and
Michael Shub\\
City University of New York\\
U.S.A.\\
e-mail: {\tt shub.michael@gmail.com}
}

\date{}

\makeatletter
\maketitle
\makeatother

\thispagestyle{empty}

\begin{quote}
{\small 
{\bf Abstract.} 
We describe algorithms for computing eigenpairs 
(eigenvalue-eigenvector pairs) of a complex $n\times n$ 
matrix $A$. These algorithms are numerically stable, 
strongly accurate, and theoretically efficient (i.e., polynomial-time). 
We do not believe they outperform in practice the algorithms 
currently used for this computational problem. The merit 
of our paper is to give a positive answer to a 
long-standing open problem in numerical linear algebra. 
}
\end{quote}
\medskip

\hfill\begin{minipage}{6cm}
{\small{\em  
So the problem of devising an algorithm [for the eigenvalue 
problem] that is numerically stable and globally (and quickly!) 
convergent remains open.}

\hfill J.~Demmel~\cite[page~139]{Demmel97}}
\end{minipage}
\bigskip\bigskip

\section{Introduction}

\subsection{The problem}

The quotation from Demmel opening this article, though 
possibly puzzling for those who day-to-day satisfactorily 
solve eigenvalue problems, summarizes a long-standing 
open problem in numerical linear algebra. The first algorithm 
that comes to mind for computing eigenvalues 
---to compute the characteristic polynomial 
$\chi_A$ of $A$ and then compute (i.e., approximate) its zeros---
has proved to be numerically unstable in practice. The so called 
Wilkinson's polynomial, 
$$
  w(x):=\prod_{i=1}^{20}(x-i) 
  = x^{20}+w_{19}x^{19}+\cdots+w_1x+w_0
$$
is often used to illustrate this fact. For a diagonal matrix $D$ 
with diagonal entries $1,2,\ldots,20$ (and therefore with 
$\chi_D(x)=w(x)$) an error of $2^{-23}$ in 
the computation of $w_{19}=-210$ produces, even if the rest of the 
computation is done error-free, catastrophic variations in the zeros 
of $\chi_D$. For instance, the zeros at 18 and 19 collide into a double 
zero close to 18.62, which will unfold into two complex conjugate 
zeros if the error in $w_{19}$ is further increased. And yet, there 
is nothing wrong in the nature of $D$ (in numerical analysis terms, 
as we will see below, 
$D$ is a well-conditioned matrix for the eigenvalue problem). 
The trouble appears to lie in the method. 

Barred from using this immediate algorithm due to its numerical 
instability, researchers devoted efforts to come up with alternate 
methods which would appear to be stable. Among those proposed, 
the one that is today's algorithm of choice is the iterated QR with 
shifts. This procedure behaves quite efficiently in general and yet, 
as Demmel pointed out in 1997~\cite[p.~139]{Demmel97}, 
\begin{quote}
{\small 
after more than 30 years of dependable service,
convergence failures of this algorithm have quite recently been
observed, analyzed, and patched [\dots]. But there is still no global
convergence proof, even though the current algorithm is considered
quite reliable. 
}
\end{quote}
Our initial quotation followed these words in Demmel's text. It 
asked 
for an algorithm which will be numerically stable and for which, 
convergence, and if possible small complexity bounds, can be 
established. Today, 17 years after Demmel's text, this demand retains 
all of its urgency: it is not known if any of the standard numerical linear algebra 
algorithms satisfies the properties above. For example:
\begin{itemize}
\item 
The unshifted QR algorithm terminates with probability 1
but is probably infinite average cost if approximations to the eigenvectors 
are to be output (see~\cite{Kostlan1988}).
\item 
The QR algorithm with Rayleigh Quotient shift fails for open sets of real 
input matrices (see~\cite{BattersonSmillie1989,BattersonSmillie1990}).
\item 
We do not know whether the Francis (double) shift algorithm converges 
generally 
on real or complex matrices, nor an estimate of its average cost.
\item Other algorithms in modern texts are analyzed but don't estimate the (necessarily infinite in the worst case) number of iterations, 
which usually relies on experimental results, see for 
example~\cite{NakatsukasaHigham2013} which uses a divide and 
conquer algorithm or~\cite{Demmel2007} which is on turn based on the 
algorithm in~\cite{Baietal}.
\end{itemize}

Algorithms which output approximate eigenvalues without accompanying 
approximate eigenvectors might be easier to analyze. The experimental 
evidence of~\cite{PfrangDeiftMenon} for symmetric matrices suggests that
many of the algorithms in use are of average finite cost and even that there 
is some universality. An informal explanation of this fact is that the 
eigenvalues of symmetric matrices are very well conditioned, see for 
example~\cite[eq. (1.5)]{vanLoan}. But eigenvectors are another matter. 
When the matrices 
are close to having multiple eigenvalues the condition of the eigenvector 
tends to infinity. For example, even for $2\times 2$ symmetric matrices, 
any pair of orthogonal vectors $(a,b)$ and $(-b,a)$ are the eigenvectors 
of a matrix
\[
\begin{pmatrix}
1+\varepsilon_1 &\varepsilon_3\\\varepsilon_3&1+\varepsilon_2
\end{pmatrix}
\]
for $|\varepsilon_i|$, $i=1,2,3$, arbitrarily small. 

The only goal of this paper is to give a positive answer to Demmel's 
question.  The set of our main results can be informally stated as 
follows. 
\medskip

\noindent
{\bf main results\quad}
{\sl We exhibit algorithms which on input a complex matrix $A$ with 
complex Gaussian entries generate (with probability $1$) 
an ``excellent'' approximation to one of (or all) the 
(eigenvalue, eigenvector) pairs of $A$. Some of these algorithms 
are deterministic while some are randomized. Their running time (expected running time for the randomized case) 
is polynomial in $n$ on average (w.r.t.~$A$) as well as under \
a standard smoothed analysis.}
\medskip

More precisely, the average complexity bounds we prove, for
$n\times n$ matrices, are $\Oh(n^7)$ ---for the 
computation of a single eigenpair with either a deterministic 
or a randomized algorithm--- and $\Oh(n^9)$ 
---for the computation  of all eigenpairs with a 
deterministic algorithm. Note that these are just upper bounds: the
practical performance of the algorithms might be better.  The
precise statements of the main results 
are in Theorems~\ref{thm:main},~\ref{thm:main2}, 
and~\ref{th:random}.%\ref{thm:str_rand}.

\subsection{A few words on approximations}

It must be said upfront that we do not think the algorithm we 
propose will outperform, in general, iterated QR with shifts.  
It nonetheless possesses some worthy features which we 
want to describe in this introduction. The key one, we already 
mentioned, is that both convergence and complexity bounds 
can be established for it. It is also numerically stable. But in 
addition, it is {\em strongly accurate}. 

A starting point to understand the meaning of this last claim, 
is the observation that there are two different obstructions to 
the exact computation of an eigenpair. Firstly, the use of finite 
precision, and the ensuing errors accumulating during the 
computational process. The expression {\em numerically stable} is
usually vested on algorithms for which this accumulated error on the 
computed quantities is not much 
larger than that resulting from an error-free computation on an 
input datum which has been read (and approximated) with machine 
precision. Secondly, the nonlinear character of the equations 
defining eigenvalues and eigenvectors in terms of the given matrix. 
For $n\geq 5$, 
we learned from Abel and Galois, we cannot write down these eigenvalues 
in terms of the matrix' entries, not even using radicals, and the same remains true for eigenvectors. Hence, 
we can only compute approximations of them and this is so 
{\em even assuming infinite precision in the computation}. 

The expression {\em strongly accurate} refers to the quality of these 
approximations. It is common to find in the literature (at least) 
three notions of approximation which we next briefly 
describe. To simplify, we illustrate with the computation of 
a value $\zeta\in\C$ from a datum $A\in\C^N$ (and the reader 
may well suppose that this computation is done with infinite 
precision). We let $\tilde{\zeta}$ be the quantity actually computed
and we consider the following three requirements on it:
\begin{description}
\item{\em Backward approximation.}
The element $\tilde{\zeta}$ is the solution of a datum $\tilde{A}$ 
close to~$A$. Given $\e>0$, we say that $\tilde\zeta$ is an 
$\e$-backward approximation when $\|A-\tilde{A}\|\leq \e$ (resp. $\|A-\tilde{A}\|\leq \e\|A\|$ if we are interested in relative errors). 
\item{\em Forward approximation.} 
The quantity $\tilde{\zeta}$ is close to $\zeta$. 
Given $\e>0$, we say that $\tilde\zeta$ is an 
$\e$-forward approximation when $|\zeta-\tilde{\zeta}|\leq \e$ (resp. $|\zeta-\tilde{\zeta}|\leq \e|\zeta|$).
\item{\em Approximation \`a la Smale.} An appropriate version of Newton's 
iteration, starting at $\tilde{\zeta}$, converges immediately, quadratically 
fast, to $\zeta$, either in absolute or relative error. 
\end{description}
These requirements are increasingly demanding. For instance, 
if $\zeta$ is an $\e$-backward approximation then the forward error 
$|\zeta-\tilde{\zeta}|$ is bounded, roughly, by $\e\,\cond(A)$.  
Here $\cond(A)$ is the condition number of $A$, 
a quantity usually greater than 1. So, in general, 
$\e$-backward approximations are not $\e$-forward approximations,  
and if $A$ is poorly conditioned $\tilde{\zeta}$ may be a 
poor forward approximation of $\zeta$. 
We also observe that if $\tilde{\zeta}$ is an approximation \`a la Smale 
we can obtain an $\e$-forward approximation from $\tilde{\zeta}$ 
by performing $\Oh(\log|\log\e|)$ Newton's steps. Obtaining an 
approximation \`a la Smale from an $\e$-forward approximation 
is a much less obvious process. 
\smallskip

When we say that our algorithm is strongly accurate, we refer to 
the fact that the returned eigenpairs are approximations \`a la Smale 
of true eigenpairs. 

In our case, we can not only efficiently compute $\e$-forward 
approximations as above but also with respect to relative error. 
These are pairs $(\zeta,w)$ satisfying 
$$
   d_{\P}(w,v)\leq\e
   \qquad\mbox{and}\qquad
   |\zeta-\lambda|\leq\e|\lambda|
$$
for some true eigenpair $(\lambda,v)$ of $A$. Here 
$d_{\P}(w,v)$ denotes the angle between $w$ and~$v$ 
(note that here the scaling of eigenvectors renders moot 
the relativization of the error).

\begin{theorem}\label{th:relativeerror}
We exhibit an algorithm that, given a matrix $A\in\Cnn$,
an approximate eigenpair $(\zeta,w)$ returned by any of 
the algorithms in the main results, and an accuracy 
$\e\in(0,1/2)$, produces an $\e$-forward approximation 
(in relative error) of the approximated true eigenpair $(\lambda,v)$. 
The algorithm terminates provided $\lambda\neq0$. 
Its average cost over Gaussian 
matrices $A$ (independently of the chosen approximate eigenpair) 
is $\Oh(n^3\log_2\log_2(n/\e))$.
\end{theorem}

Combining our main results with Theorem \ref{th:relativeerror} we can thus compute $\e$-forward approximations (in relative error) of all eigenpairs of random Gaussian matrices with average running time $\Oh(n^9+n^3\log_2\log_2(n/\e))$. See \S\ref{sec:mainrelativeepsilon} for a proof of 
Theorem~\ref{th:relativeerror}.

\subsection{A few words on complexity}

The cost, understood as the number of arithmetic operations 
performed, of computing an approximation 
of an eigenpair for a matrix $A\in\Cnn$, depends on the 
matrix $A$ itself. Actually, and this is a common feature in 
numerical analysis, it depends on the condition $\cond(A)$ of 
the matrix $A$. But this condition is not known a priori. It was 
therefore advocated by Smale~\cite{Smale97} to eliminate this 
dependency in complexity bounds by endowing data space 
with a probability distribution and estimating 
expected costs. This idea has its roots in early work of 
Goldstine and von Neumann~\cite{vNGo51}. 

In our case, data space is $\Cnn$, and a common 
probability measure to endow it with is the standard 
Gaussian. Expectations of cost w.r.t.~this measure 
yield expressions in $n$ usually referred to as 
{\em average cost}. A number of considerations, including 
the suspicion that the use of the standard Gaussian could 
result in complexity bounds which are too optimistic compared 
with ``real life'', prompted Spielman and Teng to introduce 
a different form of probabilistic analysis, called {\em smoothed 
analysis}. In this, one replaces the average analysis' goal of 
showing that 
\begin{quote}
for a random $A$ it is unlikely that the cost for $A$ will be large
\end{quote}
by the following
\begin{quote}
for all $\hA$, it is unlikely that a slight random 
perturbation $A=\hA+\Delta A$ will require a large cost. 
\end{quote}
The expectations obtained for a smoothed analysis will now be 
functions of both the dimension $n$ and some measure of 
dispersion for the random perturbations (e.g., a variance). 

Smoothed analysis was first used for the simplex method of 
linear programming~\cite{ST:04}. Some survey expositions 
of its rationale are in~\cite{ST:02,ST:09,Bu:10}. One may argue that it 
has been well accepted by the scientific community from 
the fact that Spielman and Teng were awarded the 
G\"odel 2008 and Fulkerson 2009 prizes for it (the former 
by the theoretical computer science community and the latter 
by the optimization community). Also, in 2010, Spielman was awarded 
the Nevanlinna prize, and smoothed analysis appears in the 
laudatio of his work.  
\smallskip

In this paper we will exhibit bounds for the cost of our algorithm 
both for average and smoothed analyses. 

\subsection{A few words on numerical stability}

The algorithm we deal with in this paper belongs to the class 
of homotopy continuation methods. Experience has shown 
that algorithms in this class are very stable and stability 
analyses have been done for some of 
them, e.g.~\cite{BriCuPeRo,BeltranLeykin2011,DedieuMalajovichShub}. 
Because of this, we will assume 
infinite precision all along this paper and steer clear of 
any form of stability analysis. We nonetheless observe that 
such an analysis can be easily carried out following the steps 
in the papers mentioned above.

\subsection{Previous and related work}

Homotopy continuation methods go back, at least, to the 
work of Lahaye~\cite{Lahaye}. A detailed survey of their 
use to solve polynomial equations is in~\cite{Li:03}. 
More explicit focus in eigenvalue computations is considered 
in~\cite{Chu1984,LiSauer1987,LiZeng1992,LiZengCong1992} 
but we do not know of any serious attempt to implement them.

In the early 1990s Shub and Smale set up a program to 
understand the cost of solving square systems of complex 
polynomial equations using homotopy methods. In a 
collection of articles~\cite{Bez1,Bez2,Bez3,Bez4,Bez5}, known 
as the {\em B\'ezout series}, they put in place many of the notions 
and techniques that occur in this article. The B\'ezout series 
did not, however, conclusively settle the understanding of the cost 
above, and in 1998 Smale proposed it as the 17th in his list 
of problems for the mathematicians of the 21st 
century~\cite{Smale98}. The problem is not yet considered fully 
solved by the community but 
significant advances appear 
in~\cite{BePa:09,BePa:11,BuCu11}. 

The results in these papers cannot be directly used for 
the eigenpair problem since instances of the latter are ill-posed as  
polynomial systems. But the intervening ideas can be 
reshaped to attempt a tailor-made analysis for the eigenpair 
problem. A major step in this direction was done by 
Armentano in his PhD thesis (see~\cite{Armentano:13} and its precedent~\cite{DeSh00}), 
where the condition number $\mu$ for the eigenpair problem 
was exhaustively studied. A further step was 
taken in~\cite{ArmCuc} where $\mu$ was used to analyze a 
randomized algorithm for the Hermitian eigenpair problem. 
%A difference between our paper and 
%both~\cite{Armentano:13} and~\cite{ArmCuc} is that in the latter 
%the technical development binds inputs and outputs (eigenvalues) 
%together. We have found more natural to uncouple them. 

Our paper follows this stream of research.

\subsection{Structure of the exposition}

The remainder of this paper is divided into two parts. 
In the first one, \S\ref{ss:GeoFrame} below, 
we introduce all the technical preliminaries, we describe 
with details the algorithms, and we state our main 
results (Theorems~\ref{thm:main},~\ref{thm:main2}, 
and~\ref{th:random}). 
The condition number~$\mu$, Newton's method,  
the notion of approximate eigenpair, and 
Gaussian distributions are among these 
technical preliminaries. The second part, which 
occupies us in the subsequent sections, 
is devoted to proofs. Some short proofs are included in 
\S\ref{ss:GeoFrame} as well.

\section{Preliminaries, Basic Ideas, and Main Result}\label{ss:GeoFrame}

\subsection{Spaces and Metric Structures}\label{subsec:spaces}

Let $\Cnn$ be the set of $n\times n$ complex matrices. 
We endow this complex linear space with the restriction of the 
real part of the {\em Frobenius Hermitian product} $\pes~~_F$ and 
the associated  {\em Frobenius norm} $\|\cdot\|_F$ given by  
$$
  \pes{A}{B}_F:=\mbox{trace }(B^*A)=\sum_{i,j=1}^n a_{ij}\,
  \overline {b_{ij}},\quad \|A\|_F=\pes{A}{A}^{1/2},
$$
where $A=(a_{ij})$ and $B=(b_{ij})$. The unit sphere will be denoted 
by $\SS(\Cnn)$ or simply by $\SS$. We endow the product vector space $\ambientu$ with the 
canonical Hermitian inner product structure and its 
associated norm structure and (Euclidean) distance. 

The space  $\C^n$ is equipped with the canonical Hermitian inner 
product $\pes~~$. We denote by $\P(\C^n)$ its associated 
projective space. This is a smooth manifold which carries a natural
Riemannian metric, namely, the real part of the \emph{Fubini-Study
metric} on $\Pnn$. The Fubini-Study metric is the Hermitian structure
on $\P(\C^n)$ given in the following way: for $x\in\C^n$, 
\begin{equation}\label{eq:defTip}
   \pes{w}{w'}_x:=\frac{\pes{w}{w'}}{\|x\|^2},
\end{equation}
for all $w,\,w'$ in the Hermitian complement $x^\perp$ of $x$ in
$\C^n$. We denote by $d_\P$ the associated Riemannian 
distance. An explicit formula for that distance 
(see for example~\cite[p. 226]{BlCuShSm98}) is
\begin{equation}\label{eq:distanceP}
d_\P(v,w)=\arccos\frac{|\langle v,w\rangle|}{\|v\|\cdot\|w\|}.
\end{equation}
Note that this formula makes sense for $v,w\in\C^n$ (as the distance between the respective projective classes). 

The space $\ambient$ is endowed with the Riemannian product structure. 
The resulting distance equals 
$$
   \dist((A,\la,v),(A',\la',v'))^2 := 
   \|A-A'\|_F^2+|\lambda-\lambda'|^2+\dpr(v,v')^2.
$$
We will only use this distance  on $\SS\times\C \times \P(\C^n)$.
Note that for $A,A' \in\SS$, the distance $\dist((A,\la,v),(A',\la',v'))$ is smaller than or equal to 
the natural geodesic (product) distance in $\SS\times\C \times \P(\C^n)$.
For any nonzero matrix $A\in\Cnn$ (not necessarily of unit norm) we will write  
$$
   \dist_A((\la,v),(\la',v'))^2 := \frac{|\lambda-\lambda'|^2}{\|A\|_F^2}+\dpr(v,v')^2.
$$
Note that for any nonzero $A\in\Cnn$, $\dist_A$ is a distance function in $\C\times\Pn$ and if $A\in\SS$, then $\dist_A((\la,v),(\la',v'))=\dist((A,\la,v),(A,\la',v'))$.
\subsection{The Varieties $\V$, $\W$, $\Sigma'$ and $\Sigma$}\label{sec:2.2}
\label{subsec:SolVar}

We define the \textit{solution variety} for the eigenpair 
problem as
$$ 
  \V=\V_n := \left\{ (A, \lambda, v) \in   \Cnn \times \C
  \times \P ( \C^{n}): \ (A-\lambda \Id)v=0 \right\}.
$$

\begin{proposition}\label{prop:Vsmooth}
The solution variety $\V$ is a smooth submanifold of 
$\ambient$, of the same dimension 
as $\Cnn$. 
\end{proposition}

\proof
See~\cite[Proposition~2.2]{Armentano:13}.
\eproof

The set $\V$ inherits the Riemannian structure of the ambient space.
Associated to it there are natural projections: 
\smallskip

\addtocounter{equation}{1}
\begin{equation}\tag*{\raisebox{30pt}{(\theequation)}}
\begin{tikzpicture}
\path (-0.25,0.3) node[right]{$\V$};
\draw[->] (0.15,0) -- (1.15,-1) node[above =3mm]{$\pi_2$};
\draw[->] (-0.15,0) -- (-1.15,-1) node[above =3mm]{$\pi_1$};
\path (-1.9,-1.3) node[right]{$\Cnn$};
\path (0.5,-1.3) node[right]{$\C\times\P(\C^n)$.};
\end{tikzpicture}
\end{equation}
Because of Proposition~\ref{prop:Vsmooth}, 
the derivative $D\pi_1$ at $(A,\lambda,v)$ 
is a linear operator between spaces of equal dimension. 
The \emph{subvariety $\W$ of well-posed triples} is the subset of
triples $(A,\lambda,v)\in\V$ for which $D\pi_1 (A,\lambda,v)$ is an
isomorphism. In particular, when $(A,\lambda,v)\in\W$, the 
projection $\pi_1$ has a branch of its inverse (locally defined) 
taking $A\in\Cnn$ to $(A,\lambda,v)\in\V$. %This branch of  $\pi_1^{-1}$ is called the \emph{solution map} at $(\lambda,v)$. 

For $v\in\P(\C^n)$ we denote by 
$v^\perp=\{x \in \C^n \mid \langle x,v\rangle =0 \}$ 
the tangent space to $\P(\C^n)$ at $v$. Let   
$P_{v^\perp}\colon\C^n\to v^\perp$ be the orthogonal projection. 
Given $(A,\lambda,v)\in\ambient$, we 
let $A_{\lambda,v}\colon v^\perp\to v^\perp$ be the linear operator 
given by 
\begin{equation}\label{eq:defAlv}
  A_{\lambda,v}:=P_{v^\perp} \circ (A-\lambda \Id)|_{v^\perp} .
\end{equation}
If we choose a representative such that $\|v\|=1$ and we assume that $A_{\lambda,v}$ is invertible, then we have
\begin{equation}\label{eq:defAlvbis}
i_{\C^n} A_{\lambda,v}P_{v^\perp}=(\Id-vv^*)(A-\lambda \Id)(\Id-vv^*),
\end{equation}
and
\begin{equation}\label{eq:defAlvbisbis}
i_{\C^n} A_{\lambda,v}^{-1}P_{v^\perp}=\big((\Id-vv^*)(A-\lambda \Id)(\Id-vv^*)\big)^\dagger,
\end{equation}
where $i_{\C^n}:v^\perp\to\C^n$ is the inclusion and $^\dagger$ denotes Moore--Penrose pseudoinverse.

The set of well-posed triples is exactly
\begin{equation}\label{eq:Alv}
    \W=\{(A,\lambda,v)\in\V:\;A_{\lambda,v}\,\mbox{ is invertible}\},
\end{equation}
see~\cite[Lemma 2.7]{Armentano:13}). 
We define $\Sigma':=\V \setminus \W$ to be the variety of 
\emph{ill-posed triples}, and $\Sigma=\pi_1(\Sigma')\subset\Cnn$ the
\emph{discriminant variety}, i.e., the subset of {\em ill-posed inputs}.

\begin{remark}
From (\ref{eq:Alv}) it is clear that the subset $\Sigma\rq{}$ is the
set of triples $(A,\lambda,v)\in \V$ such that $\lambda$ is an
eigenvalue of $A$ of algebraic multiplicity at least~2. 
Hence $\Sigma$ is the set of matrices $A\in\Cnn$ with multiple
eigenvalues, and for $A\in\Cnn\setminus\Sigma$, 
the eigenvalues of $A$ are pairwise different and 
$\pi_1^{-1}(A)$ is the set of triples
$(A,\lambda_1,v_1),\ldots,(A,\lambda_n,v_n)$, where 
$(\lambda_i,v_i)$, $i=1,\ldots,n$, are the eigenpairs of~$A$.
\end{remark}

\begin{proposition}\label{prop:codim}
The discriminant variety $\Sigma\subset\Cnn$ is a complex 
algebraic hypersurface. Consequently, for all $n\geq2$, 
we have $\dim_{\R}\Sigma=2n^2-2$. 
\end{proposition}

\proof
See~\cite[Proposition~20.18]{Condition}.
\eproof

\medskip

\subsection{Unitary invariance}\label{sec:unitaryaction}

Let $\Un$ be the group of $n\times n$ unitary matrices.  The group
$\Un$ naturally acts on $\prc$ by $U\cdot[w]:=[Uw]$. 
In addition, $\Un$ acts on
$\Cnn$ by conjugation (i.e., $U\cdot A:=UAU^{-1}$), and on
$\ambientu$ by $U\cdot (A,\lambda):= (UAU^{-1},\lambda)$. 
These actions
define an action on the product space $\ambient$, namely,
\begin{equation}\label{eq:UactionPP}
    U\cdot(A,\lambda,v):= (UAU^{-1},\lambda,Uv).
\end{equation}

\begin{remark}
The varieties $\V$, $\W$, $\Sigma\rq{}$, and $\Sigma$,
are invariant under the action of $\Un$ (see~\cite{Armentano:13} for details).
\end{remark}

\subsection{Condition of a triple}

In a nutshell, condition numbers measure the worst possible 
output error resulting from a small perturbation on the input data. 
More formally, a condition number is the operator 
norm of the derivative of a solution map such as the branches of $\pi^{-1}$ 
mentioned in~\S\ref{subsec:spaces} above, 
(see~\cite[\S14.1.2]{Condition} for a general exposition).  

In the case of the eigenpair problem, one can define two 
condition numbers for eigenvalue and eigenvector, respectively, 
and formulas for both of them have been known at 
least since~\cite{vanLoan}. 
Armentano has shown that one can merge the two 
in a single one (see \S3 in~\cite{Armentano:13} for details).    
Deviating slightly from~\cite{Armentano:13}, 
we define the \emph{condition number} of
$(A,\lambda,v)\in\mnc\times\C\times\C^n$ as
\begin{equation}\label{eq:defmu}
    \mu(A,\lambda,v) := \|A\|_F \|A_{\lambda,v}^{-1} \|,
\end{equation}
(or $\infty$ if $A_{\lambda,v}$ is not invertible) where $\|\cdot\|$
is the operator norm. This coincides with $\mu_v(A,\lambda,v)$ 
in~\cite{Armentano:13}. 
Note that from~\eqref{eq:defAlvbisbis}, if 
$(A,\lambda,v)$ is such that $\mu(A,\lambda,v)<\infty$,
and $\|v\|=1$ then:
\begin{equation}\label{eq:mualternative}
    \mu(A,\lambda,v)=\|A\|_F
   \left\|\left((\Id-vv^*)(A-\lambda\Id)(\Id-vv^*)\right)^\dagger\right\|.
\end{equation}

\begin{remark}\label{rem:muscaling}
The condition number $\mu$ is invariant under the action of the 
unitary group $\Un$, i.e.,
$\mu(UAU^{-1},\lambda,Uv)=\mu(A,\lambda,v)$ for all $U\in\Un$, and  scale 
invariant on the first two components, i.e.,
$\mu(sA,s\lambda,v)=\mu(A,\lambda,v)$ for all $s\in\C\setminus\{0\}$.
\end{remark}

\begin{lemma}[Lemma~3.8 in~\cite{Armentano:13}]\label{le:lb_mu}
For $(A,\la,v)\in\V$ we 
have $\mu(A,\la,v) \ge \frac{1}{\sqrt{2}}$. \hfill $\Box$ %\eproof
\end{lemma}

The essence of condition numbers is that they measure how much may  
outputs vary when inputs are slightly perturbed. The following result, 
which we will prove in \S\ref{sec:condition_property}, 
quantifies this property for $\mu$. 

\begin{proposition}\label{prop:dotzdotf-spherical} 
Let $\Gamma:[0,1]\to \V$, $\Gamma(t)=(A_t,\lambda_t,v_t)$ be a 
smooth curve such that $A_t$ 
lies in the unit sphere of $\Cnn$, for all $t$. 
Write $\mu_t:=\mu(\Gamma(t))$. Then we have, for all $t\in[0,1]$, 
$$
    |\dot{\lambda_t}|\leq \sqrt{1+\mu_t^2}\;\|\dot{A_t}\|,
    \qquad
    \|\dot{v_t}\|\leq \mu_t\;\|\dot{A_t}\|.
$$
In particular,
$$
  \big\|\dot{\Gamma}(t)\big\|\leq 
  \sqrt{6}\;\mu_t\;\|\dot{A_t}\|.
$$
\end{proposition}

Condition numbers are generally associated to input data. 
In the case of a problem with many possible solutions 
(of which returning an eigenpair of a given matrix is a 
clear case) one can derive the condition of a data from a notion 
of condition for each of these solutions. A discussion of 
this issue is given in~\cite[\S6.8]{Condition}. For the 
purposes of this paper, we will be interested in two such derivations: 
the {\em maximum condition number} of $A$,
$$
    \mum(A):=\max_{1\leq j\leq n}\mu(A,\lambda_j,v_j),
$$
and the {\em mean square condition number} of $A$,
$$
   \mu_{\av}(A):=\left(\frac1n\sum_{j=1}^n
   \mu^2(A,\lambda_j,v_j)\right)^{\frac12}=\left(\frac1n\sum_{j=1}^n
  \|A\|_F^2\|A_{\lambda_j,v_j}^{-1}\|^2\right)^{\frac12}.
$$

Condition numbers themselves vary in a controlled manner. 
The following Lipschitz property and its corollary make this 
statement precise.

\begin{theorem}\label{th:lipsch}
Let $A,A'\in\C^{n\times n}$ be such that $\|A\|_F = \|A'\|_F = 1$, let 
$v,v'\in\C^n$ be nonzero, and let  $\la,\la'\in\C$. Suppose that
$$
 \mu(A,\la,v)\;\dist((A,\la,v),(A',\la',v'))\;\le\;  
 \frac{\e}{4\sqrt{3}}  \quad \text{for some $\e\in(0,1)$.}
$$
Then we have 
$$
 \frac{1}{1+\e}\,\mu (A,\lambda,v)\ \le\ \mu(A',\la',v') \ \le\  
 (1+\e)\, \mu(A,\la,v) .
$$
\end{theorem}

\begin{corollary}\label{cor:b}
Let $A\in\SS$, $A\not\in\Sigma$, and let $A'\in\SS$ be 
such that
\[
    \|A-A'\|_F\leq \frac{\e}{50\mum^2(A)},
     \quad \text{for some $\e\in(0,1)$.}
\]
Then, $A'\not\in\Sigma$ and 
$(1+\e)^{-1}\mum(A)\leq \mum(A')\leq (1+\e)\mum(A)$.
\end{corollary}

We give the proofs of Theorem~\ref{th:lipsch} and Corollary~\ref{cor:b} 
in \S\ref{se:lipschitz} below. 

We close this paragraph observing that restricted to  
the class of normal matrices, the condition number $\mu$  
admits the following elegant characterization.

\begin{lemma}[Lemma~3.12 in~\cite{Armentano:13}]\label{lem:char}
Let $A\in\Cnn\setminus\Sigma$ be normal, and let 
$(\lambda_1,v_1),\ldots,(\lambda_n,v_n)$ be its eigenpairs. Then
\begin{equation}\tag*{\qed}
  \mu(A,\lambda_1,v_1) = 
  \frac{\|A\|_F}{\min_{j=2,\ldots,n}|\lambda_j-\lambda_1|}.
\end{equation}
\end{lemma}

%%%%%%%%%%%%%%

\subsection{Newton's method and approximate eigenpairs}

For a nonzero matrix $A\in\Cnn$, we 
define the \emph{Newton map} associated to $A$, 
$$
    N_A:\C\times(\C^n\setminus\{0\})\to \C\times(\C^n\setminus\{0\}), 
$$
by
\[
  N_A\binom{\lambda}{v}=
 \binom{\lambda}{v}-\binom{\dot\lambda}{\dot v},
 \quad\text{where}\quad 
 \binom{\dot \lambda}{\dot v}=
 \left(DF_A(\lambda,v)\mid_{\C\times v^\perp}\right)^{-1}
 F_A\binom{\lambda}{v}
\]
and $F_A(\lambda,v)=(A-\lambda \Id)v$ is the evaluation map. 
This is a rational map (it is only defined on an open subset of 
$\C\times(\C^n\setminus\{0\})$). It was introduced 
in~\cite{Armentano:13} as a 
Newton-like operator associated to 
the evaluation map $F_A$, and the following formulas were 
obtained for $\dot v$ and $\dot\lambda$ (recall the definition 
of $A_{\lambda,v}$ from~\eqref{eq:defAlv}):
\begin{equation}\label{eq:newton}
  \dot v  ={A_{\lambda,v}}^{-1}\,
  P_{v^\perp}Av, \qquad 
  \dot\lambda =\frac{\pes{(\lambda\,\Id-A)(v-\dot v)}{v}}{\pes vv}.
\end{equation}
The map $N_A$ is defined for every 
$(\lambda,v)\in\C\times(\C^n\setminus\{0\})$ such that 
$A_{\lambda,v}$ is invertible. 
\begin{equation}\label{eq:hom}
 N_{zA}^k\binom{z\lambda}{v}=\begin{pmatrix}z&0\\0&1\end{pmatrix}N_A^k\binom{\lambda}{v},
\end{equation}
where the superindex $^k$ means $k$ iterations. See~\cite[Sec. 4]{Armentano:13} for 
more details. 

The notion of approximate solution as a point where Newton's method 
converges to a true solution immediately and quadratically fast was 
introduced by Steve Smale~\cite{Smale86}. It allows to elegantly talk
about polynomial time without dependencies on pre-established 
accuracies. In addition, these approximate solutions are ``excellent 
approximations'' (as mentioned in the statement of the main results) 
in a very strong sense: the distance to the exact solution dramatically 
decreases with a single iteration of Newton's method. In the context of 
eigenpair computations this concept is settled as follows.

\begin{definition}\label{def:app_eigen}
Given $(A,\lambda,v)\in\W$  we say that 
$(\zeta,w)\in\C\times(\C^n\setminus\{0\})$ is an 
{\em approximate eigenpair} of $A$ with associated eigenpair $(\lambda,v)$ 
when for all $k\geq 1$ the $k$th iterate $N_A^k(\zeta,w)$ 
of the Newton map at $(\zeta,w)$ is well defined and satisfies 
$$
  \dist_A\big((N_A^k(\zeta,w)),(\lambda,v)\big)\leq
  \left(\frac12\right)^{2^k-1} \dist_A\big((\zeta,w),(\lambda,v)\big).
$$
% If $\|A\|_F\neq 1$ we say that $(\zeta,w)$ is an approximate eigenpair
% of $A$ with associated eigenpair $(\lambda,v)$ if $(\zeta/\|A\|_F,w)$
% is an approximate eigenpair of $A/\|A\|_F$ with associated eigenpair
% $(\lambda/\|A\|_F,v)$.
\end{definition} 
The following result estimates, in terms of the 
condition of an eigenpair, the radius of a ball of approximate 
eigenpairs associated to it. For a complete proof 
see~\cite[Theorem~5]{Armentano:13}.

\begin{theorem}\label{th15.1}
There is a universal constant $c_0>1/5$ with the following property. 
Let $(A,\lambda,v)\in\W$ with $\|A\|_F=1$ and let 
$(\zeta,w)\in\C\times(\C^n\setminus\{0\})$. If 
$$
    \dist_A\big((\lambda,v),(\zeta,w)\big)\leq
   \frac{c_0}{\mu(A,\lambda,v)},
$$ 
then $(\zeta,w)$ is an approximate eigenpair of $A$ with 
associated eigenpair $(\lambda,v)$.% One may choose $c_0=0.2881$.
\end{theorem}
It is a simple exercise to check that for any nonzero $z\in\C$, $(\zeta,w)$ is an approximate zero of $A$ with associated zero $(\lambda,v)$ if and only if $(z\zeta,w)$ is an approximate zero of $zA$ with associated zero $(z\lambda,v)$. So from the point of view of analyzing the effect of the Newton methods we may pick whatever scaling is convenient. For us it will be convenient to assume that $\|A\|_F=1$ which we will do in the following.
\medskip

\proofof{Theorem \ref{th15.1}}
 Note that \cite[Theorem~5]{Armentano:13} is the same result with $c_0=0.2881$, but the definition of the condition number in \cite{Armentano:13} is slightly different from ours. More exactly, if we denote by $\kappa(A,\lambda,v)$ the condition number defined in \cite{Armentano:13} then we have $\kappa(A,\lambda,v)=\max(1,\mu(A,\lambda,v))$. Theorem \ref{th15.1} is hence true with $\kappa$ in the place of $\mu$ and $c_0=0.2881$. However, from Lemma \ref{le:lb_mu} we know that $\mu(A,\lambda,v)\geq2^{-1/2}$ which readily implies $\kappa(A,\lambda,v)\leq\sqrt{2}\mu(A,\lambda,v)$. Theorem \ref{th15.1} now follows from the fact that $0.2881>\sqrt{2}/5$.
\eproof

\begin{remark}
We note that $N_A(\zeta,w)$ can be computed from the matrix 
$A$ and the pair $(\zeta,w)$ in $\Oh(n^3)$ operations, since 
the cost of this computation is dominated by that of 
inverting a matrix (or simply solving a linear system).
\end{remark}
\medskip

\subsection{Gaussian Measures on $\Cnn$}
\label{sec:gaussian}

Let $\sigma>0$. We say that the complex random variable 
$Z=X+\sqrt{-1}Y$ has distribution 
$\mcN_{\C}(0,\sigma^2)$ when the real part $X$ and the
imaginary part $Y$ are independent and identically distributed
(i.i.d.) drawn from $\mcN(0,\frac{\sigma^2}{2})$, i.e., they are 
Gaussian centered random variables with variance 
$\frac{\sigma^2}{2}$. 

If $Z\sim\mcN_{\C}(0,\sigma^2)$ then its
density $\varphi\colon\C\to\R$ with respect to the Lebesgue measure 
is given by 
\begin{equation*}%\label{eq:densidadcompleja}
   \varphi(z):=\frac{1}{\pi\sigma^2}e^{-\frac{|z|^2}{\sigma^2}}. 
\end{equation*}

We will write $v\sim\mcN_{\C^n}(0,\sigma^2)$ 
to indicate that the vector
$v\in\C^n$ is random with i.i.d.~coordinates drawn from
$\mcN_{\C}(0,\sigma^2)$. Also, we say that 
$A\in\C^{m\times n}$ is (isotropic) Gaussian
and we write $A\sim\mathcal{N}_{\C^{m\times n}}(0,\sigma^2)$,
if its entries are i.i.d. Gaussian random
variables. The resulting probability space is sometimes called 
the {\em{Ginibre ensemble}}. 

If $\hA\in\C^{m\times n}$ and 
$G\sim\mathcal{N}_{\C^{m\times n}}(0,\sigma^2)$,  
we say that the random matrix $A=G+\hA$ has the 
{\em Gaussian distribution centered at $\hA$}, and we write 
$A\sim\mcN_{\C^{m\times n}}(\hA,\sigma^2)$. The density of this 
distribution is given by
$$
  \varphi_{m\times n}^{\hA,\sigma}(A):=\frac{1}{(\pi\sigma^2)^{mn}}\,
  e^{-\frac{\|A-\hA\|_F^2}{\sigma^2}}. 
$$
For conciseness, we will sometimes write 
$A\sim\mathcal{N}_{\C^{m\times n}}$ in the particular case where 
$\hA=0$ and $\sigma=1$.

Crucial in our development is the following result giving 
a bound on the average condition for Gaussian matrices 
arbitrarily centered. Its statement is similar to the main technical 
result in~\cite[Thm.~3.6]{Condition}. We will prove 
it in \S\ref{se:muave}. 

For technical reasons we will be interested in the following 
variation of $\mu$:
$$
\mu_F(A,\lambda,v):=\|A\|_F\|A_{\lambda,v}^{-1}\|_F
$$
(note, we only replaced $\|A_{\lambda,v}^{-1}\|$ by 
$\|A_{\lambda,v}^{-1}\|_F$) and the corresponding  
$$
\muFa(A):=\left(\frac1n\sum_{j=1}^n
   \mu_F^2(A,\lambda_j,v_j)\right)^{\frac12}.
%=\left(\frac1n\sum_{j=1}^n
%  \|A\|_F^2\|A_{\lambda_j,v_j}^{-1}\|_F^2\right)^{\frac12}.
$$

\begin{theorem}\label{th:mu2-bound}
For $\hA\in\Cnn$ and $\s>0$ we have
$$
   \Exp_{A\sim \mcN_{\C^{n\times n}}(\hA,\s^2)}
   \Big(\frac{\muFa^2(A)}{\|A\|_F^2} \Big)\ \le\ \frac{n}{\s^2}.
$$
Moreover, for $A$ chosen with the uniform distribution 
$\msU(\SS)$ in the unit sphere $\SS$ of $\mnc$ we have:
\[
%   \Exp_{A\sim \msU(\SS)}\left(\mu_{\av}^2(A)\right)\leq  
   \Exp_{A\sim \msU(\SS)}\left(\muFa^2(A)\right)\leq n^3.
\]
\end{theorem}

\begin{remark} 
\begin{description}
\item[(i)]
We note that no bound on the norm of $\hA$ is required in the 
first claim of Theorem~\ref{th:mu2-bound}.
Indeed, using $\muFa(s A)=\muFa(A)$,
it is easy to see that the assertion for  a pair $(\hA,\s)$ implies the
assertion for $(s\hA,s\s)$, for any $s>0$.
% \item[(ii)]
% Because of Proposition~\ref{prop:codim}, with probability one, 
% matrices drawn from $\mcN_{\C^{n\times n}}(\hA,\s^2)$ have all its
% eigenvalues different. Therefore the expected value in 
% Theorem~\ref{th:mu2-bound} is well-defined.  
\item[(ii)] 
It is remarkable that if we change Gaussian matrices to some classes 
of structured matrices, the expected value of the condition number 
can be very high, see for example~\cite{George} and references 
therein.
\end{description}
\end{remark}

\subsection{Truncated Gaussians and smoothed analysis}

For $T,\sigma>0$, we define the 
{\em truncated Gaussian} $\mcN_{\C^{n\times n},T}(0,\sigma^2)$
on $\Cnn$ to be the distribution  
given by the density 
\begin{equation}\label{eq:truncated}
   \rho^{\sigma}_T(A)=
    \left\{\begin{array}{ll}
     \frac{\varphi_{n\times n}^{0,\sigma}(A)}{P_{T,\s}} & 
     \mbox{if $\|A\|_F\leq T$,}\\
     0 &\mbox{otherwise,}\end{array}\right.
\end{equation}
where $P_{T,\s}:=\Prob_{A\sim \mathcal{N}_{\C^{n\times n}}(0,\sigma^2)}
\{\|A\|_F\leq T\}$, 
and, we recall, $\varphi_{n\times n}^{0,\sigma}$ is the density of 
$\mathcal{N}_{\C^{n\times n}}(0,\sigma^2)$. For the rest of this paper we 
fix the threshold $T:=\sqrt{2}\,n$. The fact that 
$\|A\|_F^2$ is chi-square distributed with $2n^2$ degrees of freedom, 
along with~\cite[Corollary~6]{choi:94} yields the following result.

\begin{lemma}\label{lem:X}
We have $P_{T,\s} \ge \frac12$
for all $0<\sigma\leq 1$.\eproof
\end{lemma}

The space $\C^{n\times n}$ of matrices with the Frobenius norm 
and the space $\C^{n^2}$ with the canonical Hermitian product are 
isomorphic as Hermitian product spaces. Hence, the 
Gaussian $\mcN_{\C^{n\times n}}(0,\s^2)$ 
on the former corresponds to the 
Gaussian $\mcN_{\C^{n^2}}(0,\s^2)$ 
on the latter, and we can deduce invariance 
of $\mathcal{N}_{\C^{n\times n}}(0,\sigma^2)$ under the action 
of $\mcU_{n^2}$ 
(in addition to that for conjugation under $\mcU_n$ discussed 
in~\S\ref{sec:unitaryaction}), and the same is true for the truncated 
Gaussian. In particular, the pushforward 
of both distributions for the projection 
$\Cnn\setminus\{0\}\to\SS$, $A\mapsto\frac{A}{\|A\|_F}$, is 
the uniform distribution $\msU(\SS)$ 
(see~\cite[Chapter~2]{Condition} for details) and we have
\begin{equation}\label{eq:truncating}
  \Exp_{A\sim\mathcal{N}_{\C^{n\times n}}(0,\sigma^2)} F(A) \,=\,
 \Exp_{A\sim\mathcal{N}_{\C^{n\times n},T}(0,\sigma^2)} F(A)  \,=\,
  \Exp_{A\sim\msU(\SS)} F(A),
\end{equation}
for any measurable scale invariant function $F:\Cnn\to[0,\infty)$.
\medskip

Complexity analysis has traditionally been carried out either 
in the {\em worst-case} or in an {\em average-case}. More generally, 
for a function $F:\R^m\to\R_+$ (some measure for the computational 
cost of solving an instance in $\R^m$), the former amounts to the 
evaluation of $\sup_{a\in\R^m}F(a)$ and the latter to that of 
$\Exp_{a\sim\msD}F(a)$ for some probability distribution 
$\msD$ on $\R^m$. Usually, $\msD$ is taken to be 
the standard Gaussian in the input space. With the beginning of 
the century, Daniel Spielman and Shang-Hua Teng introduced 
a third form of analysis, {\em smoothed analysis}, which is meant to 
interpolate between worst-case and average-case. We won't 
elaborate here on the virtues of smoothed analysis; a defense 
of these virtues can be found, e.g., 
in~\cite{ST:02,ST:09} or in~\cite[\S2.2.7]{Condition}. We 
will instead limit ourselves to the description of what smoothed 
analysis is and which form it will take in this paper. 

The idea is to replace the two operators above (supremum 
and expectation) by a combination of the two, namely,
$$
   \sup_{\hat a\in\R^m} \Exp_{a\sim\msD(\hat a,\sigma)} F(a) 
$$
where $\msD(\hat a,\sigma)$ is a distribution ``centered'' at 
$\hat a$ having $\sigma$ as a measure of dispersion. A  
typical example is the Gaussian $\mcN(\hat a,\sigma^2)$. 
Another example, used for scale invariant functions $F$, 
is the uniform measure on a spherical cap centered at 
$\hat a$ and with angular radius $\sigma$ on the unit 
sphere $\SS(\R^m)$ (reference~\cite{Condition} exhibits  
smoothed analyses for both examples of distribution). 
In this paper we will perform a smoothed analysis with respect 
to a truncated Gaussian. More precisely, we will be interested in 
quantities of the form 
\begin{equation*}
   \sup_{\hA\in\Cnn} \;\Exp_{A\sim\mcN_{\C^{n\times n},T}(\hA,\sigma^2)} 
  F(A) 
\end{equation*}
where $F$ will be a measure of computational cost for the 
eigenpair problem. 
We note that, in addition to the usual dependence on $n$, 
this quantity depends also on~$\sigma$ and tends to $\infty$ 
when $\sigma$ tends to $0$. When $F$ is scale invariant, 
as in the case of $\mu_{\av}$ or $\mum$, it is customary 
to restrict attention to matrices of norm 1. That is, to 
study the following quantity:
\begin{equation}\label{eq:smoothed}
   \sup_{\hA\in\SS}  \;\Exp_{A\sim\mcN_{\C^{n\times n},T}
   (\hA,\sigma^2)} F(A) .
\end{equation}

\subsection{The eigenpair  continuation algorithm}
\label{subsec:ALH}

We are ready to describe the main algorithmic construct 
in this paper. When dealing with algorithms it will be more convenient 
to view the solution variety as the corresponding subset of
$\Cnn\times\C\times(\C^n\setminus\{0\})$, which, abusing notation, 
we still denote by $\V$. 

Given two matrices $B,B_0\in\SS$, $B\neq\pm B_0$, 
let $\a:=\dS(B_0,B)\in(0,\pi)$  be the spherical distance (i.e. the angle) 
from $B_0$ to $B$, and let
\begin{equation}\label{eq:mathcalL}
   \mathcal{L}_{B_0,B}=\{B_s:0\leq s\leq \a\}
\end{equation}
be the portion of the great circle in $\SS$, parametrized by 
arc-length, joining $B_0$ and $B$, so $B_\a=B$. By abuse of notation, 
for any $A_0,A\in\mnc$ such that $A_0,A$ are not $\R$-linearly 
dependent, we simply write 
\[
    d_\SS(A_0,A):=d_\SS\left(\frac{A_0}{\|A_0\|_F},\frac{A}{\|A\|_F}\right)
   \mbox{ and }\quad 
   \mathcal{L}_{A_0,A}:=\mathcal{L}_{\frac{A_0}{\|A_0\|_F},\frac{A}{\|A\|_F}}.
\]
The following neighborhood of the set of well-posed solution triples plays 
a distinguished role in the continuation algorithm. Its definition uses
a constant $c_*$ which we will later take to be $10^{-4}$ 
(cf.~\S\ref{sec:lip-beta}). 

\begin{definition}\label{def:fund_neigh}
The {\em initial neighborhood} of the set $\W$ is the set 
$$
   \wW:=\Big\{(A,\zeta,w)\mid \exists (\lambda,v) \mbox{ s.t. } 
    (A,\lambda,v)\in\W \mbox{ and } 
    \dist_A((\zeta,w),(\lambda,v))\leq 
   \frac{c_*}{\mu_{\max}(A)}\Big\}.
$$  
\end{definition}

The choice of $c_*$ implies that the hypothesis of 
Theorem~\ref{th15.1} is satisfied and, hence, that 
$(\zeta,w)$ is an approximate eigenpair of $A$ with associated 
eigenpair $(\lambda,v)$.

Suppose that we are given an an \emph{initial triple} $(B_0,\zeta_0,w_0)\in\wW$, $B_0\in\SS$ and an input matrix $B\in\SS\setminus\{\pm B_0\}$.
Let $\lambda_0,v_0\in\C\times\P(\C^n)$ be as in 
Definition~\ref{def:fund_neigh}.  
As a consequence of the inverse function theorem, if 
$(\mathcal{L}_{B_0,B}\setminus\{B_0\}) \cap \Sigma=\emptyset$, 
then the map $s\mapsto B_s$
can be uniquely extended to a continuous map
\begin{equation}\label{eq:curva}
  [0,\a]\to \V,\quad s \mapsto (B_s,\lambda_s,v_s), 
\end{equation}
We call this map the {\em lifting} of $\mathcal{L}_{B_0,B}$ with
origin $(B_0,\lambda_0,v_0)$.
We can try to approximate the eigenpair $(\lambda_\a,v_\a)$ of $B$ by
following the lifting of $\mathcal{L}_{B_0,B}$. To this end we can
differentiate $B_sv_s-\lambda_sv_s=0$ w.r.t. $s$. This produces an
Initial Value Problem (IVP) whose solution can be approximated by any
standard numerical ODE solver. The main ingredient for the complexity
estimate is the number of points in the discretization of $[0,\a]$
needed to approximate the solution of the IVP.

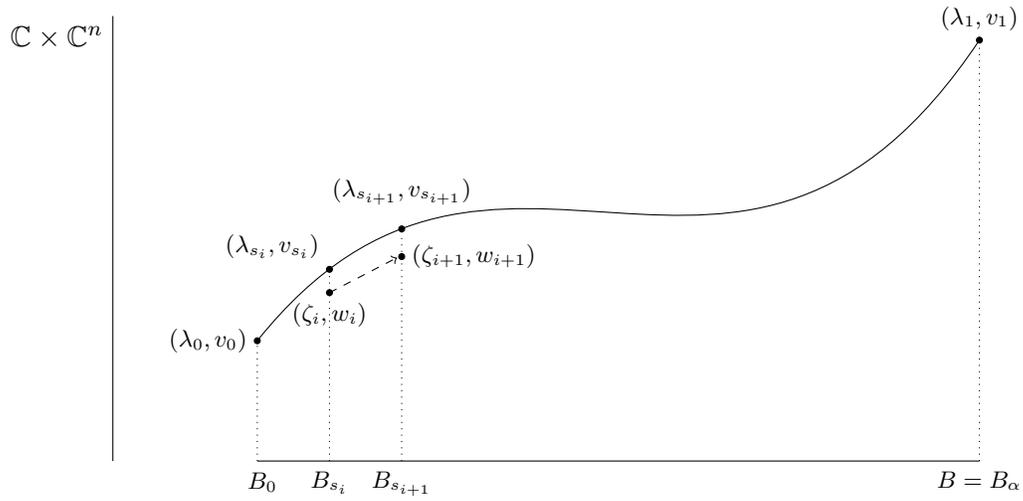
\begin{figure}[H]
\begin{center}
\begin{tikzpicture}[scale=1.6]
\def\mypoints{0.7}
\filldraw (0,0) circle(\mypoints pt); 
\filldraw (6,2.5) circle(\mypoints pt); 
\draw (0,0) node[left]{\footnotesize $(\lambda_0,v_0)$} .. controls (2,2.5) 
and (4,-0.5) .. (6,2.5) node[above]{\footnotesize $(\lambda_1,v_1)$};
\draw (0,-1)  node[below=2.8mm, left=-4mm]{\footnotesize $B_0$} -- (6,-1)  
node[below]{\footnotesize $B=B_{\a}$};
\draw[dotted] (0,-1) -- (0,0);
\draw[dotted] (0.6,-1)  node[below]{\footnotesize $B_{s_i}$}-- (0.6,0.6) 
node[above left]{\footnotesize $(\lambda_{s_i},v_{s_i})$};
\filldraw (0.6,0.4) circle(\mypoints pt)  node[below]{\footnotesize $(\zeta_i,w_i)$}; 
\filldraw (0.6,0.594) circle(\mypoints pt); 
\draw[dotted] (1.2,-1)  node[below]{\footnotesize $
B_{s_{i+1}}$} -- (1.2,0.9); 
\draw (1.2,1.05) 
node[above]{\footnotesize $(\lambda_{s_{i+1}},v_{s_{i+1}})$};
\filldraw (1.2,0.7) circle(\mypoints pt) 
node[right]{\footnotesize $(\zeta_{i+1},w_{i+1})$}; 
\filldraw (1.2,0.93) circle(\mypoints pt); 
\draw[dotted] (6,-1) -- (6,2.5);
\draw[dashed] [->](0.6,0.4) -- (1.16,0.69);
\draw (-1.2, -1) -- (-1.2,2.7) node[below left]{$\C\times\C^n$};
\end{tikzpicture}
\end{center}
  \caption{The continuation of the solution path.}
  \label{fig:homotopy}
\end{figure}

Formalizing this idea to get an actual guarantee of convergence is a
nontrivial task; only a non-constructive method has been described
in~\cite{Armentano:13} following the ideas in~\cite{Shub2007}.  We now
describe how to algorithmically construct a numerically stable method 
for this task: subdivide the interval $[0,\a]$ into subintervals with
extremities at $0=s_0<s_1<\cdots<s_K=\a$
and successively compute approximations~$(\zeta_i,w_i)$ 
of~$(\lambda_{s_i},v_{s_i})$, starting with $(\zeta_0,w_0)$ 
and then using Newton's method. 
To ensure that these are good approximations, we actually  
want to ensure that for all $i$, $(\zeta_i,w_i)$ is 
an approximate eigenpair of $B_{s_{i+1}}$.
Figure~\ref{fig:homotopy} attempts to convey the general idea. 

The algorithmic counterpart of this idea is the following.
\bigskip

\algoritmo
\begin{algorithm}\label{alg:ECsphere}
{\sf Path-follow}\\
\inputalg{$A,A_0\in \Cnn\setminus\{0\}$, 
and $(\zeta_0,w_0)\in\C\times\C^n$}\\
\specalg{$(A,\zeta_0,w_0)\in\wW$, $\|w_0\|=1$, $\R A_0\neq\R A$}\\
\bodyalg{
redefine $A_0:=A_0/\|A_0\|_F$, $A_1:=A/\|A\|_F$\\[2pt]
$\a:=\dS(A_0,A_1)$, 
$s_0:=0$, $B:=A_0$, $i=0$\\[2pt]
repeat\\[2pt]
\espacio $\dot A:=$ unit tangent vector (in the direction of the\\[2 pt]
\eeespacio parametrization) to $\mathcal{L}_{A_0,A_1}$ at $B$\\[2pt]
\espacio $\Delta s:=\;${\sf Choose\_step$(B,\dot A,\zeta_i,w_i)$}\\[2pt]
\espacio $s_{i+1}:=\min\{\a,s_i+\Delta s\}$\\[2pt]
\espacio $B:=$ point in $\mathcal{L}_{A_0,A_1}$ with 
$d_{\SS}(A_0,B)=s$\\[2pt]
\espacio $(\zeta_{i+1},w_{i+1}):=N^3_{B}(\zeta_i,w_i)$ 
(three Newton iterations)\\[2pt]
\espacio if $|\zeta_{i+1}|>1$ then $\zeta_{i+1}=\zeta_{i+1}/|\zeta_{i+1}|$.\\[2pt]
\espacio $i:=i+1$\\[2 pt]
until $s= \alpha$\\[2pt]
return $(\zeta',w')=(\|A\|_F\zeta_i,w_i)$\\
}
\Output{$(\zeta',w')\in \C\times \C^{n}$}\\
\postcond{The execution halts if the lifting of 
$\mathcal{L}_{A_0,A}$ at $(\lambda,v)$ doesn't cut $\Sigma'$. 
In this case, $(A,\zeta',w')\in\wW$ and hence $(\zeta',w')$ is an approximate 
eigenpair of~$A$.}
\end{algorithm}
\falgoritmo

Algorithm {\sf Path-follow} is unambiguous except for the 
subroutine {\sf Choose\_step}, which will be described at the 
end of \S\ref{sec:onestep} below. We next state the 
main results for it. Recall, the point $s_\ell\in[0,\alpha]$ 
is the value of $s$ generated by \EC\ at the $\ell$th 
iteration. 

\begin{theorem}\label{thm:main_path_following}
Suppose that $\mathcal{L}_{A_0,A}=\{B_s\}_{s\in[0,\alpha]}$ (so
$B_0=A_0/\|A_0\|_F$ and $B_\alpha=A/\|A\|_F$) and assume that
$\mathcal{L}_{A_0,A}\cap\Sigma=\emptyset$. Then the algorithm \EC\
stops after at most $\lceil K\rceil$ steps where 
$K$ is given by 
$$ 
  K:=K(A,A_0,\zeta_0,w_0):=C\int_0^{\a} \mu(B_s,\lambda_s,v_s)\|(\dot
  B_s,\dot\lambda_s,\dot v_s)\|\,ds. 
$$ 
Here the pairs $(\lambda_s,v_s)$ are given by~\eqref{eq:curva}. 

More generally, let $q\in\Z, q\geq 0$. Then algorithm \EC\ stops 
after at most (the smallest integer greater than or equal to)
$$
   q+C\int_{s_q}^{\a} \mu(B_s,\lambda_s,v_s)\|(\dot
  B_s,\dot\lambda_s,\dot v_s)\|\,ds
$$
steps. The returned pair $(\zeta,w)$
is an approximate eigenpair of $A$ with associated eigenpair
$(\lambda_{\a},v_{\a})$. Here, $C\leq 3000$ is a universal constant.
\end{theorem}

For $0\leq a<b\leq \alpha$, the quantity
\begin{equation}\label{eq:condlength}
L_{\mu,a,b}(B_s,\lambda_s,v_s)=\int_{a}^{b}
     \mu(B_s,\lambda_s,v_s)\|(\dot B_s,\dot\lambda_s,\dot v_s)\|\,ds
\end{equation}
in Theorem~\ref{thm:main_path_following} is the length of the curve
$\{(B_s,\lambda_s,v_s):a\leq s\leq b\}$ in the so-called {\em
condition metric}. This is the metric that is obtained by pointwise
multiplying the natural metric in $\SS\times\C\times\Pn$ by the
condition number squared. We call $ L_{\mu,a,b}(B_s,\lambda_s,v_s)$ 
the {\em condition length} of this curve.

The proof of Theorem ~\ref{thm:main_path_following} 
is given in~\S\ref{sec:homotopy}. 

\begin{remark}\label{rmk:mu2}
From Proposition~\ref{prop:dotzdotf-spherical} , we have 
\[
    L_{\mu,s_q,\alpha}(B_s,\lambda_s,v_s)\leq  \sqrt{6} \int_{s_q}^{\a}
    \mu^2(B_s,\lambda_s,v_s)\,ds.
\]
\end{remark}

The following result gives an alternative bound to the number of steps.

\begin{corollary}\label{corollary:daigual}
Let $A_0,A\in\mnc$ be $\R$-linearly independent and consider the path
$A_t=(1-t)A_0+tA$ which satisfies $A_1=A$. If
$\mathcal{L}_{A_0,A}\cap\Sigma=\emptyset$ then, for any 
$q\in\Z, q\geq 0$, the algorithm \EC\ 
stops after at most $\lceil K\rceil$ steps where 
\[
  K:=q+\sqrt{6}C\|A_0\|_F\,\|A_1\|_F\int_{t_q}^1
  \frac{\mu^2(A_t,\lambda_t,v_t)}{\|A_t\|_F^2}\,dt
\]
steps, with
$$
    t_q:=\frac{\|A_0\|_F}
     {\|A_1\|_F(\sin\alpha \cot s_q-\cos\alpha)+\|A_0\|_F}.
$$
Here, $\lambda_t$, $v_t$ are the eigenvalue and eigenvector 
of $A_t$ defined by continuation.
\end{corollary}

\proof
From Theorem~\ref{thm:main_path_following} and Remark~\ref{rmk:mu2},
letting $C'=\sqrt{6}\,C$, the number of steps is at most the smallest
integer greater than or equal to
\[
   q+C' \int_{s_q}^\a\mu^2(B_s,\lambda_{B_s},v_{B_s})\,ds
  =q+C'\int_{s_q}^\a\mu^2(B_s,\lambda_{B_s},v_{B_s})\|\dot B_s\|_F\,ds,
\]
where $\lambda_{B_s}$ and $v_{B_s}$ are the eigenvalue and eigenvector 
of $B_s$ defined by continuation. Now, reparametrizing that spherical 
segment by $C_t=A_t/\|A_t\|_F$, $t\in[0,1]$, the integral does not change. 
The quantity above is thus equal to
\[
  q+C'\int_{t_q}^1\mu^2(C_t,\lambda_{C_t},v_{C_t})\|\dot C_t\|_F\,dt
  \underset{\text{Rmk.~\ref{rem:muscaling}}}{=}
  q+C'\int_{t_q}^1\mu^2(A_t,\lambda_t,v_t)\left\|\frac{d}{dt}
  \left(\frac{A_t}{\|A_t\|_F}\right)\right\|\,dt.
\]
Substituting $\dot A_t=A_1-A_0$ and $A_t=(1-t)A_0+tA_1$ in this 
last formula 
and with some elementary computations 
(see~\cite[Lemma~17.5]{Condition}) we conclude that
\[
  \left\|\frac{d}{dt}\left(\frac{A_t}{\|A_t\|_F}\right)\right\|
  \leq\frac{\|A_0\|_F\,\|A_1\|_F}{\|A_t\|_F^2}.
\]
The corollary follows.
\eproof

The inequality of Remark~\ref{rmk:mu2} implies that (up to constants) the upper bound for the number of
steps by an algorithm in terms of the condition length as in Theorem \ref{thm:main_path_following} is smaller than the upper bound in terms of the integral of the squared condition number as in Corollary \ref{corollary:daigual}. A similar situation applies in the context of polynomial system root finding. In this case implementations exist in both contexts  
see \cite{Bez5,BuCu11,BeltranLeykin2011,DedieuMalajovichShub}. The condition length algorithm is more subtle and the proof of correctness more difficult, both for the polynomial system and eigenvalue, eigenvector cases. So the temptation is to present condition number squared algorithms. Here for the first time we present a quantitative estimate of the improvement which is significant and which justifies presenting the more complex condition length algorithm. In Theorem \ref{th:randomhomotopy} a randomized algorithm is studied. The upper bound given by the condition length algorithm is $O(n^2)$ while an algorithgm with complexity given by the condition number squared would give $O(n^3)$.

In our main results we are interested in the cost of algorithms over  
random matrices $A$. {The following quantity ---the expected number of 
iterations of \EC\ for a given initial triple $(A_0,\lambda_0,v_0)$--- 
becomes essential, 
\[
      \aviter(A_0,\lambda_0,v_0):=\E_{A\sim\mathcal{N}_{\C^{n\times n}}}
      K(A,A_0,\lambda_0,v_0).
\]
We can also consider the smoothed number of iterations  
of \EC\ that results by drawing instead 
the input matrix $A$ from $\mcN_T(\hA,\sigma^2)$ 
where $\hA\in\SS$ is arbitrary. We thus define
$$
  \siter(A_0,\lambda_0,v_0,\sigma):=\sup_{\hA\in\SS}
   \E_{A\sim\mcN_{\C^{n\times n},T}(\hA,\sigma^2)} K(A,A_0,\lambda_0,v_0).
$$}

\begin{proposition}\label{prop:fixA_0}
For $A_0\in\Cnn$ we have 
\[
   \aviter(A_0,\lambda_0,v_0)=\Oh(n^4\mu^2(A_0,\lambda_0,v_0))
\]
and, for all $\sigma\in(0,1]$,  
\[
   \siter(A_0,\lambda_0,v_0,\sigma)=\Oh\Big(\frac{n^4\mu^2(A_0,\lambda_0,v_0)}{\sigma^2}\Big).
\]

\end{proposition}

\begin{remark}\label{rem:codimension}
Note that Proposition~\ref{prop:fixA_0} ensures that 
\EC\ halts in finite time with probability 1 
even if $A_0\in\Sigma$, as long as 
$\mu(A_0,\lambda_0,v_0)$ is finite.
The algorithm does a first step that depends on 
$\mu(A_0,\lambda_0,v_0)$ only and advances to 
$B_{s_1}$ with $s_1>0$. After that, 
the fact that the real codimension of $\Sigma$ in $\Cnn$ is 2 
(shown in Proposition~\ref{prop:codim}) ensures that, 
almost surely, 
$(\mathcal{L}_{A_0,A}\setminus\{A_0\})\cap\Sigma=\emptyset$. 
Therefore, with probability 1, none of the 
matrices $B_s$ is in $\Sigma$ and the integral over 
$[s_1,\alpha]$ in Theorem~\ref{thm:main_path_following} 
is finite. 
\end{remark}

To compute all the eigenpairs from an initial matrix $A_0$ and its
eigenpairs $(\lambda^{(1)},v^{(1)}), \ldots,(\lambda^{(n)},v^{(n)})$
we may proceed by following the $n$ paths corresponding to taking 
these eigenpairs in the initial triples. In this case, we will 
be interested in the quantities 
\begin{align*}
      \aviterall(A_0):= &\sum_{i=1}^n\aviter(A_0,\lambda^{(i)},v^{(i)})\\=&\E_{A\sim\mathcal{N}_{\C^{n\times n}}}
      \sum_{i=1}^n K(A,A_0,\lambda^{(i)},v^{(i)}),
\end{align*}
and
\[
  \siterall(A_0,\sigma):=\sup_{\hA\in\SS}
   \E_{A\sim\mcN_{\C^{n\times n},T}(\hA,\sigma^2)} 
   \sum_{i=1}^n K(A,A_0,\lambda^{(i)},v^{(i)}).
\]
For these quantities we prove the following result.

\begin{proposition}\label{prop:fixA_0_all}
For $A_0\in\Cnn$ we have 
\[
   \aviterall(A_0) = \Oh(n^4\mum^2(A_0))
\]
and, for all $\sigma\in(0,1]$,  
\[
   \siterall(A_0,\sigma) = \Oh\Big(\frac{n^4\mum^2(A_0)}{\sigma^2}\Big).
\]

\end{proposition}

We prove Propositions~\ref{prop:fixA_0} and~\ref{prop:fixA_0_all} in 
Section~\ref{sec:main_proof}. Note that the latter does not inmediately follow from the former since in general
\[
 \sum_{i=1}^n\mu^2(A_0,\lambda^{(i)},v^{(i)})>>\mum^2(A_0).
\]

\subsection{Initial triples and global algorithms}

The \EC\ routine assumes an initial triple $(A_0,\lambda_0,v_0)$ 
at hand. We next deal with this issue. We first consider the 
case of computing a single eigenpair. In this case we consider the
diagonal matrix $H$ whose diagonal entries are
$(1,0,\ldots,0)$ and its well-posed eigenpair $(1,e_1)$. 
\bigskip\bigskip

\algoritmo
\begin{algorithm}\label{alg:LV}
{\sf Single\_Eigenpair}\\
\inputalg{$A\in \Cnn$}\\
\bodyalg{
$(\zeta,w):=\mbox{\EC}(A,H,1,e_1)$\\[2pt]
}
\Output{$(\zeta,w)\in \C\times \C^{n}$}\\
\postcond{The execution halts if the lifting of 
$\mathcal{L}_{H,A}$ at $(1,e_1)$ doesn't cut $\Sigma'$. 
In this case, the returned $(\zeta,w)$ is an approximate 
eigenpair of~$A$.}
\end{algorithm}
\falgoritmo

We can formally state (and prove) the first of our main results. To 
this end, we define the average cost $\avcost(n)$ of 
{\sf Single\_Eigenpair} to be the average (over the input matrix $A$) 
of the number of arithmetic operations performed by the algorithm. 
We similarly define its smoothed cost $\smcost(n,\sigma)$. 

\begin{theorem}\label{thm:main}
Algorithm~{\sf Single\_Eigenpair} returns (almost surely) an 
approximate eigenpair of its input $A\in\Cnn$. Its average 
cost satisfies
$$
    \avcost(n)=\Oh(n^7). 
$$ 
For every $0<\sigma\leq 1$, its smoothed cost satisfies
$$
    \smcost(n,\sigma)=\Oh\Big(\frac{n^7}{\sigma^2}\Big). 
$$ 
\end{theorem} 

\proof
Lemma~~\ref{lem:char} and the fact that $\|H\|_F=1$ imply that 
$\mu(H,1,e_1)=1$. The statement is then a consequence of 
Proposition~\ref{prop:fixA_0} and the fact that the average cost 
is obtained by multiplying $\aviter(H,1,e_1)$ by the cost $\Oh(n^3)$ 
of each iteration.  
\eproof

\begin{remark}
The triple $(H,1,e_1)$ is the version, in our context, of the 
initial pair proposed by Shub and Smale in~\cite{Bez5} for the 
computation of zeros of polynomial systems. In this later context, 
the problem of showing that one can efficiently follow linear homotopies with this initial pair 
remains open.
\end{remark}

The fact that any other eigenpair of $H$ is ill-posed prevent us from 
using them to compute other eigenpairs of $A$. If we 
want to compute all the eigenpairs of $A$ we will need to consider 
a different approach. 

To do so, for any fixed $n\geq2$, let
\[
D=\diag(\eta_1,\ldots,\eta_n),
\]
where the $\eta_i$ are points in the unit side hexagonal 
lattice chosen in such an order that $0=|\eta_1|\leq\cdots\leq|\eta_n|$. 

\begin{lemma}\label{lem:initial}
We have $\mum(D)\leq \frac{\sqrt{3/2}}{\pi}n+o(n)$.
\end{lemma}

\proof
The hexagonal lattice is the set of points of the form
\[
Q\binom{a}{b},\quad \text{where}\quad Q=\begin{pmatrix}1&1/2\\0&\sqrt{3}/2\end{pmatrix},\quad a,b\in\Z.
\]
We first find a real number $r>0$ with the property that the circle $\mathbb{D}(r)$ of radius $r$ contains at least $n$ such lattice points. To do so, note that a lattice point $Q\binom{a}{b}$ is in $\mathbb{D}(r)$ if and only if
\[
\binom{a}{b}\in Q^{-1}\mathbb{D}(r).
\]
Now, the singular values of $Q^{-1}$ are $\sqrt{2}$ and $\sqrt{2/3}$, so $Q^{-1}\mathbb{D}(r)$ is an ellipse of area $2\pi r^2/\sqrt{3}$ and maximal radius $\sqrt{2}r$. Dividing by the smallest integer $N$ that is greater than  $2\sqrt{2}r$ and translating the resulting ellipse to have center $(1/2,1/2)$, we look for points of the form $(a/N,b/N)$ with $a,b\in\{0,\ldots,N-1\}$ which lie inside an ellipse of area
\[
 \frac{2\pi r^2}{\sqrt{3}N^2}
\]
contained in $[0,1]^2$. This is a particular instance of the problem of counting lattice points in semialgebraic sets, a well studied problem for which a quite complete solution is for example \cite[Th. 3, p. 327]{BlCuShSm98}. We conclude that the number of points in the hexagonal lattice in $\mathbb{D}(r)$ is at least
\[
\frac{2\pi r^2}{\sqrt{3}}-2N\geq \frac{2\pi r^2}{\sqrt{3}}-4\sqrt{2}r-2.
\]
In particular, we can find $n$ lattice points in a circle of radius $r=3^{1/4}/(2\pi)^{1/2}n^{1/2}+o(n^{1/2})$. Moreover, around each lattice point we can place a circle of radius $1/2$ without overlappings. It is a simple exercise to check that for any $z\in\C$ and $\e>0$,
\[
|z|^2\leq\frac{1}{\pi\e^2}\int_{|y-z|<\e}|y|^2\,dy.
\]
We thus have
\[
\|D\|_F^2=\sum_{j=1}^n|\eta_j|^2\leq \sum_{z}|z|^2\leq\frac{4}{\pi}\sum_{z} \int_{|y-z|<1/2}|y|^2\,dy\leq \frac{4}{\pi}\int_{|y|<r+1/2}|y|^2\,dy,
\]
where $z$ runs over all lattice points contained in $\mathbb{D}(r)$. Solving the last integral yields
\[
\|D\|_F^2\leq 2r^4+o(r^4)=\frac{3n^2}{2\pi^2}+o(n^2).
\]
Finally, from Lemma~\ref{lem:char} we conclude that
\[
\mum^2(D)=\|D\|_F^2\leq \frac{3n^2}{2\pi^2}+o(n^2).
\]
\eproof

We now put together the continuation algorithm \EC\ and this 
specific initial triple.
\bigskip\bigskip

\algoritmo
\begin{algorithm}\label{alg:All}
{\sf All\_Eigenpairs}\\
\inputalg{$A\in \Cnn$}\\ 
\bodyalg{
compute $D$\\[2pt]
for $j\in\{1,\ldots,n\}$ do\\[2pt]
\espacio $(\zeta_j,w_j):=\mbox{\EC}(A,D,\eta_j,e_j)$\\[2pt]
}
\Output{$\{(\zeta_1,w_1),\ldots,(\zeta_n,w_n)\}\in (\C\times \C^{n})^n$}\\
\postcond{The algorithm halts if 
$\mathcal{L}_{D,A}\cap\Sigma=\emptyset$. 
In this case, the pairs $(\zeta_j,w_j)$ are approximate 
eigenpairs of $A$ with pairwise different associated eigenpairs.}
\end{algorithm}
\falgoritmo

We can now state (and prove) the second 
of our main results.

\begin{theorem}\label{thm:main2}
Algorithm~{\sf All\_Eigenpairs} returns (almost surely)
$n$~approximate eigenpairs of its input $A\in\Cnn$, 
with pairwise different associate eigenpairs. Its average 
cost satisfies
$$
    \avcost(n)=\Oh(n^9). 
$$ 
For every $\sigma\leq 1$ its smoothed cost satisfies
$$
    \smcost(n,\sigma)=\Oh\Big(\frac{n^9}{\sigma^2}\Big). 
$$ 
\end{theorem} 

\proof
It easily follows from Lemma~\ref{lem:initial} and 
Proposition~\ref{prop:fixA_0_all}.  
\eproof

\subsection{Randomized algorithms}

In this section we follow the ideas in~\cite{BePa:11} adapting them to
the case of eigenvalue, eigenvector computations. 
Consider the probability distribution $\mcD$ in the solution variety 
$\V$ defined via the following procedure:

\begin{equation}\label{eq:D} 
\begin{split}
&\mbox{\tt randomly choose  
$A_0\sim\mcN_{\C^{n\times n}}$}\\
&\mbox{\tt randomly choose one eigenpair $(\lambda_0,v_0)$ 
of $A_0$}
\end{split}
\end{equation}

\noindent
Next assume that we have a routine {\sf draw\_from\_{$\mcD$}} 
to draw triples $(A_0,\lambda_0,v_0)$ 
from the distribution $\mcD$ on $\V$. Then we can consider the 
following algorithm.

\bigskip\bigskip

\algoritmo
\begin{algorithm}\label{alg:Random}
{\sf Random\_initial\_triple (scheme)}\\
\inputalg{$A\in \Cnn$}\\
\bodyalg{
$(A_0,\lambda_0,v_0):={\mathsf{draw\_from\_}{\mcD}}$\\[2pt] 
 $(\zeta,w):=\mbox{\EC}(A,A_0,\lambda_0,v_0)$\\[2pt]
}
\Output{$(\zeta,w)\in (\C\times \C^{n})^n$}\\
\postcond{The execution halts if the lifting of 
$\mathcal{L}_{A_0,A}$ at $(\lambda_0,v_0)$ doesn't cut $\Sigma'$.
In this case, $(\zeta,w)$ is an approximate 
eigenpair of $A$.}
\end{algorithm}
\falgoritmo

The interest of this algorithmic scheme is that we can prove 
better bounds (than those in Theorem~\ref{thm:main}) for the number 
of iterations of \EC.
\smallskip
 
\begin{theorem}\label{th:randomhomotopy}
The expected average number of homotopy steps of 
\EC\ in  Algorithm~\ref{alg:Random} satisfies
\begin{eqnarray*}
& & \Exp_{A\sim\mcN_{\C^{n\times n}}}
 \Exp_{(A_0,\lambda_0,v_0)\sim\mcD}
  K(A,A_0,\lambda_0,v_0)
\leq 4Cn^2,
\end{eqnarray*}
where $C$ is as in Theorem~\ref{thm:main_path_following}.
\end{theorem}

Theorem~\ref{th:randomhomotopy}, 
which will be proved in Section~\ref{sec:proofofrandomisrandom}, 
brings to focus the need 
of an implementation of {\sf draw\_from\_{$\mcD$}}. 
It must be noted though that a direct implementation is not 
possible since the second line in~\eqref{eq:D} 
(choosing $(\lambda_0,v_0)$ at random) implicitly requires solving
an EVP problem, the very question that this article is attempting to solve! 
This is a similar situation to that solved 
in~\cite{BePa:09,BePa:11}, where a random polynomial system and
one of its zeros at random had to be chosen. It is also similar to that 
dealt with in~\cite{ArmCuc} for the computation of eigenpairs of Hermitian 
matrices. A version of the proof of 
Theorem~\ref{th:randomhomotopy} with the Gaussian Unitary Ensemble 
replacing the Gaussian distribution, actually yields the following 
improvement over the main result in~\cite{ArmCuc}.

\begin{corollary}
In the case of Hermitian matrices, the expected average number 
of homotopy steps of \EC\ in Algorithm~\ref{alg:Random} (with 
the randomization algorithm in~\cite{ArmCuc} in place of 
${\mathsf{draw\_from\_}{\mcD}}$) is $\Oh(n^2)$. This yields an  
expected cost of $\Oh(n^5)$ for the computation of one eigenpair
and of $\Oh(n^6)$ for the computation of all of them. 
Here the input matrix $H$ is drawn from the Gaussian Unitary Ensemble.
\eproof
\end{corollary}

Following the ideas in those papers, we note that  
Theorem~\ref{th:randomhomotopy} would yield
an algorithm with average running time $\Oh(n^5)$ if we could find some
collection of probability spaces $\Omega_n$ and functions
$\psi_n:\Omega_n\rightarrow\V_n$, $n\geq2$, such that:
\begin{enumerate}
\item Choosing $\omega\in\Omega_n$ can be done starting with a number 
of random choices of numbers with the $N_\C$ distribution, and 
performing some arithmetic operations on the results, the total expected 
running time being at most $\Oh(n^5)$.
\item Given $\omega\in\Omega_n$, $\psi_n(\omega)$ is computable 
in average time $\Oh(n^5)$, that is, the expected number of arithmetic 
operations for computing $\psi_n(\omega)$ must be $\Oh(n^5)$.
\item Choosing $\omega$ at random in $\Omega_n$ and computing 
$({A}_0,{\lambda}_0,v_0)=\psi_n(\omega)$ is equivalent to 
choosing $A_0\sim\mcN_{\C^{n\times n}}$ at random and 
choosing at random $(\lambda_0,v_0)$ such that $Av_0=\lambda_0v_0$. 
That is, for any measurable mapping $\phi:\V_n\rightarrow[0,\infty]$ 
we must have
\begin{equation}\label{eq:wish}
  \Exp_{\omega\sim\Omega_n}\left(\phi(\psi_n(\omega))\right)
 =\Exp_{A_0\sim\mcN_{\C^{n\times n}}}
  \left(\frac{1}{n}\sum_{\lambda_0,v_0:A_0v_0=\lambda_0v_0}
  \phi(A_0,\lambda_0,v_0)\right),
\end{equation}
so that we can apply this equality to 
\[
    \phi(A_0,\lambda_0,v_0)
    =\Exp_{A\sim\mcN_{\C^{n\times n}}}
    \left(K(A,A_0,\lambda_0,v_0)\right)
\]
and apply Theorem~\ref{th:randomhomotopy}.
\end{enumerate}
Unfortunately, we are not able to produce a collection of probability
spaces $\Omega_n$ and functions $\psi_n$ as described
above. However, we will prove that relaxing~\eqref{eq:wish} to the
following less restrictive situation is actually possible: instead of
demanding the equality in~\eqref{eq:wish} we can just demand an
inequality where the right-hand term is multiplied by some polynomial
in $n$. Moreover, we do not need~\eqref{eq:wish} to hold for every
measurable function $\phi$ since all the interesting functions for the
EVP problem are projective functions, invariant under the action of
the unitary group. We can thus relax~\eqref{eq:wish} to hold only with
a polynomially bounded inequality, and for unitary invariant
projective functions. Proving that this can actually be done is our
goal now. 

We start by defining $\Omega_n$ and $\psi_n$. Consider
the classical Stiefel manifold consisting of orthonormal $(n-1)$-frames 
in $\C^n$, given by
\[
   \mathcal{S}_{n-1}(\C^n)=\{Q\in\C^{n\times(n-1)}:Q^*Q=I_{n-1}\},
\]
endowed with its natural probability measure given by the restriction 
of the Frobenius Hermitian structure to the tangent bundle.

For every $n\geq2$, let 
\[
     \A_n:=\{(M,Q):\ker(M)=\ker(Q^*)\}
     \subseteq\C^{(n-1)\times n}\times \mathcal{S}_{n-1}(\C^n).
\]
In other words, $\A_n$ consists of pairs of matrices $M,Q$ such that the 
columns of $Q$ form an orthogonal basis of the complement of $\ker(M)$. 
The set $\A_n$ has a natural probability measure $\mu_{\A_n}$ given by
\[
  \mu_{\A_n}(X):=\Exp_{M\sim\mathcal{N}_{\C^{(n-1)\times n}}}
  \left(\frac{1}{\Vol(Q:(M,Q)\in\A_n)}
  \int_{Q:(M,Q)\in\A_n}\uno_X(M,Q)\,dQ\right),
\]
for measurable sets $X\subseteq\A_n$. 
\begin{definition}\label{def:OmegayVarphi}
Let
\[
  \Omega_n:=\{(\lambda,w,(M,Q)):2
  \Re(\bar{\lambda}\mathrm{tr}(MQ))\leq 1-|\lambda|^2(n-1)\}\subseteq \C\times\C^{n-1}\times\A_n
\]
be endowed with the product measure $\mu_{\Omega_n}$, normalized 
to have total unit mass (see Remark \ref{remark:details} below).
Then, let
\[
  \psi_n(\lambda,w,M,Q)=
  \left(\begin{pmatrix}
 \lambda&w^*\\0&MQ+\lambda I_{n-1}
\end{pmatrix},\lambda,e_1
\right).
\]
\end{definition}
\begin{remark}\label{remark:details}
An explicit description of the product measure $\mu_{\Omega_n}$ is as follows:
\[
   \mu_{\Omega_n}(Y):=C_n\Exp_{\lambda,w}
   \left(\mu_{\A_n}(\{(M,Q):(\lambda,w,M,Q)\in Y\})\right),
\]
for measurable sets $Y\subseteq\Omega_n$, where 
$\lambda\sim\mathcal{N}_{\C}$, 
$w\sim\mathcal{N}_{\C^{n\times 1}}$ and $C_n$ is 
a normalizing constant given by
\begin{equation}\label{eq:Cn}
 C_n^{-1}=\Prob_{\lambda,M,Q}
 \left(2\Re(\bar{\lambda}\mathrm{tr}(MQ))\leq 1-|\lambda|^2(n-1)\right).
\end{equation}
\end{remark}
Our last main result (see Section~\ref{sec:proofmainrandom} for a proof)
is the following.

\begin{theorem}\label{th:mainrandom}
Let $\Omega_n$ and $\psi_n$ be as in 
Definition~\ref{def:OmegayVarphi} for all $n\geq2$. Then:
\begin{enumerate}
\item 
Choosing $\omega\in \Omega_n$ can be done by choosing 
$2n^2-2n+1$ numbers with the $N_\C$ distribution, checking a test
which involves the computation of a Moore-Penrose inverse and
computing two QR decompositions.  This process
may be repeated as a function of the test's outcome, 
but the expectation of the number of times the test is performed
is at most $C_n$. The total expected running time is
$\Oh(n^3C_n)$.
\item 
Given $\omega\in\Omega_n$, computing $\psi_n(\omega)$ 
can be done with running time $\Oh(n^3)$. 
\item 
For any unitarily invariant measurable mapping
$\phi:\V\to[0,\infty]$ we have:
\begin{equation}\label{eq:wish2}
  \Exp_{\omega\sim\Omega_n}\left(\phi(\psi_n(\omega))\right)\leq 
 e\,nC_n \Exp_{A_0\sim\mcN_{\C^{n\times n}}}
  \left(\frac{1}{n}\sum_{\lambda_0,v_0:A_0v_0
  =\lambda_0v_0}\phi(A_0,\lambda_0,v_0)\right).
\end{equation}
\end{enumerate}
\end{theorem}
Note that~\eqref{eq:wish2} can be understood as follows: let $m_1$ be
the push-forward measure of $\psi_n$ in $\V$ and let $m_2$ be the
measure in $\V$ given by
\[
  m_2(X)=\E_{A_0\sim\mcN_{\C^{n\times n}}}
 \left(\frac{1}{n}\sharp\big\{(\lambda_0,v_0):A_0v_0
 =\lambda_0v_0,(A_0,\lambda_0,v_0)\in X\big\}\right),
\]
for any measurable set $X\subseteq\V$. Then, the Radon-Nikodim 
derivative $dm_1/dm_2$ is bounded above by $enC_n$.

\begin{problem}
Describe an alternative collection $(\Omega_n,\psi_n)$ which
satisfies a sharper version of~\eqref{eq:wish2}, with a constant in
the place of $ n C_n$. This would improve the running time of 
Algorithm~{\sf Random\_initial\_triple} below by a factor of $\Oh(n C_n)$.
\end{problem}
We are now prepared to describe our random homotopy algorithm.
\bigskip\bigskip

\algoritmo
\begin{algorithm}\label{alg:Random2}
{\sf Random\_initial\_triple}\\
\inputalg{$A\in \Cnn$}\\
\bodyalg{
Randomly choose  $\omega\in \Omega_n$\\[2pt]
$(A_0,\lambda_0,v_0):=\psi_n(\omega)$ (note that $v_0=e_1$)\\[2pt]
 $(\zeta,w):=\mbox{\EC}(A,A_0,\lambda_0,v_0)$\\[2pt]
}
\Output{$(\zeta,w)\in (\C\times \C^{n})^n$}\\
\postcond{The execution halts if the lifting of 
$\mathcal{L}_{A_0,A}$ at $(\lambda_0,v_0)$ doesn't cut $\Sigma'$. 
In this case, $(\zeta,w)$ is an approximate 
eigenpair of $A$.}
\end{algorithm}
\falgoritmo

From Theorem~\ref{th:mainrandom}, the expected running time of the
computation of $(A_0,\lambda_0,v_0)$ is $\Oh(n^3C_n)$. 
Moreover, the expected number of homotopy steps in the execution
of  $\mbox{\EC}(A,A_0,\lambda_0,v_0)$ is
\[
  S=\Exp_{A\sim\mcN_{\C^{n\times n}},\,\omega\sim\Omega_n}
  \left(K(A,\psi_n(\omega))\right)
  =\Exp_{\omega\sim\Omega_n}
  \left(\Exp_{A\sim\mcN_{\C^{n\times n}}}
  \left(K(A,\psi_n(\omega))\right)\right),
\]
where $K(A,\psi_n(\omega))$ is as in 
Theorem~\ref{thm:main_path_following}. From~\eqref{eq:wish2}, we have
\begin{eqnarray*}
  S&\leq& e\, n C_n\Exp_{A_0\sim\mcN_{\C^{n\times n}}}
  \left(\frac{1}{n}\sum_{\lambda_0,v_0:A_0v_0=\lambda_0v_0}
  \Exp_{A\sim\mcN_{\C^{n\times n}}}
  \left(K(A,A_0,\lambda_0,v_0)\right)\right)\\
  &\underset{\text{\rm Th.~\ref{th:randomhomotopy}}}{\leq}& \Oh(n^3C_n).
\end{eqnarray*}
We multiply the number of steps by $\Oh(n^3)$ to get the following 
complexity bound.

\begin{theorem}\label{th:random}
Algorithm~{\sf Random\_initial\_triple} returns (almost surely) an 
approximate eigenpair of its input $A\in\Cnn$. Its average 
cost satisfies
\begin{equation}\tag*{\qed}
    \avcost(n)=\Oh(n^6C_n)\underset{\text{\rm Lemma~\ref{lem:SnCn}}}{\leq}
    \Oh(n^7). 
\end{equation}
\end{theorem}

\section{Some properties of the condition number $\mu$}
\label{sec:condition_property}

There is a general geometric framework for defining condition numbers,
see~\cite[\S14.3]{Condition}. In our situation, 
this framework takes the following form. 

If $(A,\la,v)\in\W$, then from the inverse function theorem the projection 
$\pi_1\colon\V\to\Cnn$ (cf.~(3)), around $(A,\la,v)$, 
has a local inverse $\mathcal{U}\to\V, A\mapsto (A,G(A))$, that is defined 
on an open neighborhood $\mathcal{U}$ of $A$ in $\Cnn$.
We call $G$ the {\em solution map}. The map $G$ decomposes as 
$G=(G_{\la},G_v)$, where 
$$
  G_{\la}:\mathcal{U}\to\C
  \qquad{\mbox{and}}\qquad
  G_v:\mathcal{U}\to\P(\C^n)
$$
associate to matrices $B\in\mathcal{U}$ an eigenvalue and 
the corresponding eigenvector. Let 
$$
  DG_{\la}(A):\Cnn\to\C
  \qquad{\mbox{and}}\qquad
  DG_v(A):\Cnn\to v^\perp
$$
be the derivatives of these maps at $A$. The condition numbers for the 
eigenvalue $\la$ and the eigenvector $v$ of $A$ are defined as follows:
$$
  \mu_{\la}(A,\la,v):=\|DG_{\la}(A)\|
  \qquad{\mbox{and}}\qquad
  \mu_v(A,\la,v):=\|DG_v(A)\|,
$$
where the norms are the operator norms with respect to the  
chosen norms (on $\Cnn$ we use the Frobenius norm and on $v^\perp$ 
the norm given by~\eqref{eq:defTip}). 

The following result, Lemma~14.17 in~\cite{Condition}, 
gives explicit descriptions of $DG_{\la}$ and 
$DG_v$. Before stating it, we recall that if $\la$ is an 
eigenvalue of $A$ there exists 
$u\in\P(\C^n)$ (the {\em left eigenvector}) such that $(A-{\la}\Id)^*u=0$. 
Recall the linear map 
$A_{\lambda,v}\colon v^\perp\to v^\perp$ 
introduced in~\eqref{eq:defAlv}.

\begin{lemma}\label{lem:lefteigen}
Assume that $Av=\lambda v$ and $\lambda$ has multiplicity $1$. 
Then, the associated left eigenvector is
\begin{equation}\label{eq:lefteigen}
u=v-i_{\C^n}A_{\lambda,v}^{-*}P_{v^\perp}A^*v.
\end{equation}
Here we denoted $A_{\lambda,v}^{-*}:=\big(A_{\lambda,v}^{-1}\big)^*$. Note that $\langle u,v\rangle=\|v\|^2$.
\end{lemma}

\proof
Take a representative such that $\|v\|=1$ and let
\[
z:=i_{\C^n}A_{\lambda,v}^{-*}P_{v^\perp}A^*v\underset{\eqref{eq:defAlvbisbis}}{=}\big((\Id-vv^*)(A-\lambda \Id)^*(\Id-vv^*)\big)^\dagger A^*v.
\]
From the definition of the Moore--Penrose pseudoinverse, $z$ is
the unique element in $v^\perp$ that minimizes 
$\|(\Id-vv^*)(A-\lambda \Id)^*z-A^*v\|$, that is we have 
$(\Id-vv^*)(A-\lambda \Id)^*z=P_{v^\perp}(A^*v)=(\Id-vv^*)A^*v$ 
or equivalently
\[
  (A-\lambda \Id)^*z=A^*v+tv\quad\text{ for some }t\in\C.
\]
Premultiplying both sides by $v^*$ we have
\[
  v^*(A-\lambda \Id)^*z=v^*A^*v+t\|v\|^2\quad\then\quad 0=(\lambda v)^*v+t\|v\|^2=(\bar\lambda+t)\|v\|^2,
\]
so we have $t=-\bar\lambda$ and then
\[
(A-\lambda \Id)^*(v-z)
 =(A-\lambda \Id)^* v-A^*v+\bar\lambda v=0,
\]
that is, $v-z$ is a left eigenvector of $A$ with associated (left) eigenvalue $\bar\lambda$ as wanted.
\eproof

The following is Lemma~14.17 in~\cite{Condition}.

\begin{lemma}\label{lem:14.17}
Let $(A,\la,v)\in\W$ and let $u$ be a left eigenvector of $A$ 
with eigenvalue~$\bar{\la}$. Then:

{\bf (a)\ } We have $\langle v,u\rangle \ne 0$. 
\smallskip

{\bf (b)\ } The derivative of the solution 
map is given by 
$DG(A)(\dot{A}) = (\dot{\lambda},\dot{v})$, where 
\begin{equation}\tag*{\qed}
 \dot{\lambda} = \frac{\langle \dot{A}v,u\rangle}
  {\langle v,u\rangle}, \quad 
 \dot{v} = -A_{\lambda,v}^{-1} \, P_{v^\perp} \dot{A} v.
\end{equation}
\end{lemma}

The following result, which follows directly from 
Lemma~\ref{lem:14.17}, was already pointed out in~\cite{vanLoan} (see also \cite[Prop. 14.15]{Condition}).

\begin{proposition}\label{pro:cn-eigenv}
Choosing the Frobenius norm on $T_A\Cnn = \C^{n\times n}$ and 
$\frac1{\|v\|}\, \|\cdot\|$ on $v^\perp$, the condition numbers $\mu_v$ 
for the eigenvector problem and $\mu_{\la}$ for the eigenvalue 
problem satisfy: 
\begin{equation*}
\mu_{\la}(A,\lambda,v)\;=\;\|DG_{\la}(A)\| 
 = \frac{\|u\|\|v\|}{|\langle u ,v \rangle | }
 \underset{\eqref{eq:lefteigen}}{=}\frac{\|u\|}{\|v\|}\leq\sqrt{1+\mu^2(A,\lambda,v)}
\end{equation*}
and
\begin{equation}\tag*{\qed}
\mu_{v}(A,\la,v)\;=\;\|DG_v(A)\| 
 = \big\| A_{\lambda,v}^{-1} \big\|=\frac{\mu(A,\lambda,v)}{\|A\|_F}. 
\end{equation}
\end{proposition}
\medskip

\proofof{Proposition~\ref{prop:dotzdotf-spherical}}
The first two inequalities are immediate from 
Proposition~\ref{pro:cn-eigenv}. For the third one,  note that
\begin{equation*}
\|\dot{\Gamma}(t)\|\;=\;\|(\dot{A},\dot{\la},\dot{v})\| 
 \;\le\;\|\dot{A}\|\sqrt{1+\mu_t^2+(1+\mu_t^2)} 
  \;\leq\; \|\dot{A}\|\sqrt{6\mu_t^2} 
\end{equation*}
the last inequality since $\mu_t\geq\frac{1}{\sqrt{2}}$ 
(Lemma~\ref{le:lb_mu}).
\eproof

Our last lemma is a version of Lemma \ref{le:lb_mu} without the assumption that our point lies on $\W$.
\begin{lemma}\label{le:lb_mu2}
For $A\in\Cnn$, $w\in\C^n$ and $\zeta\in\C$ with $|\zeta|\leq\|A\|_F$ we 
have 
\[
 \mu(A,\zeta,w) \ge\frac{1}{1+\sqrt{1-\frac{\|w^*A\|^2}{\|w\|^2\|A\|_F^2}}}\geq \frac{1}{2}
\]
\end{lemma}
\proof
W.l.o.g we can assume that $w=e_1$ and write
\[
 A=\begin{pmatrix}\lambda&a^*\\b&\hat A\end{pmatrix},
\]
where $\lambda\in\C$ and $a,b\in\C^{n-1}$. Then, $A_{\zeta,w}\equiv \hat A-\zeta\Id_{n-1}$ and
\begin{align*}
 \frac{\mu(A,\zeta,e_1)}{\|A\|_F}=\|(\hat A-\zeta\Id_{n-1})^{-1}\|\geq\frac{1}{\|\hat A-\zeta\Id_{n-1}\|}\geq\;& \frac{1}{\|\hat A\|+|\zeta|}\\\geq\;&\frac{1}{\sqrt{\|A\|_F^2-\|e_1^*A\|^2}+\|A\|_F}.
\end{align*}
The lemma follows.

\eproof
\section{Proofs of Theorem~\ref{th:lipsch} and Corollary~\ref{cor:b}}
\label{se:lipschitz} 

It will be handy to use the definition of $\mu$ given 
in~\eqref{eq:mualternative}. We start with a very simple linear algebra 
lemma about the Moore-Penrose pseudoinverse.

\begin{lemma}\label{lem:moore}
Let $R,R'\in\C^{n\times n}$ be such that $R$ has rank $n-1$. Assume 
moreover that $\det(R')=0$ and
\[
  \|R-R'\|\leq \frac{\e}{\|R^\dagger\|},\quad \text{ for some } 0\leq\e<1.
\]
Then, $R'$ has rank $n-1$ and
\[
  \frac{\|R^\dagger\|}{1+\e}\leq \|R'^\dagger\|
  \leq \frac{\|R^\dagger\|}{1-\e}.
\]
\end{lemma}

\proof
Let $\sigma$ and $\sigma'$ be the $(n-1)$th singular values of $R$ 
and $R'$, respectively. 
Note that $\sigma-\|R-R'\|\leq \sigma'\leq \sigma+\|R-R'\|$ 
(this is a classical fact proved for the first time in \cite{Weyl1912}, see also ~\cite[Cor. 8.6.2]{GoVa96}). In particular, $\sigma'\geq \sigma-\|R-R'\|=\|R^\dagger\|^{-1}-\|R-R'\|>0$, so $R'$ 
has rank at least $n-1$ and by hypothesis it has rank $n-1$. 
Moreover, we have
\[
   \|R^\dagger\|=\frac{1}{\sigma}
   =\frac{1}{\sigma'}\frac{\sigma'}{\sigma}
   \leq\|R'^\dagger\|\frac{\sigma+\|R-R'\|}{\sigma}
   \leq(1+\e)\|R'^\dagger\|.
\]
The upper bound follows from a similar argument.
\eproof
\medskip

\proofof{Theorem~\ref{th:lipsch}} 
Choose representatives such that $\|v\|=\|v'\|=1$ and let
\[
  Q:=(\Id-vv^*)(A-\lambda\Id)(\Id-vv^*),\quad 
  Q':=(\Id-v'v'^*)(A'-\lambda'\Id)(\Id-v'v'^*).
\]
We have $\rank(Q)=n-1$ since $\mu(A,\lambda,v)<\infty$ by 
our assumption, cf.~\eqref{eq:mualternative}. 
We claim that 
\begin{equation}\label{eq:qqprime}
\|Q-Q'\|\leq \frac{\e}{\|Q^\dagger\|}, 
\end{equation}
which from 
Lemma~\ref{lem:moore} implies that $Q'$ has rank $n-1$ and
\[
 \frac{\|Q^\dagger\|}{1+\e}\leq \|Q'^\dagger\|
  \leq \frac{\|Q^\dagger\|}{1-\e},
\]
that is (recall $\|A\|_F=\|A'\|_F=1$),
\[
   \frac{\mu(A,\lambda,v)}{1+\e}\leq \mu(A',\lambda',v')
   \leq \frac{\mu(A,\lambda,v)}{1-\e},
\]
as wanted. 

It remains to prove the claim. To do so, 
let $(A_t,\lambda_t,v_t)$ be a geodesic in $\Cnn\times\C\times\P(\C^n)$, parametrized by arc-length, 
joining $(A,\lambda,v)$ and $(A',\lambda',v')$, so by hypothesis we 
have $t\in[0,\e/(4\sqrt{3}\|Q^\dagger\|)]$ and we can choose representatives $v_t$ in such a way that
\[
\dot v_t\perp v_t,,\quad\|v_t\|=1, \quad \|\dot A_t\|_F^2+|\dot\lambda_t|^2+\|\dot v_t\|^2=1,\quad\text{ for all } t\in[0 ,\e/(4\sqrt{3}\|Q^\dagger\|)].
\]
Let 
$Q_t:=(\Id-v_tv_t^*)(A_t-\lambda_t\Id)(\Id-v_tv_t^*)$. In order to bound $\|\dot Q_t\|$ we first note that for $x\in\C^n$,
\begin{align*}
 \|(\dot v_t v_t^*+v_t\dot v_t^*)(x)\|=&\|\dot v_t\langle x,v_t\rangle+v_t\langle x,\dot v_t\rangle\| =\sqrt{\|\dot v_t\|^2|\langle x,v_t\rangle|^2+\|v_t\|^2|\langle x,\dot v_t\rangle|^2}\\
 =&\|\dot v_t\|\,\sqrt{\left|\langle x,v_t\rangle\right|^2+\left|\langle x,\frac{\dot v_t}{\|\dot v_t\|}\rangle\right|^2}\leq \|\dot v_t\|\,\|x\|,
\end{align*}
that is $\|\dot v_t v_t^*+v_t\dot v_t^*\|\leq \|\dot v_t\|$. On the other hand, 
\[
\|A_t-\lambda_t\Id\|\leq \|A_t\|+|\lambda_t|\leq \|A_t\|+\|A_t\|\leq 2\|A_t\|_F= 2, 
\]
so computing the derivative of $Q_t$ we see that
\begin{align*}
   \|\dot Q_t\|\;\leq\;&2\|\dot v_t v_t^*+v_t\dot v_t^*\|\|A_t-\lambda_t\Id\|+\|\dot A_t-\dot \lambda_t\Id\|\\\leq\;& 4\|\dot v_t\|+\|\dot A_t-\dot\lambda_t\Id\|
   \leq 4(\|\dot v_t\|+\|\dot A_t\|+|\dot\lambda_t|)\\
   \leq\; & 4\sqrt{3}(\|\dot v_t\|^2+\|\dot A_t\|^2+|\dot\lambda_t|^2)^{1/2}
   =4\sqrt{3}.
\end{align*}
Thus, we have 
$\|Q-Q'\|\leq \int_0^{\e/(4\sqrt{3}\|Q^\dagger\|)}4\sqrt{3}\,dt\leq 4\sqrt{3}\e/(4\sqrt{3}\|Q^\dagger\|)$, finishing the proof of \eqref{eq:qqprime} and that of Theorem \ref{th:lipsch}.
\eproof
\medskip

\proofof{Corollary~\ref{cor:b}} 
Consider the portion $\mathcal{L}_{A,A'}=\{A_t\}$ of great circle 
defined for $t\in[0,\alpha]$ where $\alpha=d_\SS(A,A')$. 
Let $(\lambda,v)$ be any eigenpair of $A$. From the inverse
function theorem, the map $t\mapsto \Gamma(t)=(A_t,\lambda_t,v_t)\in\W$ (given by the local inverse of $\pi_1$) is
defined and, denoting $\dist(t):=\dist((A,\lambda,v),(A_t,\lambda_t,v_t))$, satisfies
\begin{equation}\label{eq:igualdad}
  \dist(t)<\frac{\e}{4\sqrt{3}\mu(A,\lambda,v)},\quad t\in[0,\delta)
\end{equation}
for some maximal $\delta\in(0,\alpha]$. 
From Theorem~\ref{th:lipsch} this implies
\begin{equation}\label{eq:olv}
  \mu(\Gamma(t))\leq 
  (1+\e)\mu(A,\lambda,v),\quad t\in[0,\delta),
\end{equation}
giving a global bound for the norm of the derivative of the solution map. 
Thus, $\Gamma(t)$ can be continuously extended to $[0,\delta]$ 
and from maximality of $\delta$, either $\delta=\alpha$ or $\delta<\alpha$ 
and~\eqref{eq:igualdad} 
must become an equality when changing $t$ to $\delta$. In this last case, which we will rule out by contradiction, note that from Proposition~\ref{prop:dotzdotf-spherical} 
and~\eqref{eq:olv}
\[
  \dist(\delta) \leq \int_0^\delta\|\dot \Gamma(t)\|\,dt
  \leq\sqrt{6}\int_0^\delta\mu(\Gamma(t))\,dt
  \leq \sqrt{6}(1+\e)\delta\mu(A,\lambda,v).
\]
Now, in our range of values we have
$\alpha=d_\SS(A,A')=2\arcsin(\|A-A'\|_F/2)<
\frac{1001}{1000}\|A-A'\|_F$ which, together with 
$\delta<\alpha$, implies
\[
   \dist(\delta)<\frac{1001\|A-A'\|_F\sqrt{6}(1+\e)\mu(A,\lambda,v)}
  {1000}\leq\frac{1001\e\sqrt{6}(1+\e)}{50000\mu(A,\lambda,v)}
 < \frac{\e}{4\sqrt{3}\mu(A,\lambda,v)},
\]
which is a contradiction (note that we have used the bound  $\|A-A'\|_F\, \mu^2(A,\la,v) \le \e/50$). We thus conclude that
$\delta=\alpha$ and the corollary follows.
\eproof

\section{Condition-length homotopy continuation}\label{sec:conditionlength}

The goal of this section is to fully describe the routine 
{\sf Choose\_step} (with which algorithm \EC\ will be complete) 
and to prove Theorem~\ref{thm:main_path_following}. 

To do so, it will be useful to denote by $\beta(A,\zeta,w)$ the
length (in the tangent space) of the Newton step with matrix
$A\in\SS$ and input $(\zeta, w)\in\C\times\P(\C^n)$. 
That is, if we take a representative such that $\|w\|=1$,
\[
  \beta(A,\zeta,w)^2:=
  \left\|\left(DF_{A}(\zeta,w)|_{\C\times w^\perp}\right)^{-1} 
  F_{A}(\zeta,w) \right\|^2
  =\left\|\begin{pmatrix}\dot\lambda\\\dot v\end{pmatrix} \right\|^2
  =|\dot\lambda|^2+ \|\dot v\|^2,
\]
where $\dot \lambda,\dot v$ are given in~\eqref{eq:newton}. 
When $\|\dot{v}\|$ is small $\beta$ approximates the length of 
the Newton step also under the distance $\dist$. Indeed,  
\begin{eqnarray*}
    \dist_A((\zeta,w),(N_A(\zeta,w)))^2
    &=& |\dot \lambda|^2+d_\P(w,w-\dot v)^2\\
     &\underset{\eqref{eq:distanceP}}{=}&
      |\dot\lambda|^2+\left(\arccos\frac{|\langle w,w
    -\dot v\rangle|}{\|w\|\|w-\dot v\|}\right)^2\\
   &=&|\dot \lambda|^2+\left(\arccos\frac{1}{\sqrt{1+\|\dot v\|^2}}\right)^2.
\end{eqnarray*}
It is easy to check that
\[
 \frac{9}{10}\,x^2 \ \le\
\Big(\arccos\frac{1}{\sqrt{1+x^2}}\Big)^2 \ \le\ x^2 \quad \mbox{ for
$x\le \frac13$,}
\]
so we have proved that whenever $\|\dot v\|\leq 1/3$,
\begin{equation}\label{eq:betadist}
   \frac{9}{10}\beta(A,\zeta,w)\leq \dist_A((\zeta,w),(N_A(\zeta,w)))
   \leq \beta(A,\zeta,w).
\end{equation}
(The upper bound inequality is valid regardless of the value of $\|\dot v\|$). The knowledge of $\beta$ allows us to bound the distance from an 
approximate eigenpair to the associated exact eigenpair.

\begin{lemma}\label{lem:beta1}
Assume that $(\zeta,w)$ is an approximate eigenpair of $A\in\SS$
with associated eigenpair $(\lambda,v)$. Then,
$\dist_A((\zeta,w),(\lambda,v))\leq 2\beta(A,\zeta,w)$. Moreover, if
$\beta(A,\zeta,w)\leq 1/3$ we also have
\[
    \frac{1}{2}\beta(A,\zeta,w)\leq \dist_A((\zeta,w),(\lambda,v))
    \leq 2\beta(A,\zeta,w).
\]
\end{lemma}

\proof
We have
\begin{align*}
 \dist_A((\zeta,w),(\lambda,v))\leq\; 
 &\dist_A((\zeta,w),(N_A(\zeta,w)))+\dist_A((N_A(\zeta,w)),(\lambda,v))\\
 \underset{\eqref{eq:betadist}}{\leq}&\beta(A,\zeta,w)+\frac{1}{2}
 \dist_A((\zeta,w),(\lambda,v)),
\end{align*}
the last from the definition of approximate eigenpair. 
The upper bound in the statement follows. Using a similar argument,
\begin{align*}
  \dist_A((\zeta,w),(\lambda,v))
  \geq\; &\dist_A((\zeta,w),(N_A(\zeta,w)))
  -\dist_A((N_A(\zeta,w)),(\lambda,v))\\
  \underset{\eqref{eq:betadist}}{\geq}\;&\frac{9}{10}
  \beta(A,\zeta,w)-\frac{1}{2}\dist_A((\zeta,w),(\lambda,v)),
\end{align*}
and the lower bound follows as well.
\eproof

\subsection{The Lipschitz property of $\beta$}\label{sec:lip-beta}

To describe {\sf Choose\_step} we need a set of constants 
satisfying a few relations. Not all of them are used in the description of 
{\sf Choose\_step}. Some only occur in the proof of the 
correctness of \EC. We will consider constants 
$c_1,c_1',c_u,c_u',c_*,c_4,c_5,c_6,c_7,K$. These are any collection 
of positive numbers satisfying the following:
\[
   \sqrt{3} c_1'\leq c_1<\frac12,\qquad 
   \sqrt{3} c_u'\leq c_u-\frac{\frac32c_1^2(\sqrt{3}-1)}{1-3c_1}
\]
\[
  4\sqrt{3}c_*< 1,\qquad c_4=c_*+(1+4\sqrt{3}c_*)(c_1+2c_u),
  \qquad 4\sqrt{3}c_4< 1,
\]
\[
   2(1+4\sqrt{3}c_*)(1+4\sqrt{3}c_4)c_u<Kc_*<\frac15,\quad 
  c_5=  \frac{c_u'}{(1+4\sqrt{3}c_*)}
  -2\frac{2c_*+\frac32c_1^2(1+4\sqrt{3}c_*)}{(1-3c_1)},
\]
\[
  c_6=\frac{c_5(1-3c_1)-2(1+3c_1)c_*}{2(1+3c_1)(1+4\sqrt{3}c_4)},
  \qquad c_7=\min\left(\frac{c_1'}{(1+4\sqrt{3}c_4)(1+4\sqrt{3}c_*)},c_6\right).
\]
A collection of parameters satisfying these constraints is 
shown in Table~1.

\begin{table}[h]
\centering
    \begin{tabular}{|l|l||l|l|}
    \hline
    $c_1'$ & $10^{-3}$ & $c_4$  & $0.005306\ldots$  \\
        \hline
    $c_1$  & $\sqrt{3}c_1'$     & $c_5$  & $0.00099\ldots$ \\
        \hline
    $c_u'$ & $10^{-3}$ & $c_6$   & $0.00038\ldots$   \\
        \hline
    $c_u$  & $\sqrt{3}c_u'+\frac{3c_1^2(\sqrt{3}-1)}{2(1-3c1)}$ & $K$ & 64\\
        \hline
    $c_*$  & $ 10^{-4}$ & $c_7$ &       $0.00038\ldots$   \\
        \hline
    \end{tabular}
    \caption{Our choice of constants for 
 \EC}
\end{table}

In addition to these constants, for the sake of clarity, we spell out a 
working hypothesis that we will repeatedly use. 

\begin{hypothesis}\label{hyp}
We have $(A,\lambda,v)\in\V$, $\|A\|_F=1$, and $(\zeta,w)$ an 
approximate eigenpair of $A$ satisfying $|\zeta|\leq1$ and
\[
  \mu(A,\lambda,v)\dist_A((\lambda,v),(\zeta,w))<c_*.
\]
Also, $\theta\colon[0,\pi)\rightarrow \SS$, $\theta(s)= A_s$ is 
some arc-length parametrized half-great circle with $A_0=A$.
\end{hypothesis}

The main goal of this section is to prove the following result.
 
\begin{lemma}\label{lem:betaestimates1}
Under Hypothesis~\ref{hyp}, let $s\leq c_1/\mu(A_0,\zeta,w)$ where 
$0\leq c_1<1/3$. Let
\[
   \Phi= \left\|\left(DF_{A_0}(\zeta,w)|_{\C\times w^\perp}\right)^{-1}
   \dot A_0 w\right\|.
\]
Then,
$$
   \beta^{-}(s)\leq \beta(A_s,\zeta,w)\leq \beta^{+}(s),
$$
 where 
$$
 \beta^{-}(s)=\frac{s\,\Phi- \beta(A_0,\zeta,w)  -\frac32c_1^2/\mu(A_0,\zeta,w)}{
 1+3c_1} ,\qquad 
\beta^{+}(s)=\frac{s\,\Phi+ \beta(A_0,\zeta,w)+\frac32c_1^2/\mu(A_0,\zeta,w)  }{
 1-3c_1}.
$$
\end{lemma}
 
For the proof we need some stepping stones.
 
\begin{lemma}\label{lem:muvsdf}
Let $A\in\SS$, $\zeta\in\C$, $|\zeta|\leq1$, $w\in\C^n$, $\|w\|=1$. Then,
$$
  \left\|\left(DF_{A}(\zeta,w)|_{\C\times w^\perp}\right)^{-1} \right\|
  \leq 3\mu(A,\zeta, w).
$$
\end{lemma}

\proof
This result is similar to ~\cite[Proposition 6.6]{Armentano:13}, but our assumptions are weaker and our definition of condition number is slightly different (see the proof of Theorem \ref{th15.1}). We can assume that $w=e_1$ for the proof and write
\[
 A=\begin{pmatrix}\lambda&a^*\\b&\hat A\end{pmatrix},
\]
where $\lambda\in\C$ and $a,b\in\C^{n-1}$. Then, for $\dot \zeta\in\C$ and $\dot w\in\C^{n-1}$,
\begin{align*}
 DF_{A}(\zeta,e_1)\left(\dot\zeta,\binom{0}{\dot w}\right)=\;&(A-\zeta\Id)\binom{0}{\dot w}-\dot\zeta e_1\\=\;&\begin{pmatrix}\lambda-\zeta&a^*\\b&\hat A-\zeta\Id_{n-1}\end{pmatrix}\binom{0}{\dot w}-\dot\zeta e_1
 \\=\;&\begin{pmatrix}-1&a^*\\0&\hat A-\zeta\Id_{n-1}\end{pmatrix}\binom{\dot\zeta}{\dot w}.
\end{align*}
We thus have
\begin{align*}
  \left\|\left(DF_{A}(\zeta,w)|_{\C\times w^\perp}\right)^{-1} \right\|=\;& \left\|\begin{pmatrix}-1&a^*\\0&\hat A-\zeta\Id_{n-1}\end{pmatrix}^{-1} \right\|\\
  =\;&\left\|\begin{pmatrix}-1&a^*(\hat A-\zeta\Id_{n-1})^{-1}\\0&(\hat A-\zeta\Id_{n-1})^{-1}\end{pmatrix} \right\|\\
  \leq\;&\left\|\begin{pmatrix}0&0\\0&(\hat A-\zeta\Id_{n-1})^{-1}\end{pmatrix} \right\|+\left\|\begin{pmatrix}-1&a^*(\hat A-\zeta\Id_{n-1})^{-1}\\0&0\end{pmatrix} \right\|\\
  \leq \;&\|A_{\zeta,w}^{-1}\|+\sqrt{1+\|a\|^2\|A_{\zeta,w}^{-1}\|^2}
  \\
  = \;&\mu(A,\zeta,w)+\sqrt{1+\|a\|^2\mu(A,\zeta,w)^2},
  %   \\
%   =\;&\sqrt{1+\|a\|^2}\max(1,\|A_{\zeta,w}^{-1}\|).
%   \\=\;& \sqrt{1+\|a\|^2}\max(1,\mu(A,\zeta,w)),
\end{align*}
the last two claims in this chain by definition of $A_{\zeta,w}$ and the fact that $\|A\|_F=1$. Now, from Lemma \ref{le:lb_mu2} we know that $\mu(A,\zeta,w)\geq(1+\sqrt{1-\|a\|^2})^{-1}$ which implies
\begin{align*}
 \frac{\mu(A,\zeta,w)+\sqrt{1+\|a\|^2\mu(A,\zeta,w)^2}}{\mu(A,\zeta,w)}=\;&1+\sqrt{\|a\|^2+\frac{1}{\mu(A,\zeta,w)^2}}\\
 \leq\;& 1+\sqrt{\|a\|^2+(1+\sqrt{1-\|a\|^2})^2}\\=\;&1+\sqrt{2+2\sqrt{1-\|a\|^2}}\leq3,
\end{align*}
the last inequality due to $0\leq \|a\|\leq \|A\|_F=1$. We have proved that
\[
 \left\|\left(DF_{A}(\zeta,w)|_{\C\times w^\perp}\right)^{-1} \right\|\leq\mu(A,\zeta,w)+\sqrt{1+\|a\|^2\mu(A,\zeta,w)^2}\leq  3\mu(A,\zeta,w),
\]
and the bound claimed in the lemma follows.

% The result 
% in~\cite{Armentano:13} assumes that $Aw=\zeta w$ but this is actually 
% not used in the proof. Due to the slight difference between our condition number and that of \cite{Armentano:13} we point out that the result is still true with our definition. Indeed, following the proof of ~\cite[Proposition 6.6]{Armentano:13}, it suffices to take $v=1$ and
% \[
%  A=\begin{pmatrix}\lambda& a\\0&\hat A\end{pmatrix},
% \]
% which yields
% \[
%  \left\|\left(DF_{A}(\zeta,w)|_{\C\times w^\perp}\right)^{-1} \right\|=\left\|\left(DF_{A}(\zeta,w)|_{\C\times w^\perp}\right)^{-1} \right\|
% \]

\eproof
 
\begin{lemma}\label{lem:auxauxbeta}
Under Hypothesis~\ref{hyp}, let $s>0$ satisfy $3\mu(A_0,\zeta,w)\,s<1$. 
Then,
\begin{equation*}
  \left\|\left[\left(DF_{A_s}(\zeta,w)|_{\C\times w^\perp}\right)^{-1}
  \left(DF_{A_0}(\zeta,w)|_{\C\times w^\perp}\right)\right]\right\| 
  \leq \frac{1}{1-3\mu(A_0,\zeta,w)\,s}  
\end{equation*}
and
\begin{equation*}
  \left\|\left[\left(DF_{A_0}(\zeta,w)|_{\C\times w^\perp}\right)^{-1}
  \left(DF_{A_s}(\zeta,w)|_{\C\times w^\perp}\right)\right]\right\|
  \leq 1+3\mu(A_0,\zeta,w)\,s.  
\end{equation*}
\end{lemma}

\proof
Note that
$$
   DF_{A_0}(\zeta,w)|_{\C\times w^\perp}(\dot\eta,\dot x)
  =DF_{A_s}(\zeta,w)|_{\C\times w^\perp}(\dot\eta,\dot x)+(A_0-A_s)\dot x,
$$
hence 
\begin{align*}
 \Big\|\Id - & 
 \left(DF_{A_0}(\zeta,w)|_{\C\times w^\perp}\right)^{-1}
\left(DF_{A_s}(\zeta,w)|_{\C\times w^\perp}\right)
 \Big\|\\
& \leq 
  \left\|\left(DF_{A_0}(\zeta,w)|_{\C\times w^\perp}\right)^{-1}\right\|
  \|A_s-A_0\|\underset{\text{Lemma~\ref{lem:muvsdf}}}{\leq}
  3\mu(A_0,\zeta,w)\,s
\end{align*}
where we used that $\|A_s-A_0\|\leq \|A_s-A_0\|_F\leq s$.  
The second claim in the statement is now obvious and the first one 
follows from the Banach lemma, $\|\Id+\Delta\|\leq (1-\|\Delta\|)^{-1}$, 
valid for $\|\Delta\|<1$. 
\eproof

\begin{lemma}\label{lem:betaestimatesaux}
Under Hypothesis~\ref{hyp}, let $s>0$ satisfy $3\mu(A_0,\zeta,w)\,s<1$ 
and let
\[
  \Phi_{\mathsf{aux}}(s)= 
  \left\|\left(DF_{A_0}(\zeta,w)|_{\C\times w^\perp}\right)^{-1}
  (A_0-A_s)w\right\|.
\]
Then,
$$
   \frac{\Phi_{\mathsf{aux}}(s)- \beta(A_0,\zeta,w) }
    {1+3\mu(A_0,\zeta,w)\,s}\leq \beta(A_s,\zeta,w)
    \leq \frac{\Phi_{\mathsf{aux}}(s)+\beta(A_0,\zeta,w)}
    {1-3\mu(A_0,\zeta,w)\,s}.
$$
\end{lemma}
  
\proof
From the definition,
 \begin{eqnarray}\label{eq:betas}
  \beta(A_s,\zeta,w) &=& 
  \left\|\left(DF_{A_s}(\zeta,w)|_{\C\times w^\perp}\right)^{-1} 
  F_{A_{s}}(\zeta,w) \right\|\notag\\
  &=& \Big\|\left[\left(DF_{A_s}(\zeta,w)|_{\C\times w^\perp}\right)^{-1}
  \left(DF_{A_0}(\zeta,w)|_{\C\times w^\perp}\right)\right]\\
  & &\quad\left(DF_{A_0}(\zeta,w)|_{\C\times w^\perp}\right)^{-1}
  (F_{A_{0}}(\zeta,w) + (A_s-A_0)w )\Big\|.\notag
\end{eqnarray}
Thus,
\[
  \beta(A_s,\zeta,w) \leq 
  \left\|\left(DF_{A_s}(\zeta,w)|_{\C\times w^\perp}\right)^{-1}
   \left(DF_{A_0}(\zeta,w)|_{\C\times w^\perp}\right)\right\| 
  (\Phi_{\mathsf{aux}}(s)+\beta(A_0,\zeta,w)).
\]
Similarly, one shows 
\[
  \beta(A_s,\zeta,w)\geq 
  \frac{\Phi_{\mathsf{aux}}(s)-\beta(A_0,\zeta,w)}
  {\left\|\left(DF_{A_0}(\zeta,w)|_{\C\times w^\perp}\right)^{-1}
  \left(DF_{A_s}(\zeta,w)|_{\C\times w^\perp}\right)\right\|.}
\]
The statement now follows from Lemma~\ref{lem:auxauxbeta}.
\eproof
\medskip

\proofof{Lemma~\ref{lem:betaestimates1}}
We claim that we can write
\[
A_s=A_0+s\dot A_0+\frac{s^2}{2}B,
\]
for some $B$ with $\|B\|_F\leq1$. Indeed, this is an elementary observation which follows from the fact that $A_s$ is a great circle in the sphere (w.l.o.g. one can choose $A_s$ to be the circle parametrized by $(\cos s,\sin s,0\ldots,0)$ to prove it). As a consequence we have
\begin{align*}
  |\Phi_{\mathsf{aux}}(s)-s\Phi |\leq\; &  
  \left\|\left(DF_{A_0}(\zeta,w)|_{\C\times w^\perp}\right)^{-1}
  (A_0-A_s+s\dot A_0)w \right\|\\
  \leq\; & \frac{s^2}{2}
  \left\|\left(DF_{A_0}(\zeta,w)|_{\C\times w^\perp}\right)^{-1} \right\|
 \underset{\text{Lemma~\ref{lem:muvsdf}}}{\leq}\frac32s^2\mu(A_0,\zeta,w).
\end{align*}
The upper and lower bounds for $\beta(A_s,\zeta,w)$ in 
Lemma~\ref{lem:betaestimates1} now follow from this last estimate and 
Lemma~\ref{lem:betaestimatesaux}.
\eproof

\subsection{The step's length in the homotopy continuation}
\label{sec:onestep}

The following result is crucial for the understanding of 
the homotopy algorithm. Its proof follows a logic which is similar to that of the proof of Corollary~\ref{cor:b}.

\begin{proposition}\label{prop:republica}
Under Hypothesis~\ref{hyp}, let $s'>0$ be any number such that
\[
   c'_1\leq \mu(A,\zeta,w)s'\leq c_1
\]
and let $s''$ be any number such that
\[
  c_u'\leq \mu(A,\zeta,w) \beta^{+}(s'')\leq c_u,
\]
where $\beta^{+}$ is as in Lemma~\ref{lem:betaestimates1}.  
Let $\bar{s}=\min(s',s'')$.
Then, 
\begin{enumerate}
\item[I.] 
The (continuous) branch of the solution map 
$\pi_1^{-1}\circ \theta\colon [0,\bar{s}]\rightarrow\V$ 
with $\pi_1^{-1}\circ\theta(0)=(A,\lambda,v)$ is well defined.
\item[II.] 
For every $s\in[0,\bar{s}]$, let 
$(A_s,\lambda_s,v_s):=\pi_1^{-1}\circ\theta(s)$. Then, 
\begin{equation}\label{eq:nuevaaux}
  \frac{1}{1+4\sqrt{3} c_4} \mu(A,\lambda,v)
 < \mu(A_s,\lambda_s,v_s)< (1+4\sqrt{3} c_4) \mu(A,\lambda,v)
\end{equation}
and
\begin{equation}\label{eq:nuevaaux3}
  \dist_{A_s}((\lambda_s,v_s),(\zeta,w))
 <\frac{Kc_*}{\mu(A_s,\lambda_s,v_s)}.
\end{equation}

\end{enumerate}
In particular, $(\zeta,w)$ is an approximate eigenpair of $A_s$ 
with associated eigenpair $(\lambda_s,v_s)$ for every $s\in[0,\bar{s}]$. 

Finally,  the condition length $L_{\mu,0,\bar{s}}(\pi_1^{-1}(A_s))$ 
(see~\eqref{eq:condlength}) of the curve $\pi_1^{-1}(A_s)$ is 
at least $c_7$.
\end{proposition}
Before we prove Proposition \ref{prop:republica} we make a few comments on how the proof of the proposition and Theorem \ref{thm:main_path_following} proceed. The hypotheses $\beta^+(s'')$ give us a bound on the distance from $(\zeta,w)$ to $(\lambda_s,v_s)$ and ultimately $(\lambda,v)$ to $(\lambda_s,v_s)$. Together with the bound on $s'$ this allows us to invoke Theorem \ref{th15.1} and Theorem \ref{thm:main_path_following} to prove conclussions I and II of the proposition. The tricky part will be to prove the last statement that $L_{\mu,0,\bar{s}}(\pi_1^{-1}(A_s))\geq c_7$. This last statement gives us the upper bound on the number of steps in Theorem \ref{thm:main_path_following} as the condidion length divided by $c_7$. Now to see that
\[
 \int_0^{\bar s}\mu(A_s,\lambda_s,v_s)\|(\dot A_s,\dot \lambda_s,\dot v_s)\|\,ds=L_{\mu,0,\bar{s}}(\pi_1^{-1}(A_s))\geq c_7
\]
it will suffice to prove that
\[
 \int_0^{\bar s}\|(\dot A_s,\dot \lambda_s,\dot v_s)\|\,ds\geq\frac{constant}{\mu(A,\lambda,v)},
\]
since $\mu$ is almost constant on the interval. For this it will suffice to prove that one of
\[
 \int_0^{\bar s}\|\dot A_s\|\,ds\quad\text{ or }\quad \int_0^{\bar s}\|\dot v_s\|\,ds
\]
is greater than $constant/\mu(A,\lambda,v)$. The first integral is $\bar{s}$. We will see that if $\bar{s}$ is small then $\dist_A(v_s,v)$ is greater than or equal to some constant over $\mu(A,\lambda,v)$, and so is the integral of $\|\dot v_s\|$ which is the length of a path between $v$ and $v_s$.
\bigskip

\proofof{Proposition \ref{prop:republica}}
We prove the first part of the proposition. From the inverse function
theorem and the continuity of $\mu$, there exists a maximal $s^*\leq
\bar s$ such that I and II hold changing $[0,\bar{s}]$ to $[0,s^*)$. % We just need to prove that $s^*>\bar{s}$. 
The global upper bound for $\mu(\pi_1^{-1}(A_s))$ shown in \eqref{eq:nuevaaux} is in turn an upper bound for the derivative of the solution
map, for $s\in[0,s^*)$. A standard limit argument in compact sets then
implies that $\pi_1^{-1}\circ\theta$ can be extended in a continuous
manner to $[0,s^*]$, and because \eqref{eq:nuevaaux} and \eqref{eq:nuevaaux3} are open conditions we must have one of the two following scenarios:
\begin{enumerate}
 \item[i)] $s_*=\bar s$ and both \eqref{eq:nuevaaux} and \eqref{eq:nuevaaux3} hold changing $s$ to $s_*$, or
 \item[ii)] At least one of  \eqref{eq:nuevaaux} and \eqref{eq:nuevaaux3} does not hold changing $s$ to $s_*$.
\end{enumerate}
We now discard the second option. Note that from Hypotheses \ref{hyp} and Theorem~\ref{th:lipsch},
\begin{equation}\label{eq:one*}
  \frac{1}{1+4\sqrt{3} c_* }\mu(A,\lambda, v)< \mu(A,\zeta,w)
<(1+4\sqrt{3} c_*)\mu(A,\lambda, v).
\end{equation}
Then, we have for every $s\in[0,s^*]$
\[
  \dist((A,\zeta,w),(A_s,\zeta,w))=\|A-A_s\|_F< s^*\leq s'
  \leq\frac{c_1}{\mu(A,\zeta,w)}\leq 
  \frac{(1+4\sqrt{3} c_*)c_1}{\mu(A,\lambda, v)},
\]
and because $s^*\leq s''$, from Lemma~\ref{lem:betaestimates1} we have
\[
   \beta(A_s,\zeta,w)\leq \beta^{+}(s)\leq 
   \frac{c_u}{\mu(A,\zeta,w)}< 
   \frac{(1+4\sqrt{3} c_*)c_u}{\mu(A,\lambda, v)}.
\]
From Lemma~\ref{lem:beta1} (recall that from II and Theorem \ref{th15.1},
$(\zeta,w)$ is an approximate eigenpair of $A_s$ with associated 
eigenpair $(\lambda_s,v_s)$ for $s\leq s_*$), this last inequality implies
\begin{equation}\label{eq:nuevaaux4}
  \dist_{A_s}((\zeta,w),(\lambda_s,v_s))\leq 
  2\beta(A_s,\zeta,w)< 
  \frac{2(1+4\sqrt{3} c_*)c_u}{\mu(A,\lambda, v)}.
\end{equation}
We have thus proved that for every $s\in[0,s^*]$,
\begin{align*}
  \dist((A,\lambda,v),(A_s,\lambda_s,v_s))\leq\;&
  \dist_{A}((\lambda,v),(\zeta,w))+\dist((A,\zeta,w),(A_s,\zeta,w))\\
  &\qquad+\dist_{A_s}((\zeta,w),(\lambda_s,v_s))\\
 <\;&\frac{c_*}{\mu(A,\lambda,v)}+\frac{(1+4\sqrt{3} c_*)c_1}
  {\mu(A,\lambda, v)}+\frac{2(1+4\sqrt{3} c_*)c_u}{\mu(A,\lambda, v)}\\
  =\;&\frac{c_4}{\mu(A,\lambda, v)}.
\end{align*}
Then, from Theorem~\ref{th:lipsch} (note the strict inequality in the displayed formula above),
\begin{equation}\label{eq:nuevaaux2}
   \frac{1}{1+4\sqrt{3} c_4} \mu(A,\lambda,v)
  < \mu(A_{s^*},\lambda_{s^*},v_{s^*})< (1+4\sqrt{3} c_4) \mu(A,\lambda,v).
\end{equation}
Thus, \eqref{eq:nuevaaux} holds at $s^*$ and moreover
\[
\dist_{A_{s^*}}((\lambda_{s^*},v_{s^*}),(\zeta,w))\underset{\eqref{eq:nuevaaux4},\eqref{eq:nuevaaux2}}<\frac{2(1+4\sqrt{3} c_*)(1+4\sqrt{3}c_4)c_u}{\mu(A_{s^*},\lambda_{s^*},v_{s^*})} <\frac{Kc_*}{\mu(A_{s^*},\lambda_{s^*},v_{s^*})}. 
\]
That is, \eqref{eq:nuevaaux3} holds at $s^*$ and we can discard option ii), proving the first part of the proposition.

For the last claim of the proposition, note that the condition length of the curve $\pi_1^{-1}(A_s)$ is
\begin{eqnarray*}
  L_{\mu,0,\bar{s}}(\pi_1^{-1}(A_s))&=&
    \int_{0}^{\bar{s}} \mu(\pi_1^{-1}(A_s))
    \left\|D\pi_1^{-1}(A_s)\dot A_s\right\| \,ds\\
 &\underset{\eqref{eq:nuevaaux}}{\geq}& 
 \frac{\mu(A,\lambda,v)}{1+4\sqrt{3}c_4}\int_{0}^{\bar{s}} 
 \left\|D\pi_1^{-1}(A_s)\dot A_s\right\| \,ds.
\end{eqnarray*}

If $\bar{s}=s'$ then using that the integrand is greater than or equal to 
$1$ and~\eqref{eq:one*} we have
\[
  L_{\mu,0,\bar{s}}(\pi_1^{-1}(A_s))\geq 
\frac{\mu(A,\lambda,v)}{1+4\sqrt{3}c_4}s'
  \geq\frac{c_1'}{1+4\sqrt{3}c_4}\frac{\mu(A,\lambda,v)}{\mu(A,\zeta,w)}
  \geq\frac{c_1'}{(1+4\sqrt{3}c_4)(1+4\sqrt{3}c_*)}.
\]

Assume now that $\bar{s}=s''$. Then we have
\begin{eqnarray*}
  L_{\mu,0,\bar{s}}(\pi_1^{-1}(A_s))&\geq& 
 \frac{\mu(A,\lambda, v)}{1+4\sqrt{3}c_4}
  \int_{0}^{\bar{s}} \left\|D\pi_1^{-1}(A_s)\dot A_s\right\| \,ds\\
  &\geq& \frac{\mu(A,\lambda, v)}{1+4\sqrt{3}c_4}
  \dist((A_0,\lambda_0,v_0),(A_{\bar{s}},\lambda_{\bar{s}},v_{\bar{s}})).
\end{eqnarray*}
Now, note that (recall $(A_0,\lambda_0,v_0)=(A,\lambda,v)$)
\begin{align*}
  \dist((A_0,\lambda_0,v_0),&(A_{\bar{s}},\lambda_{\bar{s}},v_{\bar{s}})) \\
  \geq\;& \dist((A_0,\zeta,w),(A_{\bar{s}},\lambda_{\bar{s}},v_{\bar{s}}))
  -\dist_{A_0}((\lambda_0,v_0),(\zeta,w))\\
  \geq\;& \dist_{A_{\bar{s}}}((\zeta,w),(\lambda_{\bar{s}},v_{\bar{s}}))
  -\frac{c_*}{\mu(A,\lambda,v)}.
\end{align*}
We need a lower bound for this last term. We first note that from 
Lemma~\ref{lem:beta1} and Hypothesis~\ref{hyp}, 
\begin{equation}\label{eq:two*}
  \beta(A,\zeta,w)\leq \frac{2c_*}{\mu(A,\lambda,v)}\leq  
  \frac{2(1+4\sqrt{3}c_*)c_*}{\mu(A,\zeta,w)},
\end{equation}
the last by~\eqref{eq:one*}. 
Using this last bound and, again, Lemmas~\ref{lem:beta1} 
and~\ref{lem:betaestimates1}, we have
\begin{align*}
  2\frac{1+3c_1}{1-3c_1}&
  \dist_{A_{\bar{s}}}((\zeta,w),(\lambda_{\bar{s}},v_{\bar{s}}))
  \geq \frac{1+3c_1}{1-3c_1}\beta(A_{\bar{s}}, \zeta,w)\\
  \geq\;& \frac{1+3c_1}{1-3c_1}\beta^{-}(\bar{s})\;=\;
  \beta^{+}(\bar{s})-2\frac{\beta(A,\zeta,w)+\frac32c_1^2/\mu(A,\zeta,w)}{1-3c_1}\\
  \geq\;& \frac{c_u'}{\mu(A,\zeta,w)}-
  2\frac{\beta(A,\zeta,w)+\frac32c_1^2/\mu(A,\zeta,w)}{1-3c_1}\\
  \geq\;&
  \frac{c_u'}{(1+4\sqrt{3}c_*)\mu(A,\lambda, v)}
  -2\frac{2c_*+\frac32c_1^2(1+4\sqrt{3}c_*)}{(1-3c_1)\mu(A,\lambda,v)}
  \;=\;\frac{c_5}{\mu(A,\lambda, v)}.
\end{align*}
For the inequality in the third line we used the assumption in 
the statement and $s=s''$. For the inequality in the fourth line,  
we used~\eqref{eq:one*} and~\eqref{eq:two*}. For the 
equality in the fourth line we used the definition of $c_5$.

We have thus shown that if $\bar{s}=s''$, then
\begin{eqnarray*}
  L_{\mu,0,\bar{s}}(\pi_1^{-1}(A_s))&\geq&
  \frac{\mu(A,\lambda, v)}{1+4\sqrt{3}c_4}
  \left(\frac{c_5(1-3c_1)}{2(1+3c_1)
  \mu(A,\lambda, v)}-\frac{c_*}{\mu(A,\lambda, v)}\right)\\
 &=& \frac{c_5(1-3c_1)-2(1+3c_1)c_*}{2(1+3c_1)(1+4\sqrt{3}c_4)}
 \;=\;c_6.
\end{eqnarray*}
Since $c_7\leq c_6$ the proof is complete.
\eproof

We can finally describe the subroutine {\sf Choose\_step}. It is
important to note that given any matrix $A\in\mnc$, one can compute in
$\Oh(n^3)$ arithmetic operations, for example by first reducing $A$ to
tridiagonal Hessenberg form and then using the main result
of~\cite{Kahan1966}, a number $r$ such that $\|A\|\leq r\leq
\sqrt{3}\|A\|$. That is, we can compute operator norms within a factor
of $\sqrt{3}$ and, consequently, we can compute $\mu$ within a 
factor of $\sqrt{3}$. 
\bigskip

\algoritmo
\begin{algorithm}\label{alg:choosestep}
{\sf Choose\_step}\\
\inputalg{$B,\dot A\in \SS$,  and 
$(\zeta,w)\in\C\times\C^n$, $\|w\|=1$}\\
\bodyalg{
compute $r>0$ such that 
$\mu(B,\zeta,w)\leq r\leq \sqrt{3}\mu(B,\zeta,w)$\\[2pt]
$s':=c_1/r$\\[2pt]
$\Phi:= \left\|\left(DF_{B}(\zeta,w)|_{\C\times  w^\perp}\right)^{-1}
 \dot Aw\right\|$\\[2pt]
compute $s''$, the solution of the following linear equation
\[
  \frac{\Phi\,s''+ \beta(B,\zeta,w)+\frac32c_1^2\sqrt{3}/r  }{ 1-3c_1 }
 =\frac{c_u}{r}
\]
\\[2pt]
$\bar{s}:=\min(s',s'')$\\[2pt]
}
\Output{$\Delta s=\bar{s}\in[0,\pi]$}
\end{algorithm}
\falgoritmo

The step size computed by {\sf Choose\_step} cannot be too small, as we show now.

\begin{proposition}\label{prop:step-size}
The value $\Delta s$ returned by 
{\sf Choose\_step}$(B,\dot{A},\zeta,w)$ satisfies
\[
   \Delta s\geq \frac{R}{\mu^2(B,\zeta,w)}
\] 
with $R=c_7(6(1+4\sqrt{3}c_4)^2)^{-1}>0$.
\end{proposition}

\subsection{Proof of Theorem~\ref{thm:main_path_following} 
and Proposition \ref{prop:step-size}}\label{sec:homotopy}

We now prove Theorem \ref{thm:main_path_following}, and the 
proof of  Proposition \ref{prop:step-size} will follow straighforward 
from our arguments. 

From the definition of \EC\ it is clear that we can assume that
$\|A_0\|_F=\|A\|_F=1$. We further assume that the 
constants $c'_1,c_1,c'_u,c_uc_*,c_4,c_5,c_6,c_7$ and $K$ take the
values in Table~1 and denote by
$\pi_1^{-1}(\mathcal{L}_{A_0,A})$ the lift of $\mathcal{L}_{A_0,A}$
with origin $(A_0,\lambda_0,v_0)$. Note that this lift is well-defined 
since, by hypothesis, $\mathcal{L}_{A_0,A}\cap\Sigma=\emptyset$.

Let $B_i$ be the matrix $B$ at the beginning of the $i$th iteration 
of {\sf Path-follow}. Also, let $(\lambda_i,v_i)$ be such that 
$(B_i,\lambda_i,v_i)$ is the (only) triple in 
$\pi_1^{-1}(\mathcal{L}_{A_0,A})$ above $B_i$.  

We first prove that for all $i\geq 0$, $(\zeta_{i},{w}_i)$ is an approximate zero of $B_i$  
with associated eigenpair $(\lambda_i,v_i)$ and satisfies
\[
 \dist_{B_{i}}\left((\zeta_{i},w_i),(\lambda_{i},v_{i})\right)<\frac{c_*}{\mu(B_{i},\lambda_{i},v_{i})}.
\]
We reason by induction. 
The step $i=0$ is true by hypothesis (recall Definition \ref{def:fund_neigh}). For the induction step, 
note that the $s'$ defined by {\sf Choose\_step} satisfies 
(we omit the subindices $i$ in $(B_i,\zeta_i,w_i)$ in the next few lines) 
\[
  \frac{c_1'}{\mu(B,\zeta,w)}\leq 
  \frac{c_1}{\sqrt{3}\mu(B,\zeta,w)}\leq s'\leq\frac{c_1}{\mu(B,\zeta,w)}.
\]
Moreover, $s''$ satisfies
\begin{eqnarray*}
  \beta^+(s'')&=&
  \frac{\Phi\,s''+ \beta(B,\zeta,w)+\frac32c_1^2/\mu(B,\zeta,w)  }{ 1-3c_1 }\\
  &\leq& \frac{\Phi\,s''+ \beta(B,\zeta,w)+\frac32c_1^2\sqrt{3}/r  }{ 1-3c_1 }
  \,=\,\frac{c_u}{r}\leq \frac{c_u}{\mu(B,\zeta,w)},
\end{eqnarray*}
and
\begin{align*}
\beta^+(s'') =\; &
\frac{\Phi\,s''+ \beta(B,\zeta,w)+\frac32c_1^2/\mu(B,\zeta,w)  }{ 1-3c_1 }
\;\geq\; \frac{\Phi\,s''+ \beta(B,\zeta,w)+\frac32c_1^2/r  }{ 1-3c_1 }\\
=\;&
  \frac{c_u-\frac32c_1^2(\sqrt{3}-1)/(1-3c_1)}{r}\;\geq\; 
  \frac{c_u-\frac32c_1^2(\sqrt{3}-1)/(1-3c_1)}{\sqrt{3}\mu(B,\zeta,w)}
  \;\geq\;\frac{c_u'}{\mu(B,\zeta,w)}.
\end{align*}
We are thus under the hypothesis of Proposition~\ref{prop:republica} 
with $\bar{s}=\Delta s$ (by construction in {\sf Choose\_step}),   
which guarantees that $(\zeta_i,w_i)$ 
is an approximate eigenpair of the matrix $B_{i+1}$. Moreover, 
Proposition~\ref{prop:republica} also implies
\[
  \dist_{B_{i+1}}((\zeta_i,w_i),(\lambda_{i+1},v_{i+1}))
 <\frac{Kc_*}{\mu(B_{i+1},\lambda_{i+1},v_{i+1})}.
\]
Since $Kc_*<1/5$ we deduce (using 
Theorem~\ref{th15.1}) that $(\zeta_i,w_i)$ is an approximate eigenpair 
of $B_{i+1}$ with associated eigenpair $(\lambda_{i+1},v_{i+1})$. 
Consequently, after three steps of Newton iteration, we have 
that $(\zeta_{i+1},w_{i+1})$ satisfies (recall, $K=64$) 
\begin{eqnarray}\label{eq:finalista}
  \dist_{B_{i+1}}((\zeta_{i+1},w_{i+1}),(\lambda_{i+1},v_{i+1}))
  &<& \frac{1}{2^{2^3-1}}\frac{64c_*}{\mu(B_{i+1},\lambda_{i+1},v_{i+1})}\nonumber\\
  &=&\frac{c_*}{2\mu(B_{i+1},\lambda_{i+1},v_{i+1})}.
\end{eqnarray}
If $|\zeta_{i+1}|\leq 1$, this finishes the induction step. Otherwise, the algorithm divides $\zeta_{i+1}$ by its norm, in that case we have
\begin{align*}
\dist_{B_{i+1}}\left(\left(\frac{\zeta_{i+1}}{|\zeta_{i+1}|},w_{i+1}\right),(\lambda_{i+1},v_{i+1})\right)\leq \;&
 \dist_{B_{i+1}}\left(\left(\frac{\zeta_{i+1}}{|\zeta_{i+1}|},w_{i+1}\right),(\zeta_{i+1},w_{i+1})\right)+\\&\qquad\dist_{B_{i+1}}((\zeta_{i+1},w_{i+1}),(\lambda_{i+1},v_{i+1}))
 \\
 \leq \;&|\zeta_{i+1}|-1+\frac{c_*}{2\mu(B_{i+1},\lambda_{i+1},v_{i+1})}
\end{align*}
On the other hand, from \eqref{eq:finalista} we have (use $|\la_{i+1}| \le \|B_{i+1}\|_F=1$)
\[
 |\zeta_{i+1}-\lambda_{i+1}|<\frac{c_*}{2\mu(B_{i+1},\lambda_{i+1},v_{i+1})}\Rightarrow |\zeta_{i+1}|< 1+\frac{c_*}{2\mu(B_{i+1},\lambda_{i+1},v_{i+1})},
\]
so we have
\[
 \dist_{B_{i+1}}\left(\left(\frac{\zeta_{i+1}}{|\zeta_{i+1}|},w_{i+1}\right),(\lambda_{i+1},v_{i+1})\right)<\frac{c_*}{\mu(B_{i+1},\lambda_{i+1},v_{i+1})},
\]
and the induction step is finished also in the case that $|\zeta_{i+1}|>1$.

The induction step is complete. In particular, this shows the last part 
of the statement. 

To show the complexity bounds, assume \EC\ has performed $q+\ell$ 
iterations and let $0=s_0<s_1<\ldots<s_q<\ldots<s_{q+\ell}$ 
be the corresponding values of $s$. Then we have
\begin{eqnarray*}
  L_{\mu,s_q,s_{q+\ell}}(\pi_1^{-1}(B_s))&=&  
  \sum_{i=1}^\ell \int_{s_{q+i-1}}^{s_{q+i}}\mu(A_s,\lambda_s,v_s)
  \|(\dot A_s,\dot\lambda_s,\dot v_s)\|\,ds\\ 
  &=& \sum_{i=1}^\ell L_{\mu,s_{q+i-1},s_{q+i}}(\pi_1^{-1}(A_s))
  \geq \ell c_7
\end{eqnarray*}
the inequality by the last claim of
Proposition~\ref{prop:republica}. But the algorithm 
halts as soon as $s_{q+\ell}=\alpha$, i.e., as soon as 
$$
   L_{\mu,s_q,s_{q+\ell}}(\pi_1^{-1}(A_s))
   =L_{\mu,s_q,\alpha}(\pi_1^{-1}(A_s)),
$$
which occurs as soon as 
$\ell\geq c_7^{-1}\int_{s_q}^{\alpha}\mu(A_s,\lambda_s,v_s)
  \|(\dot A_s,\dot\lambda_s,\dot v_s)\|\,ds$,  
as claimed in the theorem (note that $C:=c_7^{-1}\leq3000$).

We finally prove Proposition \ref{prop:step-size}. From 
Proposition~\ref{prop:republica} we have
\begin{eqnarray*}
  \mu^2(B,\zeta,w)\Delta s 
 &\underset{\eqref{eq:nuevaaux}}{\geq}&
  \frac{1}{(1+4\sqrt{3}c_4)^2}\int_0^{\Delta s}
  \mu^2(A_s,\lambda_s,v_s)\,ds\\
 &\underset{\text{Prop. \ref{prop:dotzdotf-spherical} }}{\geq}&
  \frac{1}{6(1+4\sqrt{3}c_4)^2}L_{\mu,0,\Delta_s}(\pi_1^{-1}(A_s)) 
  \geq \frac{c_7}{6(1+4\sqrt{3}c_4)^2},
\end{eqnarray*}
so Proposition \ref{prop:step-size} follows.
\eproof

\section{Integration in the solution variety}

\subsection{The coarea formula}

On a Riemannian manifold~$M$ 
there is a well-defined measure $\vol_M$ obtained 
by integrating the indicator functions $\uno_A$ 
of Borel-measurable subsets $A\subseteq M$ 
against the volume form $dM$ of~$M$,
$$
    \vol_M(A) := \int_M  \uno_A \, dM. 
$$  
Dividing $\uno$ by $\vol_M(M)$ if $\vol_M(M)<\infty$, 
this leads to a natural notion of {\em uniform distribution}
on $M$, which we will denote by $\msU(M)$. 
More generally, we will call any measurable function 
$f\colon M\to [0,\infty]$ such that $\int_M f\, dM = 1$
a {\em probability density} on~$M$.

The coarea formula (a modern classical formula due 
to Federer~\cite{Fe69}, see the Appendix of~\cite{howa:93}  
for a smooth version)  is an extension of the transformation formula  
to not necessarily bijective smooth maps 
between Riemannian manifolds. 
In order to state it, we first need to generalize 
the notion of Jacobians.

Suppose that $M,N$ are Riemannian manifolds of 
dimensions~$m$, $n$,
respectively such that $m\ge n$. Let $\psi\colon M\to N$ be
a smooth map. By definition, the derivative $D\psi(x)\colon
T_xM\to T_{\psi(x)} N$ at a regular point $x\in M$ is surjective.
Hence the restriction of $D\psi(x)$ to the orthogonal complement of
its kernel yields a linear isomorphism. The absolute value of its
determinant is called the
{\em normal Jacobian} (sometimes called {\em normal determinant} 
in the context of Linear Algebra, see~\cite{PeterAmeluxen}) 
of $\psi$ at $x$ and is denoted by $\NJ\psi(x)$.
We set $\NJ\psi(x):=0$ if $x$~is not a regular point.  

If $y$ is a regular value of~$\psi$, then 
the fiber $F_y:=\psi^{-1}(y)$ is a Riemannian submanifold 
of~$M$ of dimension $m-n$, and it makes sense to integrate functions on $F_y$. Moreover,   
Sard's lemma states that almost all $y\in N$ are regular values.

We can now state the {\em coarea formula}.

\begin{theorem}[Coarea formula]\label{pro:coarea}
Suppose that $M,N$ are Riemannian manifolds of 
dimensions~$m$, $n$, 
respectively, and let $\psi\colon M\to N$ be a surjective 
smooth map such that $D\psi$ is surjective a.e. Put $F_y=\psi^{-1}(y)$.
Then we have for any function $\chi\colon M\to\R$ that is integrable
with respect to the volume measure of $M$ that
\begin{equation*}
  \int_M \chi \, dM= \int_{y\in N} \left( \int_{F_y}
  \frac{\chi}{\NJ\psi}\, dF_y\right) dN ,
\end{equation*}
and the integrals involved are well-defined.
\eproof
\end{theorem}

It should be clear that this result contains 
the change of variables formula as a special case. 
Moreover, if we apply the coarea formula to the projection 
$\pi_2\colon M\times N \to N,\, (x,y)\mapsto y$, 
we retrieve Fubini's theorem since $\NJ\pi_2=1$. 

\subsection{Coarea formula and double fibrations}

The coarea formula can be readily applied to the following
situation. Assume that three Riemannian manifolds $M$, $N_1$, $N_2$
are equipped with surjective smooth mappings 
$\pi_1\colon M\rightarrow N_1$ and $\pi_2\colon M\rightarrow N_2$ whose derivatives are a.e. surjective,  so $\NJ\pi_1$ and $\NJ\pi_2$ are a.e. nonzero. Let
$\chi\colon M\rightarrow[0,\infty)$ be a measurable mapping. From 
Theorem~\ref{pro:coarea} applied to $\pi_1$ we have 
(here $dx$ and $dy$ stand for the volume forms in $M$ and $N_1$, 
respectively)
\[
   \int_{x\in M}\chi(x)\NJ\pi_1(x)\,dx=\int_{y\in N_1}
   \int_{x\in \pi_1^{-1}(y)}\chi(x)\,dx\,dy.
\]
On the other hand, Theorem~\ref{pro:coarea} applied to $\pi_2$ yields
\[
   \int_{x\in M}\chi(x)\NJ\pi_1(x)\,dx
  =\int_{z\in N_2}\int_{x\in \pi_2^{-1}(z)}\frac{\NJ\pi_1(x)}{\NJ\pi_2(x)}
   \chi(x)\,dx\,dz.
\]
We thus have the following result.

\begin{theorem}\label{th:doubleprojectioncoarea}
Let $M$, $N_1$, $N_2$ be Riemannian manifolds equipped with 
surjective smooth mappings $\pi_1\colon M\rightarrow N_1$ and 
$\pi_2:M\rightarrow N_2$ whose derivatives are a.e. surjective.  Let $\chi\colon M\rightarrow[0,\infty)$ be a 
measurable mapping. Then,
\begin{equation*}
   \int_{y\in N_1}\int_{x\in \pi_1^{-1}(y)}\chi(x)\,dx\,dy
  \,=\,\int_{z\in N_2}\int_{x\in \pi_2^{-1}(z)}
   \frac{\NJ\pi_1(x)}{\NJ\pi_2(x)}\chi(x)\,dx\,dz.
\end{equation*}
(Note that when $M$ and $N_1$ have the same dimension, one 
can replace the inner integral in the left-hand side by a 
(finite or enumerable) sum).\eproof
\end{theorem}

In the next sections we shall apply this result in two different
contexts. A linear algebra argument simplifies the computation of the
quotient of Normal Jacobians. Let $E$ and $F$ be finite dimensional,
complex Euclidean vector spaces and let $\varphi\colon E\to F$ be a
surjective linear mapping. Consider the graph
$\Gamma:=\{(x,\varphi(x)) \mid x \in E\}$ of~$\varphi$. Then, $\Gamma$
is a linear subspace of $E\times F$ and the two projections $$
p_1\colon\Gamma \to E,\ (x,\varphi(x)) \mapsto x,\quad p_2\colon\Gamma
\to F,\ (x,\varphi(x)) \mapsto \varphi(x) $$ are linear maps.  Note
that $p_1$ is an isomorphism and $p_2$ is surjective as $\varphi$ is.
\begin{lemma}\label{le:detquot}
Under the above assumptions, we have 
$$
   \frac{\NJ p_1}{\NJ p_2} = |\det(\varphi\varphi^{*})|^{-1} .
$$
\end{lemma}

\proof
This result is~\cite[Lemma~3b), p.~242]{BlCuShSm98}, although we
rewrite it for complex vector spaces here (note the comment in the
proof of~\cite[Theorem~5, p.~243]{BlCuShSm98}).
\eproof

\subsection{The solution variety for the eigenpair problem}
\label{sec:interationV}

Recall from \S\ref{sec:2.2}, we have the two projections 
$$
 \pi_1\colon \V\to \Cnn, (A,\lambda,v) \to A \qquad\mbox{and}\qquad
 \pi_2\colon \V\to \C\times\Pn
 , (A,\lambda,v) \to (\lambda,v)
$$
and, for $(A,\lambda,v)\in\V$, the linear operator 
$A_{\lambda,v}\colon v^\perp \to v^\perp$ given by 
$P_{v^\perp}(A-\lambda\,\Id)_{|v^\perp}$.   
In order to apply~Theorem~\ref{th:doubleprojectioncoarea} we first need 
to compute the quotient of normal Jacobians there.

\begin{proposition}\label{pro:NJquot}
Let  $p:=(A,\lambda,v) \in \W$ and choose a representative such 
that $\|v\|=1$. Then, the derivative
$D\pi_1(p)\colon T_p\W \to T_A\Cnn$ is an isomorphism, the derivative 
$D\pi_2(p)\colon T_p\W \to T_{(\lambda,v)}(\C\times\Pn)$ is surjective, and we have 
$$
  \frac{\NJ\pi_1(p)}{\NJ\pi_2(p)} = | \det A_{\lambda,v} |^{2} 
  = \det(A_{\lambda,v}A_{\lambda,v}^*).
$$
\end{proposition}

\proof
By orthogonal  invariance, we may assume without loss of generality 
that $v=e_1=(1,0,\ldots,0)$, and then
$$
  A = \begin{pmatrix} \lambda & c^* \\ 0 & B \end{pmatrix}, \quad 
  c \in\C^{n-1}, B \in \C^{(n-1)\times (n-1)} ,
$$ 
so $A_{\lambda,v}=B-\lambda \Id_{n-1}$. Let $\Gamma=\{(\dot
A,DG(A)\dot A):\dot A\in\mnc\}\subseteq\mnc\times\C\times
v^\perp$ where $G$ is the appropriate branch of the solution map defined in some open
neighborhood of $A$. We are under the hypotheses of
Lemma~\ref{le:detquot} so we have
\[
  \frac{\NJ\pi_1(p)}{\NJ\pi_2(p)} = \frac{\NJ(D\pi_1(p))}{\NJ(D\pi_2(p))} 
 =\det(DG(A,\lambda,v)DG(A,\lambda,v)^*)^{-1}.
\]
From Lemmas~\ref{lem:lefteigen} and~\ref{lem:14.17}, we have
\[
  DG(A,\lambda,v)\dot A
  =\binom{\langle \dot A v,v-i_{\C^n}A_{\lambda,v}^{-*}
  P_{v^\perp}A^*v\rangle}{-A_{\lambda,v}^{-1}P_{v^\perp}\dot Av}
  =\binom{v^*-v^*A i_{\C^n}A_{\lambda,v}^{-1}
  P_{v^\perp}}{-A_{\lambda,v}^{-1}P_{v^\perp}}\dot A v=R\dot A v,
\]
where, in the last displayed formula, we are denoting by 
$R\colon \C^n\rightarrow\C\times v^\perp\equiv\C^n$
the linear operator multiplying by $\dot A v$ in the previous formula. 
A standard linear algebra argument then shows that
$\det(DG(A,\lambda,v)DG(A,\lambda,v)^*)=\det(RR^*)=|\det(R)|^2$. 
Now, we can identify 
$$ 
   i_{\C^n}\equiv\binom{0}{\Id_{n-1}},\quad
   P_{v^\perp}\equiv\begin{pmatrix}0&\Id_{n-1}\end{pmatrix}, 
$$ 
which implies that in standard basis we have
\[
  R=\begin{pmatrix}1&*\\
       0&-A_{\lambda,v}^{-1}\end{pmatrix},
\]
thus showing that $|\det(R)|^2=|\det A_{\lambda,v}^{-1}|^2$, and the 
proposition follows.
\eproof

We are now ready to rewrite Theorem~\ref{th:doubleprojectioncoarea} 
in this setting. The following is an important technical result that we will 
use several times.

\begin{proposition}\label{th:doubleEVP}
Let $\chi\colon \V\rightarrow[0,\infty)$ be a measurable mapping. Then,
\begin{equation*}
  \int_{A\in\C^{n\times n}}\sum_{\lambda,v:Av=\lambda v}
  \chi(A,\lambda,v) dA  
  = \int_{(\lambda,v)\in\C\times\mathbb{P}(\C^n)} 
  \int_{A:Av=\lambda v}\chi(A,\lambda,v)\,|\det(A_{\lambda,v})|^2
  \,dA\,d(\lambda,v).
\end{equation*}
Moreover, assume that $\chi$ is unitarily invariant in the sense that
$\chi(A,\lambda,v)=\chi(UAU^*,\lambda,Uv)$ for any unitary matrix
$U$. Fix any a.e.\ continuous mapping 
$\C^n\setminus\{0\}\to {\cal U}_n,\, v\mapsto U_v$ 
such that $U_ve_1=v/\|v\|$ for all $v$. Then,  
for every $\hA\in\C^{n\times n}$ and $\sigma>0$,
\begin{align}\label{eq:citame1}
  & \Exp_{A\sim\mathcal{N}_{\C^{n\times n}}(\hA,\s^2)}
  \left(\sum_{\lambda,v:Av=\lambda v}\chi(A,\lambda,v)\right)\\ 
=\;&\frac{1}{\Gamma(n)\sigma^{2(n-1)}}
\Exp_v \Exp_{(\lambda,w,B)}\left( e^{-\frac{\|\hat{y}_v\|^2}{\sigma^2}}
 \chi\left(\begin{pmatrix}\lambda&w^*\\0&B\end{pmatrix},\lambda,e_1\right)\,
 \big|\det(B-\lambda I_{n-1})\big|^2\right),\notag
\end{align}
where $v\in\Pn$ has the uniform distribution,
$\hat{y}_v=P_{v^\perp}\hA v/\|v\|$, and
$\lambda\sim\mathcal{N}_\C(\hat{\lambda},\sigma^2)$,
$w\sim\mathcal{N}_{\C^{n-1}}(\hat{w},\sigma^2)$, and
$B\sim\mathcal{N}_{\C^{(n-1)\times (n-1)}}(\hat B,\sigma^2)$ are
independent Gaussian random variables centered at 
$$
  \hat{\lambda}:=\frac{\pes{\hA v}{v}}{\|v\|^2};\quad \hat{w}:=
  J^*U_v^*\hat A^* U_v e_1;\quad \hat{B}:=J^*U_v^*\hA U_vJ.  
$$
Here, $J$ is the $n\times(n-1)$ matrix whose columns are
$e_2,\ldots,e_n$ (and, hence, $J^*=[0\;\Id_{n-1}]$ 
is the matrix  of $P_{e_1^\perp}\colon\C^n\to e_1^{\perp}$). 
In particular, if $\hA=0$ then $\hat{y}_v=0$,
$\hat{\lambda}=0$, $\hat{w}=0$ and $\hat{B}=0$, so $v$ can be removed
from the expected value in~\eqref{eq:citame1}.
\end{proposition}

\proof
The first claim of the theorem follows directly from 
Theorem~\ref{th:doubleprojectioncoarea} and 
Proposition~\ref{pro:NJquot}. 

For the second claim, let $I$ be the left-hand side 
of~\eqref{eq:citame1}. Change $\chi(A,\lambda,v)$ to
$(\sigma^2\pi)^{-n^2}\chi(A,\lambda,v)e^{-\frac{\|A-\hA\|_F^2}{\sigma^2}}$
in the first formula to get
\[
  I=\frac{1}{(\sigma^2\pi)^{n^2}}\int_{(\lambda,v)\in\C\times\mathbb{P}(\C^n)} 
  \int_{A:Av=\lambda v}\chi(A,\lambda,v)\,|\det(A_{\lambda,v})|^2
  e^{-\frac{\|A-\hA\|_F^2}{\sigma^2}}\,dA\,d(\lambda,v).
\]
Note that $\{A:Av=\lambda v\}$ can be parametrized by
\[
  (w,B)\mapsto A=U_v\begin{pmatrix}\lambda&w^*\\0&B\end{pmatrix}U_v^*, 
\]
where $w\in\C^{n-1}$, $B\in\C^{(n-1)\times(n-1)}$. This 
parametrization preserves distances; moreover
$|\det(A_{\lambda,v})|=|\det(\lambda I_{n-1}-B)|$ and from the fact
that $\chi$ is unitarily invariant, we have
\[
 I=\frac{1}{(\sigma^2\pi)^{n^2}}\int_{(\lambda,v,w,B)} 
 \chi\left(\begin{pmatrix}\lambda&w^*\\0&B\end{pmatrix},\lambda,e_1\right)\,
 |\det(\lambda I_{n-1}-B)|^2e^{-\frac{\|A-\hA\|_F^2}{\sigma^2}}
 \,d(\lambda,v,w,B),
\]
where $A$ is given by the formula above. Note now that we can write
\begin{align*}
 \|A-\hA\|_F^2
 =&\left\|\begin{pmatrix}\lambda&w^*\\0&B\end{pmatrix}-U_v^*
   \hA U_v\right\|_F^2\\
 =&|\lambda-e_1^*U_v^*\hA U_ve_1|^2
     +\|w^*-e_1^*U_v^*\hA U_vJ\|^2
     +\|B-J^*U_v^*\hA U_vJ\|_F^2\\
  &\quad+\|J^*U_v^*\hA U_ve_1\|^2\\
 =&\left|\lambda-\frac{v^*\hA v}{\|v\|^2}\right|^2
     +\left\|w-J^*U_v^*\hat A^* U_v e_1\right\|^2
     +\|B-J^*U_v^*\hA U_vJ\|_F^2+\|\hat y_v\|^2,
\end{align*}
and the second claim of the theorem follows noting that the volume of 
$\Pn$ is $\pi^{n-1}/\Gamma(n)$.
\eproof

\subsection{The linear solution variety}

It will be useful to consider a geometrical scheme similar to that of
\S\ref{sec:interationV} for the case of solving linear systems:
we consider
\[
   \Vlin=\{(M,v)\in\C^{(n-1)\times n}\times\Pn: \,Mv=0\}.
\]
The linear solution variety $\Vlin$ is a $n(n-1)$-dimensional smooth
submanifold of $\C^{(n-1)\times n}\times\Pn$, and again it inherits
the Riemannian structure of the ambient space 
(cf.~\cite[(17.14)]{Condition}).

The linear solution variety is equipped with two projections
\begin{equation}\label{eq:pilin}
\begin{array}{rccc}\pil\colon &\Vlin&\rightarrow&\C^{(n-1)\times n}\\
&(M,v)&\mapsto&M\end{array}\qquad  
\raisebox{6.5pt}{\mbox{and}}\qquad
\begin{array}{rccc}\pil_2\colon &\Vlin&\rightarrow&\Pn\\&(M,v)&\mapsto&v.
\end{array}
\end{equation}
For $M\in\C^{(n-1)\times n}$, $(\pil)^{-1}(M)$ is a copy of the 
projective linear subspace corresponding to the kernel 
of $M$ in $\Pn$, and for $v\in\Pn$, $(\pil_2)^{-1}(v)$ is a copy of
the linear subspace of $\C^{(n-1)\times n}$ consisting of the matrices
$A$ such that $Av=0$.

We can apply Theorem~\ref{th:doubleprojectioncoarea} for integrating
functions in $\Vlin$ using the projections in~\eqref{eq:pilin}. In
this case, the tangent space to $\Vlin$ at $(M,v)$ can be identified 
with 
\[
  \{(\dot M,\dot v):\dot Mv +M\dot v=0,\;v^*\dot v=0\}
  =\{(\dot M,\dot v):\dot v=\varphi(\dot M)\},\quad 
   \varphi( \dot M)=-M^\dagger\dot M v.
\]
Note that $\varphi$ is a linear mapping defined from 
$\C^{(n-1)\times n}$ to $v^\perp$.
A routine computation shows that, if $\|v\|=1$, then 
$\varphi\varphi^*=M^\dagger (M^{\dagger})^*$. Writing down the 
singular value decomposition of $M$, it follows that 
$\det(\varphi\varphi^*)=\det(MM^*)^{-1}$. From 
Lemma ~\ref{le:detquot} it follows that
\begin{equation}\label{eq:rodelu}
   \frac{\NJ(\pil)(M,v)}{\NJ(\pil_2)(M,v)}=|\det(MM^*)|.
\end{equation}

\begin{proposition}\label{prop:trick2}
Let $\phil\colon \Vlin\to[0,\infty]$ be a measurable unitarily invariant
function in the sense that $\phil(M,v)=\phil(MU^*,Uv)$ for any unitary
matrix $U\in\CU_n$. Then, 
\begin{equation*}
 \E_{M\sim\mathcal{N}_{\C^{(n-1)\times n}}}(\phil(M,\ker(M))
  =\frac{1}{\Gamma{(n)}}
 \E_{B\sim\mathcal{N}_{\C^{(n-1)\times (n-1)}}}
 (\phil((0\;B),e_1)\,|\det(B)|^2).
\end{equation*}
\end{proposition}

\proof
Let $\chi(M)=\phil(M,\ker(M))e^{-\|M\|_F^2}$. 
Theorem~\ref{th:doubleprojectioncoarea} and~\eqref{eq:rodelu} 
imply that 
\[
   \int_{M\in\C^{(n-1)\times n}}\chi(M)\,dM
  \,=\,\int_{v\in\Pn}\int_{M:Mv=0}
   |\det(M)|^2\chi(M)\,dM\,dv.
\]
Now, $|\det(M)|^2\chi(M)=|\det(MU)|^2\chi(MU^*)$ for all $U\in\CU_n$ 
by hypothesis. Hence, by parametrizing $\{M:Mv=0\}$ by 
$\{(0\;B)U_v^*:B\in\C^{(n-1)\times (n-1)}\}$ where $U_v\in\CU_n$ 
is any matrix satisfying $U_ve_1$ (we are assuming $\|v\|=1$), 
we conclude that the inner integral in the right-hand side of the 
formula above does not depend on $v$. We thus have
\[
   \int_{M\in\C^{(n-1)\times n}}\chi(M)\,dM
  \,=\,\Vol(\Pn)\int_{B\in\C^{(n-1)\times (n-1)}}
   |\det(B)|^2\chi((0,B))\,dB.
\]
The proposition follows from the form of the Gaussian density 
(recall~\S\ref{sec:gaussian}) by noting that 
$\Vol(\Pn)=\pi^{n-1}/\Gamma(n)$.
\eproof

Assume now that we are given a measurable nonnegative function
$\alpha\colon \C^{(n-1)\times(n-1)}\to[0,\infty]$. We can produce a
unitarily invariant function defined on $\Vlin$ as follows
\[
  \phil(M,v)=\E_{Q:(M,Q)\sim\A_n}(\alpha(MQ)),
\]
where $\A_n$ is given in Definition~\ref{def:OmegayVarphi} (observe that $Q$ follows the natural, uniform distribution on the Stiefel manifold). Note that
\[
\phil((0\;B),e_1)=\E_{U\in\CU_{n-1}}(\alpha(BU)).
\]
It is a simple exercise
to check that $\phil$ is unitarily invariant in the sense of 
Proposition~\ref{prop:trick2}.
Applying Proposition~\ref{prop:trick2} to $\phil$ then yields 
$$
\E_{M}\; \E_{Q:(M,Q)\in\A_n}(\alpha(MQ))
=\frac{1}{\Gamma{(n)}}\;
\E_{B}\;\E_{U\sim\CU_{n-1}}(\alpha(BU)|\det(B)|^2).  
$$ 
Finally, using Fubini's theorem, we can interchange the integration
order in the right-hand term, and then note that the isometry
$B\mapsto BU$ preserves the value of the integral inside. 
We obtain the following corollary.

\begin{corollary}\label{cor:trick2}
Let $\alpha\colon \C^{(n-1)\times(n-1)}\to[0,\infty]$ 
be an a.e. continuous function. 
Then,
\begin{equation}\tag*{\qed}
 \E_{M}\; \E_{Q:(M,Q)\in\A_n}(\alpha(MQ))=\frac{1}{\Gamma{(n)}}\;
 \E_{B\sim\mathcal{N}_{\C^{(n-1)\times (n-1)}}}(\alpha(B)|\det(B)|^2).
\end{equation}
\end{corollary}

\section{Proof of Theorem~\ref{th:mu2-bound}}\label{se:muave} 

We begin with the following result.

\begin{proposition}\label{prop:smooth} The following inequality holds
for every $\hA\in\C^{n\times n}$ and $\sigma>0$.
$$
  \Exp_{A\sim \mathcal{N}_{\C^{m\times m}}(\hA,\sigma^2)} 
  \big(\|A^{-1}\|_F^2\,|\det(A)|^2\big)\leq \frac{m}{\sigma^2}\, 
  \Exp_{A\sim \mathcal{N}_{\C^{m\times m}}(\hA,\sigma^2)} (|\det(A)|^2).
$$
Furthermore, equality holds if and only if $\hA=0$. In particular, 
$$
 \Exp_{A\sim\mathcal{N}_{\C^{m\times m}}(0,\sigma^2)} 
 \big(\|A^{-1}\|_F^2\,|\det(A)|^2\big)=m!\,m \sigma^{2m-2}.
$$
\end{proposition}

\proof
Expanding the determinant of $A$ by the $k$th column we have
$$
 \det(A)=\sum_{j=1}^m (-1)^{j+k}{a}_{j,k} \det{A}^{j,k},
$$
where $A^{j,k}$ denotes the matrix that results from the matrix $A$ by 
removing the $j$th row and $k$th column. Hence,
$$
 |\det(A)|^2=\det(A)\overline{\det(A)}=\sum_{j,j'=1}^m (-1)^{j+j'+2k}{a}_{j,k}\,
 \overline{{a}_{j',k}}\, \det{A}^{j,k}\,\overline{\det{A}^{j',k}}.
$$
Observe that the random variables ${a}_{j,k}$ and ${a}_{j',k}$ are 
independent of $\det{A}^{j,k}$ and $\det{A}^{j',k}$. Then, 
$$
  \Exp_{A\sim \mathcal{N}_{\C^{m\times m}}(\hA,\sigma^2)}|\det(A)|^2
  =\sum_{j,j'=1}^m (-1)^{j+j'+2k}\Exp({a}_{j,k}\,\overline{{a}_{j',k}} )\,
  \Exp(\det{A}^{j,k}\,\overline{\det{A}^{j',k}}).
$$
Now observe that
\[ 
  \Exp({a}_{j,k}\,\overline{{a}_{j',k}} )= \left\{ \begin{array}{ll}
          \widehat{a}_{j,k}\,\overline{\widehat{a}_{j',k}} & \mbox{if $j\neq j'$};\\
        \sigma^2+|\widehat{a}_{j,k}|^2 & \mbox{otherwise}.\end{array} \right. 
\] 
We conclude that for $k=1,\ldots,m$, 
\begin{equation}\label{saconstant2}
    \Exp_{A\sim \mathcal{N}_{\C^{m\times m}}(\hA,\sigma^2)}|\det(A)|^2
  =\Exp|\det([A;k;\hA_k])|^2+\sigma^2\sum_{j=1}^m \Exp |\det{A}^{j,k}|^2,
\end{equation}
where $[A;k;\hA_k]$ is the matrix formed by replacing the
(random) $k$th column of $A$ by the (deterministic) $k$th column of
$\hA$.
 
On the other hand, from a direct application of Cramer's Rule 
and~\eqref{saconstant2}, we deduce that 
\begin{eqnarray*}
 \sigma^2\Exp_{A\sim \mathcal{N}_{\C^{m\times m}}(\hA,\sigma^2)}
 \|A^{-1}\|^2_F\,|\det(A)|^2&=&\sigma^2\sum_{j,k=1}^m 
 \Exp|\det{A^{j,k}}|^2\\
&=& m\, 
 \Exp|\det(A)|^2- \sum_{k=1}^m\Exp|\det([A;k;\hA_k])|^2,
\end{eqnarray*}
and the first claim of the proposition follows. Moreover, when
$\hA =0$, the last term in the sum above is zero. We 
leave the proof of the converse to the reader.  
Using (\ref{saconstant2}) and the fact that the
matrices $A^{j,k}$ are 
$\mcN_{\C^{(m-1)\times(m-1)}}(0,\sigma^2)$-distributed, 
one can prove working by induction the equality 
$$
  \Exp_{A\sim \mcN_{\C^{m\times m}}(0,\sigma^2)}
  |\det(A)|^2=\sigma^{2m}m!.
$$
The second claim of the proposition follows.
\eproof

\begin{corollary}\label{cor:exactpinvM}
\[
  \E_{M\sim\mcN_{\C^{(n-1)\times n}}}(\|M^\dagger\|_F^2)=n-1.
\]
\end{corollary}

\proof
From Proposition~\ref{prop:trick2} with
$\phil(M,\zeta)=\|M^\dagger\|_F^2$, we have
\[
  \E_{M\sim\mcN_{\C^{(n-1)\times n}}}(\|M^\dagger\|_F^2)
  =\frac{1}{\Gamma{(n)}}
  \E_{B\sim\mcN_{\C^{(n-1)\times (n-1)}}}(\|B^{-1}\|_F^2|\det(B)|^2)
  \underset{\text{Prop.~\ref{prop:smooth}}}{=}n-1
\]
as claimed.
\eproof

\begin{corollary}\label{cor:inter} 
For any $\hB\in\C^{(n-1)\times (n-1)}$, $\sigma>0$, and $\lambda\in\C$, 
we have 
$$
 \Exp_{B}\big(\big\|(B-\lambda I_{n-1})^{-1}\big\|_F^2\,
 |\det(B-\lambda I_{n-1})|^2\big)\leq \frac{n-1}{\sigma^2}
 \Exp_{B}\big(|\det(B-\lambda I_{n-1})|^2\big),
$$
where $B\sim\mathcal{N}_{\C^{(n-1)\times (n-1)}}(\hB,\sigma^2)$.
\end{corollary}

\proof
Note that  
\[
  \Exp_{B}(\|(B-\lambda I_{n-1})^{-1}\|_F^2\,|\det(B-\lambda I_{n-1})|^2)
 = \Exp_{C\sim\mathcal{N}_{\C^{(n-1)\times (n-1)}}
   (\widehat{C},\sigma^2)}(\|C^{-1}\|_F^2\,|\det C|^2),
\]
where $\widehat C=\widehat B-\lambda I_{n-1}$. The proof readily 
follows from Proposition~\ref{prop:smooth}.
\eproof
\medskip

\proofof{Theorem~\ref{th:mu2-bound}} 
Fix any a.e. continuous mapping $v\mapsto U_v$ such that for 
$v\in\Pn$, $U_v$ is a unitary matrix with $Ue_0=v/\|v\|$. 
From~\eqref{eq:citame1} applied to 
$\chi(A,\lambda,v)=\frac{1}{n}\mu_F^2(A,\lambda,v)/\|A\|_F^2
=\frac1n\|A_{\lambda,v}^{-1}\|_F^2$ we have
\begin{eqnarray}\label{eq:citame}
& & \Exp_{A\sim\mathcal{N}_{\C^{m\times m}}(\hA,\sigma^2)}
  \left(\frac{\muFa^{2}(A)}{\|A\|_F^{2}}\right)\\
&=&
  \frac{1}{n\Gamma(n)\sigma^{2(n-1)}}\Exp_v\Exp_{( \lambda,w,B)}
  (e^{-\frac{\|\hat{y}_v\|^2}{\sigma^2}}\|(B-\lambda I_{n-1})^{-1}\|_F^{2}\,
  |\det(B-\lambda I_{n-1})|^2)\nonumber,
\end{eqnarray}
where $\hat{y}_v=P_{v^\perp}\hA v/\|v\|$, $v\in\Pn$ has the
uniform distribution, and
$\lambda\sim\mathcal{N}_\C(\hat{\lambda},\sigma^2)$,
$w\sim\mathcal{N}_{\C^{n-1}}(\hat{w},\sigma^2)$, and
$B\sim\mathcal{N}_{\C^{(n-1)\times (n-1)}}(\hat{B},\sigma^2)$ for some
$\hat\lambda,\hat{w},\hat{B}$ which depend uniquely on $\hA$ and $v$.

From~\eqref{eq:citame} and Corollary~\ref{cor:inter} we have 
\begin{align*}
  &\Exp_{A\sim\mathcal{N}_{\C^{n\times n}}(\hA,\sigma^2)}
  \left(\frac{\muFa^2(A)}{\|A\|_F^2}\right)
%   \\
% \leq\; &
%   \frac{1}{\Gamma(n+1)}\Exp_v\Exp_{( \lambda,w,B)}(e^{-\|\hat{y}_v\|^2}
%   \|(B-\lambda I_{n-1})^{-1}\|_F^{2}\,|\det(B-\lambda I_{n-1})|^2)
  \\
\leq\; &
  \frac{n-1}{\sigma^2\sigma^{2(n-1)}\Gamma(n+1)}
  \Exp_v\Exp_{( \lambda,w,B)}(e^{-\frac{\|\hat{y}_v\|^2}{\sigma^2}}
  |\det(B-\lambda I_{n-1})|^2).
\end{align*}
Now, if we apply again~\eqref{eq:citame1} to the constant function 
$\chi\equiv1/n$ we get
\[
  1=\Exp_{A\sim\mathcal{N}_{\C^{m\times m}}(\hA,\sigma^2)}
  \left(1\right)=
  \frac{1}{\Gamma(n+1)\sigma^{2(n-1)}}\Exp_v\Exp_{( \lambda,w,B)}
  (e^{-\frac{\|\hat{y}_v\|^2}{\sigma^2}}
  |\det(B-\lambda I_{n-1})|^2),
\]
and we have thus proved that
\[
  \Exp_{A\sim\mathcal{N}_{\C^{n\times n}}(\hA,\sigma^2)}
  \left(\frac{\muFa^2(A)}{\|A\|_F^2}\right)\leq\frac{n-1}{\sigma^2}
  \leq\frac{n}{\sigma^2},
\]
as claimed.

This proves the first part of Theorem~\ref{th:mu2-bound}. 
For the second part of the theorem,  let
\[
   I=\Exp_{A\sim\msU(\SS)}\muFa^2(A)
\]
be the quantity we want to compute. From the first part of the theorem 
(with $\hA=0$ and $\sigma=1$) we have
\begin{equation}\label{eq:201}
\frac{1}{\pi^{n^2}}\int_{A\in\mnc}\frac{\muFa^2(A)}{\|A\|_F^2}e^{-\|A\|_F^2}\,dA
\leq n.
\end{equation}
On the other hand,
\begin{equation}\label{eq:202}
\frac{1}{\pi^{n^2}}\int_{A\in\mnc}\frac{\muFa^2(A)}{\|A\|_F^2}e^{-\|A\|_F^2}\,dA 
=\frac{1}{\pi^{n^2}}\int_0^\infty\frac{e^{-\rho^2}}{\rho^2}\int_{A:\|A\|_F=\rho}
\muFa^2(A)\,dA\,d\rho.
\end{equation}
Now, because $\muFa(A)$ is invariant under multiplication of $A$ 
by nonzero complex numbers, denoting $\nu_\rho=\Vol(A:\|A\|_F=\rho)$, 
we have
\begin{equation}\label{eq:203}
 \frac{1}{\nu_\rho}\int_{A:\|A\|_F=\rho}\muFa^2(A)\,dA
=I,\quad 0<\rho<\infty.
\end{equation}
We deduce from~(\ref{eq:201}--\ref{eq:203}) that
\[
  \frac{I}{\pi^{n^2}}\int_0^\infty
  \frac{\nu_\rho e^{-\rho^2}}{\rho^2}\,d\rho\leq n.
\]
Note now that
\[
 \nu_\rho=\frac{2\pi^{n^2}}{\Gamma(n^2)}\rho^{2n^2-1}
\]
to conclude that 
\[
  I\leq  \frac{n\Gamma(n^2)}{2\int_0^\infty\rho^{2n^2-3}e^{-\rho^2}\,d\rho}
  =\frac{n\Gamma(n^2)}{\Gamma(n^2-1)}=n(n^2-1)\leq n^3.
\]
The theorem follows. 
\eproof

\section{Proof of Propositions~\ref{prop:fixA_0} 
and~\ref{prop:fixA_0_all}}\label{sec:main_proof}

\subsection{Proof of Proposition~\ref{prop:fixA_0}}

First note that \EC\ starts by 
normalizing the input, so from~\eqref{eq:truncating} we can assume that 
$\|A_0\|_F=1$ and $A\sim\mcN_{\C^{n\times n},T}(0,1)$  
where $T=\sqrt{2} n$. From 
Remark~\ref{rem:codimension}, for integration purposes 
we can also assume that 
$(\mathcal{L}_{A_0,A}\setminus\{A_0\})\cap\Sigma=\emptyset$. 
Corollary~\ref{corollary:daigual} with $q=1$ implies that
\begin{align}\label{eq:step1}
  \aviter(A_0,\lambda_0,v_0)=\;& \E_{A\sim\mcN_{\C^{n\times n},T}}
   K(A,A_0,\lambda_0,v_0)\notag\\
  \le\;& 2+\sqrt{6}C\E_{A\sim\mcN_{\C^{n\times n},T}}
            \|A\|_F\int_{t_1}^1\
            \frac{\mu^2(A_t,\lambda_t,v_t)}{\|A_t\|_F^2}\,dt\\
  \le\;& 2+\sqrt{6}C\E_{A\sim\mcN_{\C^{n\times n},T}}
         \|A\|_F\int_{t_1}^1\;\sum_{j=1}^n
         \frac{\mu^2(A_t,\lambda_t^{(j)},v_t^{(j)})}{\|A_t\|_F^2}\,dt,\notag
\end{align}
where $A_t=(1-t)A_0+tA$ and the pairs $(\lambda_t^{(j)},v_t^{(j)})$ 
are defined by continuation for all the 
eigenpairs of $A_t$. Note that we need to write $2+$ instead of $1+$ as in Corollary \ref{corollary:daigual} because we need to bound the smallest integer which is greater than the integral, rather than the integral itself. 

We therefore have (use $T=\sqrt{2}n$) 
\begin{eqnarray}\label{eq:step2}
  \aviter(A_0,\lambda_0,v_0)
&\leq&  2+\sqrt{6}CnT\E_{A\sim \mcN_{\C^{n\times n},T}}
  \left( \int_{t_1}^1\frac{\mu_{\av}^2(A_t)}{\|A_t\|_F^2}\,dt\right)\notag\\
  &\underset{\text{\eqref{eq:truncated}}}{\leq}&
  2+\sqrt{48}Cn^2\E_{A\sim \mcN_{\C^{n\times n}}}
 \left( \int_{t_1}^1\frac{\mu_{\av}^2(A_t)}{\|A_t\|_F^2}\,dt\right).
\end{eqnarray}
In order to bound the last term in the previous expression, we interchange 
the order of integration, 
\[
   \E_{A\sim \mcN_{\C^{n\times n}}}
   \left(\int_{t_1}^1\frac{\mu_{\av}^2(A_t)}{\|A_t\|_F^2}\,dt\right)
   =\int_{t_1}^1\E_{A\sim \mcN_{\C^{n\times n}}}
   \left(\frac{\mu_{\av}^2(A_t)}{\|A_t\|_F^2}\right)\,dt.
\]
Now, for fixed $t$, if $A\sim \mcN_{\C^{n\times n}}$ 
then $A_t=(1-t)A_0+tA$ satisfies 
$A_t\sim \mathcal{N}_{\C^{n\times n}}((1-t)A_0,t^2)$ 
and from Theorem~\ref{th:mu2-bound} we have
\[
   \E_{A\sim \mcN_{\C^{n\times n}}}
  \left(\frac{\mu_{\av}^2(A_t)}{\|A_t\|_F^2}\right)
   \leq\frac{n}{t^2},
\]
which implies
\begin{equation}\label{eq:step3}
   \E_{A\sim \mcN_{\C^{n\times n}}}
   \left(\int_{t_1}^1\frac{\mu_{\av}^2(A_t)}{\|A_t\|_F^2}\,dt\right)
   \leq \int_{t_1}^1\frac{n}{t^2}\,dt\leq\frac{n}{t_1}.
\end{equation}
We are thus left with the task of evaluating $t_1$. But we have (recall from 
Corollary~\ref{corollary:daigual}) that 
$$
  t_1=\frac{1}{\|A_1\|_F(\sin\alpha \cot s_1-\cos\alpha)+1}
      \geq \frac{1}{T(\cot s_1+1)+1}.
$$
We now use that $s_1={\sf Choose\_step}(A_0,\dot A,\lambda_0,v_0)$ 
(the length of the first step in the execution of \EC) and this is 
at least $\frac{R}{\mu^2(A_0,\lambda_0,v_0)}$ by 
Proposition~\ref{prop:step-size}. Hence $\cot s_1 \leq \frac1{s_1} \leq 
\frac {\mu^2(A_0,\lambda_0,v_0)}{R}$ 
and it follows that 
$$
   t_1\geq \Omega\left(\frac{1}{n \mu^2(A_0,\lambda_0,v_0)}\right).
$$
Putting together this bound and inequalities~\eqref{eq:step2} 
and~\eqref{eq:step3} we deduce the claimed bound 
for $\aviter(A_0,\lambda_0,v_0)$.

We next prove the smoothed analysis bounds. Reasoning as 
in~\eqref{eq:step1} we see that the smoothed number of iterations 
$\siter(A_0,\lambda_0,v_0,\sigma)$ is bounded by 
\[
     2+\sqrt{6}C\,\sup_{\hA\in\SS}\E_{A\sim\mcN_{\C^{n\times n},T}
            (\hA,\sigma^2)}\sum_{j=1}^n \|A\|_F
    \int_{t_1}^1\frac{\mu^2(A_t,\lambda_t^{(j)},v_t^{(j)})}{\|A_t\|_F^2}\,dt.
\]
The rest of the argument is almost exactly as above, the only
difference being the bound $\|A\|_F\leq T+\|\hA\|_F = \sqrt{2}n+1$. 
\eproof

\subsection{Proof of Proposition~\ref{prop:fixA_0_all}}

We are now following all the $n$ paths (each starting with a different
eigenpair of $A_0$). Applying Corollary ~\ref{corollary:daigual} with 
$q=1$ to each of them we obtain 
\begin{align}\label{eq:step11}
  \aviter(A_0)=& \E_{A\sim\mcN_{\C^{n\times n},T}}
   \sum_{j=1}^n K(A,A_0,\lambda^{(j)},v^{(j)})\notag\\
 \le &\; 2n+\sqrt{6}C\E_{A\sim\mcN_{\C^{n\times n},T}}
            \|A\|_F \sum_{j=1}^n \int_{t^{(j)}_1}^1\
            \frac{\mu^2(A_t,\lambda_t^{(j)},v_t^{(j)})}{\|A_t\|_F^2}\,dt\\
 \le &\; 2n+\sqrt{6}C\E_{A\sim\mcN_{\C^{n\times n},T}}
            \|A\|_F \int_{t^{*}_1}^1 \sum_{j=1}^n 
            \frac{\mu^2(A_t,\lambda_t^{(j)},v_t^{(j)})}{\|A_t\|_F^2}\,dt,\notag
\end{align}
where $t_1^* =\min\{t^{(1)}_1, \ldots,t^{(n)}_1\}$. We can now reason 
as in the preceding proof to deduce that 
\[
    \aviter(A_0)=\Oh\Big(\frac{n^3}{t_1^*}\Big)
\]
as well as the fact that, for $j=1,\ldots,n$, 
$t_1^{(j)}\geq \Omega\Big(\frac{1}{n\mu^2(A_0,\lambda^{(j)},v^{(j)})}\Big)$. 
It follows from these bounds that 
\[
   t_1^*\geq \Omega\left(\frac{1}{n\mum^2(A_0)}\right).
\]
The rest of the proof follows as in the preceding proposition.
\eproof

\section{Proof of Theorem~\ref{th:randomhomotopy}}
\label{sec:proofofrandomisrandom}

We begin with an auxiliary result. For simplicity, in what follows we write 
$\SS:=\SS(\Cnn)$. We also consider the manifold
\[
 \mathcal{S}=\{(A,\dot A)\in\SS\times\SS:\dot A\in T_A\SS\}
\]
and denote by $\msU(\mcS)$ the uniform distribution on it. 
Given any measurable mapping 
$\phi\colon \mathcal{S}\rightarrow[0,\infty]$, let
\begin{equation}\label{eq:defI}
I_\phi:=\E_{A_0,A\sim\msU(\mcS)}\left(\int_0^{d_\SS(A_0,A)}
     \phi(A_s,\dot A_s)\,ds\right),
\end{equation}
where, as usual, the $A_s$ are such that 
$\{A_s:0\leq s\leq d_\SS(A_0,A)\}=\mathcal{L}_{A_0,A}$, and 
$\dot{A}_s$ is the unit tangent vector (in the direction of the 
parametrization) to $\mathcal{L}_{A_0,A}$ at $A_s$. 
\begin{lemma}\label{lemn:sphere}
 For any measurable mapping 
$\phi\colon \SS\times\SS\rightarrow[0,\infty]$ we have
\begin{equation}\label{eq:Iphi}
  I_\phi=\frac{\pi}{2}\E_{(A,\dot A)\sim\msU(\mcS)}
  \phi(A,\dot A).
\end{equation}

\end{lemma}
\proof
We consider the manifold $\mathcal{R}=  \{(A_0,A,s)\in \SS\times\SS\times(0,\pi):s<d_\SS(A_0,A)\}$ and let
\[
\begin{matrix}
\psi\colon &\mathcal{R}&\rightarrow&\mathcal{S}\\
&(A_0,A,s)&\mapsto&(A_s,\dot A_s).
\end{matrix}
\]
We can then write $I_\phi=\Vol(\SS)^{-2}\int_{\mathcal{R}}\phi\circ \psi$. 
Applying the coarea formula (Theorem~\ref{pro:coarea}) yields
\begin{equation}\label{eq:auxIphi}
  I_\phi=\frac{1}{\Vol(\SS)^2}\int_{(A,\dot A)\in\mathcal{S}}
  \phi(A,\dot A)q(A,\dot A)\,d(A,\dot A),
\end{equation}
where
\[
q(A,\dot A)=\int_{(A_0,A,s)\in\psi^{-1}(A,\dot A)}
  \frac{1}{\NJ(\psi)(A,A_0,s)}\,d(A_0,A,s).
\]
Our goal now is to prove that $q(A,\dot A)$ is a constant (independent
of $A,\dot A$). It shall be useful to consider the two diagonal
matrices $\Delta=E_{11}$ and $\dot \Delta=E_{22}$ where $E_{ij}$
denotes the standard basis of $\Cnn$. Note that
$(\Delta,\dot\Delta)\in\mathcal{S}$.

We now fix $(A,\dot A)\in\mathcal{S}$. Let
$\sigma\colon \mnc\rightarrow\mnc$ 
be an isometric change of basis such that
$\sigma(A)=\Delta$ and $\sigma(\dot A)=\dot \Delta$. Denoting
$\sigma_\mathcal{S}=\sigma\times\sigma$ and
$\sigma_{\mathcal{R}}=\sigma\times\sigma\times I_{\R}$, where $I_{\R}$
is the identity mapping in $\R$, it is easy to check by writing down
the formula for $A_s$ that
$\psi\circ\sigma_{\mathcal{R}}=\sigma_{\mathcal{S}}\circ\psi$. We
are under the hypothesis of~\cite[Lemma~4, p.~244]{BlCuShSm98} 
(which holds for surjective maps in general, not only for 
projections), proving that
\begin{equation}\label{eq:Njequals}
   \NJ(\psi)(\sigma_\mathcal{R}(C,C_0,s))
  =\NJ(\psi)(C,C_0,s),\quad \forall\, (C,C_0,s)\in\mathcal{R}.
\end{equation}
Moreover, the mapping $\sigma_\mathcal{R}\mid_{\psi^{-1}(A,\dot{A})}$ 
is an isometry from $\psi^{-1}(A,\dot A)$ to
$\psi^{-1}(\Delta,\dot \Delta)$, which from the change of variables
theorem implies
\begin{align*}
 q(A,\dot A)=&\int_{(B_0,B,s)\in\psi^{-1}(A,\dot A)}
  \frac{1}{\NJ(\psi)(\sigma_\mathcal{R}^{-1}(B,B_0,s))}\,d(B_0,B,s)\\
  \underset{\eqref{eq:Njequals}}{=}
 &\int_{(B_0,B,s)\in\psi^{-1}(\Delta,\dot \Delta)}
  \frac{1}{\NJ(\psi)(B,B_0,s)}\,d(B_0,B,s)=q(\Delta,\dot\Delta),
\end{align*}
which proves that $q(A,\dot A)$ is equal to some constant $\hat{C}$.

Since this holds for any measurable function $\phi$, we can take the latter to be constant with value $1$ 
to derive the value of $\hat{C}$. Then \eqref{eq:auxIphi} becomes 
$$
 I_\phi= \frac{1}{\Vol(\SS)^2}\int_{A,\dot A\in\mathcal{S}}
  \hat{C}\,d(A,\dot A)
 = \frac{\Vol(\mathcal{S})}{\Vol(\SS)^2}\hat{C}. 
$$
And it follows from~\eqref{eq:defI} (always with $\phi\equiv 1$) that 
$$
I_\phi=\E_{A,A_0\sim \msU(\SS)}d_\SS(A_0,A).
$$
Hence,
\[
  \E_{A,A_0\sim \msU(\SS)} d_\SS(A,A_0)
  = \hat{C}\frac{\Vol(\mathcal{S})}{\Vol(\SS)^2}.
\]
Note that the change of variables $A_0\mapsto -A_0$ does not change
the expected value in this last formula. Moreover,
$d_\mathbb{S}(A,A_0)+d_\mathbb{S}(A,-A_0)=\pi$ for all
$A,A_0\in\mathbb{S}$. Thus,
\begin{eqnarray*}
\frac{\Vol(\mathcal{S})}{\Vol(\SS)^2}\,2\hat{C}
  &=&\E_{A,A_0\sim\msU(\SS)} d_\SS(A,A_0)
  +\E_{A,A_0\sim\msU(\SS)} d_\SS(A,-A_0)\\
  &=&\E_{A,A_0\sim\msU(\SS)}\left(d_\SS(A,A_0)+d_\SS(A,-A_0)\right)
     \;=\;\pi.
\end{eqnarray*}
proving that 
\[
  \hat{C}=\frac{\Vol(\SS)^2\pi}{2\Vol(\mathcal{S})},
\]
From \eqref{eq:auxIphi} we conclude that for any measurable 
nonnegative function $\phi$ we have
\[
  I_\phi=\frac{\pi}{2}\E_{(A,\dot A)\sim\msU(\mcS)}
  \phi(A,\dot A),
\]
as wanted.
\eproof

\proofof{Theorem~\ref{th:randomhomotopy}}
Consider the measurable function $\phi\colon\mathcal{S}\to [0,\infty]$ defined by 
$$
 \phi(A,\dot{A}) = \frac{1}{n}\sum_{(\lambda,v):Av=\lambda v}
  \mu(A,\lambda,v)\big\|(\dot A,\dot\lambda,\dot v\big)\|,
$$
for $A\in\SS$ and $\dot{A}\in T_A\SS$ such that $A\not\in\Sigma$, where $\dot \lambda,\dot v$ are the functions of $(A,\dot{A})$ 
and $(\lambda,v)$ given in Lemma~\ref{lem:14.17}. (If $A\in\Sigma$ we set $\phi(A,\dot{A})=\infty$.)

From Theorem~\ref{thm:main_path_following}, denoting 
$I=I_\phi$, we have 
(for some constant $C>1$) that  
\begin{equation}\label{eq:ExpI}
  \Exp_{A,A_0\sim\mcN_{\C^{n\times n}}}
  \left(\frac{1}{n}\sum_{\lambda_0,v_0:A_0v_0=\lambda_0v_0}
 K(A,A_0,\lambda_0,v_0)\right)\leq C\,I.
\end{equation}
Note that the left--hand side of \eqref{eq:ExpI} is the quantity to be bounded in Theorem~\ref{th:randomhomotopy}. 
It is therefore enough for us show that $I\leq 4n^2$.  To do so, write 
$\mathcal{S}_A:=\{A':(A,A')\in\mathcal{S}\}\subseteq T_A\SS$, for $A\in\SS$. 
First note that $\cal S_A$ is just the unit sphere in $T_A\SS$, so it has a natural volume form inherited from
  $\Cnn$ and $\vol({\cal S}_A)$ is independent of $A$. Moreover, the Normal Jacobian of the projection $\mathcal{S}\rightarrow\SS$, $(A,\dot A)\mapsto A$, 
is constant and equal to $1/\sqrt{2}$ (this is easy to prove: check that for $\dot A\in T_A\SS$ the pair of the form $(\dot A,\dot A')$ in $T_{(A,A')}\mathcal{S}$ 
which is orthogonal to the kernel of the derivative of the projection is $(\dot A,-\R e\langle A',\dot A\rangle A)$, then note that the vectors of that form obtained 
from any o.g. basin of $T_A\SS$ whose first element is $A'$ are orthogonal and only one of them, $(A',-A')$, changes its norm by $\sqrt{2}$), 
so we have $\vol({\cal S}) = \sqrt{2}\vol(\SS) \vol({\cal S}_A)$.
From Lemma \ref{lemn:sphere} and Theorem \ref{pro:coarea}, we then have 
\begin{align}\label{eq:anotherI}
 I=\;&\frac{\pi}{2}\E_{(A,\dot A)\in\mathcal{S}}
  \phi(A,\dot A)\nonumber\\
 =\;&\frac{\pi}{\sqrt{2}}\E_{A\in\SS}\left(\frac{1}{n}\sum_{(\lambda,v):Av
  =\lambda v}\mu(A,\lambda,v)\E_{\dot A\in \mathcal{S}_A}
 (\big\|(\dot A,\dot\lambda,\dot v)\big\|)\right).
\end{align}
In order to estimate this last quantity we shall use Lemma~\ref{lem:averageperturbation} below. Note first that from Cauchy-Schwartz,
\begin{eqnarray*}
  \E_{\dot A\in \mathcal{S}_A}
 (\big\|(\dot A,\dot\lambda,\dot v)\big\|)
  &\leq &\left(1+\E_{\dot A\in \mathcal{S}_A}
 (|\dot\lambda|^2+\|\dot v\|^2)\right)^{1/2}\\
 &\underset{\text{Lemma \ref{lem:averageperturbation}}}{\leq}& 
 \left(1+\frac{1}{n^2-\frac12}
 \left(1+2\mu_{F}(A,\lambda,v)^2\right)\right)^{1/2}\\
 &\underset{n\geq2}{\leq}& 
\frac{1}{\sqrt{7}}\left(9+\frac{16}{n^2}\mu_{F}(A,\lambda,v)^2\right)^{1/2}.
\end{eqnarray*}
It is a simple exercise to check that for positive $x\in\R$ we have 
$x(9+16x^2/n^2)^{1/2}\leq 9n/8+4x^2/n$. Using this inequality 
we get from the equations above:
\[
   I\leq \frac{\pi}{\sqrt{14}}\E_{A\in\SS}\left(
   \frac1n\sum_{(\lambda,v):Av=\lambda v}
  \left(\frac{9n}{8}+\frac{4}{n}\mu_{F}(A,\lambda,v)^2\right)\right).
\]
We next use Theorem~\ref{th:mu2-bound} (averaging over $A\in\SS$) and bound this last quantity by
\begin{equation}\label{eq:nonsharp}
    I\leq \frac{\pi}{\sqrt{14}}\left(\frac{9n}{8}+\frac{4}{n}n^3\right)
   = \frac{\pi}{\sqrt{14}}\left(\frac{9n}{8}+4n^2\right)\underset{n\geq2}{\leq} 4n^2,
\end{equation}
which finishes the proof.
\eproof

We have used the following technical lemma which is in the spirit 
of~\cite{Armentano:10}. Note that as pointed out in the proof of Theorem~\ref{th:randomhomotopy}, 
$\cal S_A$ is just the unit sphere in the tangent space $T_A\SS$ and thus has a natural measure inherited from the inner product in $\Cnn$.

\begin{lemma}\label{lem:averageperturbation}
Let $(A,\lambda,v)\in\V$. Define  
$\mathcal{S}_A:=\{A':(A,A')\in\mathcal{S}\}\subseteq T_A\SS$ as in the proof of Theorem~\ref{th:randomhomotopy} (that is, $\mathcal{S}_A$ is the unit sphere of $T_A\SS$)
and, for $\dot A\in\mathcal{S}_A$, let $\dot\lambda,\dot v$ 
be as in Lemma~\ref{lem:14.17}. Then, 
 \[
  \E_{\dot A\in \mathcal{S}_A}(|\dot\lambda|^2)
  =\frac{1}{n^2-\frac12}
  \left(\mu_\lambda(A,\lambda,v)^2-\frac{|\lambda|^2}{2}\right),
\]
and
\[
  \E_{\dot A\in \mathcal{S}_A}(
    \|\dot v\|^2)
  =\frac{1}{n^2-\frac12}\|A_{\lambda,v}^{-1}\|_F^2.
\]
In particular, from Proposition \ref{pro:cn-eigenv},
\[
 \E_{\dot A\in \mathcal{S}_A}(|\dot\lambda|^2+\|\dot v\|^2)\leq \frac{1}{n^2-\frac12}\left(1+2\|A_{\lambda,v}^{-1}\|_F^2\right) 
    =\frac{1}{n^2-\frac12}\left(1+2\mu_F(A,\lambda,v)^2\right).
\]
\end{lemma}

\proof
Note that $T_A\SS$ coincides with the (real) orthogonal complement, with respect to $\Re\langle\cdot,\cdot\rangle_F$, of $A\in\C^{n\times n}$. 
Thus $\dim_\R(T_A\SS)=2n^2-1$. On this space we consider the push-forward measure of the standard Gaussian distribution on $\C^{n\times n}$ 
by the orthogonal  projection $\C^{n\times n}\rightarrow T_A\SS$.

Since $T_A\SS$ may be split in the (real) orthogonal decomposition $ \R\sqrt{-1}A\oplus A^\perp$, where $\R\sqrt{-1}A$ is 
the linear real subspace generated by $\sqrt{-1}A$, in particular we conclude that the Gaussian distribution on $T_A\SS$ 
coincides with the distribution $t\sqrt{-1}A+\dot B\in T_A\SS$ where $t\sim \mathcal{N}(0,\frac12)$ and $\dot B\sim \mcN_{A^\perp}$ are independent.

\vspace{7pt}

\noindent \underline{Claim I:} Given a linear operator $L:T_A\SS\to \C^k$, we have
$$
\E_{\dot A\in T_A\SS}(\|L(\dot A)\|^2)={\left(n^2-\frac12\right)}\E_{\dot A\in\mathcal{S}_A} (\|L(\dot A)\|^2);
$$
The claim follows integrating in polar coordinates. More precisely,
\begin{align*}
\E_{\dot A\in T_A\SS}(\|L(\dot A)\|^2) &
=\frac{1}{\pi^{n^2-\frac12}}\int_{\dot A\in T_A\SS}\|L(\dot A)\|^2\,e^{-\|\dot A\|^2}\,d\dot A\\
&=\frac{1}{\pi^{n^2-\frac12}}\int_0^{+\infty}\rho^{2n^2}e^{-\rho^2}\,d\rho\cdot\int_{\dot A\in\mathcal{S}_A}\|L(\dot A)\|^2\,d\dot A\\
&= {\left(n^2-\frac12\right)}\E_{\dot A\in\mathcal{S}_A} (\|L(\dot A)\|^2),
\end{align*}
where we have used that $\int_0^{+\infty}\rho^{2n^2}e^{-\rho^2}\,d \rho=\frac12\Gamma(n^2+\frac12)$, and $\Vol(\mathcal{S}_A)=2\pi^{n^2-\frac12}/{\Gamma(n^2-\frac12)}$.

\vspace{5pt}

\vspace{5pt}

\noindent \underline{Claim II:} The push-forward measure of the Gaussian distribution on $A^\perp$ by the map $f:A^\perp\to v^\perp$,  $\dot A\mapsto P_{v^\perp}(\dot Av),$ is the standard Gaussian on $v^\perp$.

\vspace{5pt}
Note that for all $B\in\C^{n\times n}$, we have $\langle uv^*,B\rangle_F=\mbox{tr}(B^*uv^*)=v^*B^*u=\langle u, Bv\rangle$.
Then, the set $F=\{\dot wv^*:\,\dot w \in v^\perp\}$ is a linear subspace of $A^\perp$, and the kernel of $f$ is the Hermitian complement of $F$ as subset of $A^\perp$. Since $f\mid_F:F\to v^\perp$ is a linear isometry, the claim follows from the characterization of the standard Gaussian distribution.

\vspace{5pt}

\noindent \underline{Claim III:}
Let $m\in\N$. If $\eta \sim\mathcal{N}_{\C^m}$, and $x\in\C^m$ then
$$
\E_{\eta \sim\mathcal{N}_{\C^m}}(|\langle \eta,x\rangle|^2)=\|x\|^2.
$$
The proof of this claim is a standard exercise and is left to the reader.

Now we are ready to prove the lemma. We chose a representative of $v$ such that $\|v\|=1$ and a representative of the left eigenvector $u$ such that $\langle u,v\rangle=1$. Note that this implies by Proposition~\ref{pro:cn-eigenv} and Lemma~\ref{pro:cn-eigenv}:
\[
 \mu_\lambda(A,\lambda,v)=\frac{\|u\|\,\|v\|}{|\langle u,v\rangle}=\|u\|,\quad \dot \lambda=\langle \dot Av,u\rangle,\quad \dot v=-A_{\lambda,v}^{-1}P_{v^\perp}\dot Av.
\]
For the first statement, we have
\begin{align*}
\E_{\dot A\in T_A\SS} (|\dot\lambda|^2)=\E_{\dot A\in T_A\SS} (|\langle \dot Av,u\rangle|^2)
&=\E_{\dot A\in T_A\SS} (|\langle \dot A,uv^*\rangle_F|^2)\\
&=\E_{\dot B\in\mcN_{A^\perp}}
\E_{t\in\mathcal{N}(0,\frac12)}(|\langle t\sqrt{-1}A+\dot B,uv^*\rangle_F|^2).
\end{align*}
Since the mixed term of the expansion of $|\langle t\sqrt{-1}A+\dot B,uv^*\rangle_F|^2$ is linear in $t$, its expected value is zero. Hence,
\begin{align*}
\E_{\dot A\in T_A\SS} (|\langle \dot Av,u\rangle|^2)
&=\E_{\dot B\in\mcN_{A^\perp}}
\E_{t\in\mathcal{N}(0,\frac12)}( t^2|\langle A,uv^*\rangle_F|^2+|\langle\dot B,uv^*\rangle_F|^2)\\
&= \frac{|\lambda|^2}{2}+\E_{\dot B\in\mcN_{A^\perp}}(|\langle\dot B,\pi_{A^\perp}(uv^*)\rangle_F|^2.
\end{align*}
where we have denoted by $\pi_{A^\perp}(uv^*)$ the orthogonal projection of $uv^*$ onto $A^\perp$. With the identification of $A^\perp$ and $\C^{n^2-1}$ as Hermitian spaces, from Claim III (with $m=n^2-1$) we conclude that
\begin{align*}
 \E_{\dot B\in\mcN_{A^\perp}}(|\langle\dot B,\pi_{A^\perp}(uv^*)\rangle_F|^2 = \|\pi_{A^\perp}(uv^*)\|^2_F=\|uv^*\|_F^2-|\lambda|^2=\|u\|^2-|\lambda|^2.
\end{align*}
We have then proved:
\[
 \E_{\dot A\in T_A\SS} (|\langle \dot Av,u\rangle|^2)=\mu_\lambda(A,\lambda,v)^2-\frac{|\lambda|^2}{2},
\]
and from claim I we conclude:
\[
 \E_{\dot A\in\mathcal{S}_A} (|\langle \dot Av,u\rangle|^2)=\frac{1}{\left(n^2-\frac12\right)}\left(\mu_\lambda(A,\lambda,v)^2-\frac{|\lambda|^2}{2}\right),
\]
as claimed in the lemma. The second statement in the lemma is proved in a very similar fashion. This time we have
\begin{align*}
\E_{\dot A\in T_A\SS} (|\dot v|^2)=\;&\E_{\dot A\in T_A\SS} (\|A_{\lambda,v}^{-1}P_{v^\perp}\dot Av\|^2)=\\
=\;&\E_{\dot B\in\mcN_{A^\perp}}
\E_{t\in\mathcal{N}(0,\frac12)}\left(\left\|A_{\lambda,v}^{-1}P_{v^\perp}\left(t\sqrt{-1}A+\dot B\right)v\right\|^2\right)=\\=\;&\E_{\dot B\in\mcN_{A^\perp}}
\left(\left\|A_{\lambda,v}^{-1}P_{v^\perp}\dot Bv\right\|^2\right)=\E_{\dot w\in\mcN_{v^\perp}}(\|A_{\lambda,v}^{-1}\dot w\|^2)=\|A_{\lambda,v}^{-1}\|^2_F,
\end{align*}
the previous to last equality from claim II and the last coming from the fact that for any matrix $B\in\Cnn$,
$\E_{x\sim\mathcal{N}_{\C^n}}\|Bx\|^2=\|B\|_F^2$ (note the use of
Frobenius instead of operator norm in the last equality: that is a
crucial point). The last statement of the lemma then follows from claim I.

\eproof

\section{Proof of Theorem~\ref{th:mainrandom}}\label{sec:proofmainrandom}

\subsection{Proof of $(1)$ and $(2)$ in Theorem~\ref{th:mainrandom}}
Note that $(2)$ is trivial. We thus prove $(1)$. The procedure we
suggest to choose $\omega\in\Omega_n$ at random is the following (note
that each step requires $\Oh(n^3)$ arithmetic operations or random
choices):
\begin{enumerate}
\item 
Choose $B\sim\mathcal{N}_{\C^{(n-1)\times (n-1)}}$ and let $U$ be
the $Q$ factor in the QR decomposition of $B$, then multiply $Q$ by
the diagonal matrix with entries $r_{ii}/|r_{ii}|$ where the $r_{ii}$
are the diagonal elements of the $R$ factor. This produces a unitary
matrix $U$ uniformly distributed in $\CU_{n-1}$, see for
example~\cite{Mezzadri2007}.
\item 
Choose $\lambda\sim \mathcal{N}_\C$ and
$M\sim\mathcal{N}_{\C^{(n-1)\times n}}$. Let $H\in\CU_n$ be any
unitary matrix such that its last column is in $\ker(M)$ (it is trivial to produce such an $H$ by computing the QR decomposition of the matrix whose columns are $\ker(M)$ and the columns of $M$). Compute
$Q$ as the product of the first $n\times(n-1)$ submatrix of $H$ times
$U$. This produces an element with the uniform distribution in the set
of $Q\in\mathcal{S}_{n-1}(\C^n)$ such that $(M,Q)\in\A_n$.
\item 
If $2\Re(\bar{\lambda}\mathrm{tr}(MQ))> 1-|\lambda|^2(n-1)$ then
discard $\lambda,M,Q$ and repeat $(1)$ and $(2)$.
\item Choose $w\sim\mathcal{N}_{\C^{n-1}}$.
\end{enumerate}
The only subtle point is that steps $(1)$ and $(2)$ might have to be
repeated an arbitrary number of times. The expected number of times
that steps $(1)$ and $(2)$ will be repeated is related to $C_n$
defined in~\eqref{eq:Cn} by
\begin{align*}
\sum_{k=1}^\infty \Prob(\text{step $k$ is reached})=&
\sum_{k=1}^\infty  \Prob\left(2\Re(\bar{\lambda}\mathrm{tr}(MQ))
> 1-|\lambda|^2(n-1)\right)^{k-1}\\=&
\sum_{k=1}^\infty \left(1-C_n^{-1}\right)^{k-1}
=\frac{1}{\left(1-(1-C_n^{-1})\right)}=C_n.
\end{align*}

\subsection{Proof of $(3)$ in Theorem~\ref{th:mainrandom}}

We are now prepared for proving~\eqref{eq:wish2}. Let $\uno$ be the
characteristic function of the set
\[
  \{(\lambda,B):2\Re(\bar{\lambda}\mathrm{tr}(B))
  \leq 1-|\lambda|^2(n-1)\}\subseteq\C\times\C^{(n-1)\times(n-1)}.
\]
From the definition
and Fubini's theorem, for any measurable nonnegative function
$\phi$ defined on $\V$, the expected value 
$\mathbb{E}_{\omega\sim\Omega_n}\left(\phi(\psi_n(\omega))\right)$ 
equals:
\[
  C_n\,\E_{M}\; \E_{Q:(M,Q)\in\A_n}
  \left(\E_{\lambda,w}\left(\phi\left(\begin{pmatrix}
 \lambda&w^*\\0&MQ+\lambda I_{n-1}
\end{pmatrix},\lambda,e_1
\right)\;  \uno(\lambda,MQ)\right)\right)=
\]
\[
  C_n\,\E_{M}\; \E_{Q:(M,Q)\in\A_n}\left(\alpha(MQ)\right),
\]
where $\lambda\sim\mathcal{N}_\C$, $M\sim
\mathcal{N}_{\C^{(n-1)\times n}}$,
$w\in\mathcal{N}_{\C^{n-1}}$ and
$\alpha\colon \C^{(n-1)\times(n-1)}\rightarrow[0,\infty]$ is defined by
\[
  \alpha(B)=\E_{\lambda,w}\left(\phi\left(\begin{pmatrix}
  \lambda&w^*\\0&B+\lambda I_{n-1}
  \end{pmatrix},\lambda,e_1
\right)\;  \uno(\lambda,B)\right).
\]
We are then under the hypotheses of Corollary~\ref{cor:trick2}. 
Using this result we obtain 
\begin{align*}
  \E_{w\sim\Omega_n}\left(\phi(\psi_n(w))\right)=&\frac{C_n}{\Gamma{(n)}}
  \E_{B\sim\mathcal{N}_{\C^{(n-1)\times (n-1)}}}(\alpha(B)|\det(B)|^2).
\end{align*}
With the change of variables $B+\lambda I_{n-1}=D$, which implies
\[
 \|B\|_F^2=\|D\|_F^2+(n-1)|\lambda|^2-2\Re(\bar{\lambda}\mathrm{tr}(D)), 
\]
this last expression equals
\[
  \frac{C_n}{\Gamma{(n)}}
  \E_{\lambda,w,D}\left(\phi\left(\begin{pmatrix}
  \lambda&w^*\\0&D
\end{pmatrix},\lambda,e_1
\right)\;  
 |\det(D-\lambda I_{n-1})|^2e^{-|\lambda|^2(n-1)
 +2\Re(\bar{\lambda}\mathrm{tr}(D))}\uno(\lambda,D-\lambda I_{n-1})\right),
\]
where $D\sim\mathcal{N}_{\C^{(n-1)\times (n-1)}}$.
Now, note that
\[
  \uno(\lambda,D-\lambda I_{n-1})\neq0\Longleftrightarrow 
 e^{-|\lambda|^2(n-1)+2\Re(\bar{\lambda}\mathrm{tr}(D))}\leq e.
\]
We have thus proved
\begin{align*}
  \E_{w\sim\Omega_n}\left(\phi(\psi_n(w))\right)
  \leq& \frac{e\,C_n}{\Gamma{(n)}}
  \E_{\lambda,w,D}\left(\phi\left(\begin{pmatrix}
\lambda&w^*\\0&D
\end{pmatrix},\lambda,e_1
   \right)\;  |\det(D-\lambda I_{n-1})|^2\right)\\
  \underset{\eqref{eq:citame1}}{=}& 
  e\,nC_n\E_{A\sim\mcN_{\C^{n\times n}}}
  \left(\frac{1}{n}\sum_{\lambda,v:Av=\lambda v}
  \phi(A,\lambda,v)\right).
\end{align*}
This proves claim $(3)$ in Theorem~\ref{th:mainrandom}.
\eproof

We prove the following (non-sharp) bound for the value of $C_n$.

\begin{lemma}\label{lem:SnCn}
With the notations above,
\[
  C_n\leq 4n.
\]
\end{lemma}

\proof
Note that if $0<|\lambda|\leq (n-1)^{-1/2}$, then for any nonzero
$M\in\C^{(n-1)\times n}$ we have
\[
  \Prob_{Q}\left(2\Re(\bar{\lambda}\mathrm{tr}(MQ))
  \leq 1-|\lambda|^2(n-1)\right)\geq 
  \Prob_{Q}\left(2\Re(\bar{\lambda}\mathrm{tr}(MQ))\leq 0\right)=\frac{1}{2},
\]
the last equality coming from the linearity of the trace. We thus have
\[
  \Prob_{\lambda,M,Q}\left(2\Re(\bar{\lambda}\mathrm{tr}(MQ))
  \leq 1-|\lambda|^2(n-1)\right)\geq\frac{1}{\pi}\int_{|\lambda|
  <(n-1)^{-1/2}}\frac{e^{-|\lambda|^2}}{2}\,d\lambda
  =\frac{1-e^{-\frac{1}{n-1}}}{2}.
\]
We thus have
\[
  C_n\leq \frac{2}{1-e^{-\frac{1}{n-1}}}\leq 4(n-1)\leq 4n,
\]
as claimed.
\eproof

\section{Proof of Theorem \ref{th:relativeerror}
\label{sec:mainrelativeepsilon}}

\noindent
Consider the following algorithm.
\bigskip\bigskip

\algoritmo
\begin{algorithm}\label{alg:AllRE}
{\sf Relative\_Error}\\
\inputalg{$\e\in(0,1/2)$, $A\in \Cnn$, 
$(\zeta,w)\in\C\times\C^n$}\\
\specalg{$(A,\zeta,w)\in\wW$ so $(\zeta,w)$ is an approximate 
eigenpair of $A$ with associated eigenpair 
$(\lambda,v)$, $\|w\|=1$, $\|A\|_F=1$}\\
\bodyalg{
$k:=0$\\[2pt]
$(\zeta',w'):=(\zeta,w)$\\[2pt]
repeat\\[2pt]
\espacio $(\zeta',w'):=N_A(\zeta',w')$ (one Newton iteration)\\[2pt]
\espacio $k:=k+1$\\[2pt]
until $k\geq\log_2\log_2(\frac{4}{\e|\zeta'|})$\\[2pt]
return $(\zeta',w')$\\[2pt]
}
\Output{$(\zeta',w')\in \C\times \C^{n}$}\\
\postcond{The algorithm halts if $\lambda\neq0$. 
In this case, $(A,\zeta',w')\in\wW$ so $(\zeta',w')$ is an approximate 
eigenpair of $A$ with associated eigenpair 
$(\lambda,v)$, and $d_{\SS}(w',v)\leq \e$, and moreover
$|\zeta'-\lambda|\leq\e|\lambda|$.}
\end{algorithm}
\falgoritmo
\bigskip

By hypothesis,
\begin{equation}\label{eq:apeig1}
  \dist_A((\zeta,w),(\lambda,v))\leq \frac{c_*}{\mum(A)}<1.
\end{equation}
Hence, $|\zeta-\lambda|\leq \dist_A((\zeta,w),(\lambda,v))<1$ 
and the same bound holds with $\zeta$ replaced by $\zeta'$ 
at all the iterations of the algorithm (by 
Definition~\ref{def:app_eigen}). Using that 
$|\lambda|\leq \|A\|_F=1$, we deduce that $|\zeta'|\leq2$. 
Hence, at the end of the {\tt repeat} loop the value $k$ satisfies 
\begin{equation}\label{eq:apeig2}
    2^{2^k}\geq \frac{4}{\e|\zeta'|} 
    \geq \frac{2}{\e}.  
\end{equation}
This inequality implies that after $k$ iterations of the loop we 
have, from the definition of approximate eigenpair 
and the bound~\eqref{eq:apeig1}, that 
\begin{equation}\label{eq:apeig3}
   \dist_A((\zeta',w'),(\lambda,v))\leq
   \frac{c_*}{\mum(A)}\,\left(\frac12\right)^{2^k-1}\leq
   \left(\frac12\right)^{2^k-1}=
   \frac2{2^{2^k}}\leq\e.  
\end{equation}
In particular, $d_{\SS}(w',v)\leq \e$ as we wanted. 
On the other hand, the first inequality 
in~\eqref{eq:apeig2} implies 
$$
   2^{2^k-1}|\zeta'|-1 \geq \frac{2}{\e}-1 \geq \frac{1}{\e}
$$
the last since $\e<1$. We now use this inequality together 
with~\eqref{eq:apeig3} to obtain
\[
  \frac{|\zeta'-\lambda|}{|\lambda|}\leq  \frac{1}{|\lambda|2^{2^k-1}}
  \leq  \frac{1}{\left(|\zeta'|-\frac{1}{2^{2^k-1}}\right)2^{2^k-1}}
  = \frac{1}{2^{2^k-1}|\zeta'|-1}\leq \e,
\]
i.e., $|\zeta'-\lambda|\leq\e|\lambda|$. 

It remains to show that {\sf Relative\_Error} halts provided 
$\lambda\neq0$ and to estimate its average running time 
when $A$ is drawn from $\mathcal{N}_{\C^{n\times n}}$. 
For this, we note that as soon as
\[
  k\geq \log_2\log_2\left(\frac{8}{|\lambda|\e}\right)
\]
we shall have (using~\eqref{eq:apeig3}) 
\[
   |\zeta'|\geq|\lambda|-\left(\frac12\right)^{2^k-1}\geq 
   |\lambda|-\frac{|\lambda|\e}{4}= |\lambda|(1-\e/4).
\]
Therefore, we will also have
\[
    \log_2\log_2\left(\frac{4}{\e|\zeta'|}\right)\leq 
    \log_2\log_2\left(\frac{4}{\e|\lambda|(1-\e/4)}\right)
    \leq \log_2\log_2\left(\frac{8}{|\lambda|\e}\right)\leq k.
\]
Hence, the stopping condition will hold after at most
\[
  \log_2\log_2\left(\frac{8}{|\lambda|\e}\right)
  \leq \log_2\log_2\left(\frac{8\|A^{-1}\|}{\e}\right)
\]
iterations (we have used that $|\lambda|^{-1}\leq \|A^{-1}\|$ 
for $\lambda$ is an eigenvalue of $A$). 

We finally estimate the average cost of {\sf Relative\_Error}. 
Since each iteration of the {\tt repeat} loop requires $\Oh(n^3)$  
operations, this cost is at most $\Oh(n^3)$ times 
\[
     \E_{A\sim\mathcal{N}_{\C^{n\times n}}}
   \left(\log_2\log_2\left(\frac{8\|A^{-1}\|}{\e}\right)\right)\leq 
   \log_2\log_2\left(\frac{8}{\e}\E_{\mathcal{N}_{\C^{n\times n}}}
   (\|A^{-1}\|)\right),
\]
where we have used Jensen's inequality. Bounds for the expected value of
$\|A^{-1}\|$ when $A\sim\mcN_{\Cnn}$ are known, 
see for example~\cite[Prop. 4.22]{Condition} which, together 
with the general inequality $\E(Z)^2\leq \E(Z^2)$, implies
\[
   \E_{A\sim\mathcal{N}_{\C^{n\times n}}}(\|A^{-1}\|)\leq\sqrt{\frac{e(n+1)}{2}}
   \leq\sqrt{e\,n}.
\]
The statement follows.
\eproof

\noindent
{\bf Acknowledgments.\quad}
We want to thank Jim Demmel for helpful comments and 
Froil\'an Dopico for suggesting us to work out the case of relative error. P. B\"urgisser, F. Cucker and M. Shub want to acknowledge the Simons institute where they spent some time while this work was being carried out.
\bigskip

{\small
%\bibliography{../../book/book}

}

\end{document}